\documentclass[10pt]{article}
\title{On representation-finite selfinjective algebras, coverings, multiplicative bases ...}
\author{Klaus Bongartz\thanks{klausbongartz@online.de}\\ Universität Wuppertal}
\usepackage{amsthm}
\usepackage[utf8]{inputenc}

\newtheorem{theorem}{Theorem}
\newtheorem{lemma}{Lemma}
\newtheorem{corollary}{Corollary}
\newtheorem{proposition}{Proposition}

\newcommand\me{k(\Z\Delta)}

\newcommand\Z{\mathbf{Z}}
\newcommand\N{\mathbf{N}}

\begin{document}
\maketitle
\begin{abstract}{ We give a simplified complete proof for 
the classification of the selfinjective representation-finite algebras. We explain the relations between the two different approaches and also to further developments.
Some historical remarks are made.}
\end{abstract}

\section{Introduction}

This article concerns results published  between 1978 and 1986 and some supplements added after 2009. It consists in two parts of a quite different character and length. 
The first part deals with
the classification of selfinjective representation-finite  algebras   and  the second part contains comments on this classification and 
 parallel or later
work on mild algebras i.e. on algebras where all proper quotients are representation-finite.  

There are two different approaches to the classification: via coverings of Auslander-Reiten quivers or via coverings of quivers with relations resp. of trivial extensions. Following most
 of the time the first approach  we give a proof 
for the classification which is complete apart from basics on covering theory and some lists obtained by computer. We simplify many  arguments, 
correct some errors and include also the results on trivial extensions of the second approach. Our thorough treatment should help  to understand this interesting subject 
 not yet treated in a textbook.

The second part starts in section 12. It contains no proofs  but many mathematical and historical comments. We compare the two approaches to the classification thereby mentioning 
an essential gap in the second approach. Then we discuss multiplicative bases, coverings of ray categories, the finiteness criterion, the Brauer-Thrall conjectures and the intimate relations 
between these themes.  We emphasize facts and points of view that are often formulated sloppily or omitted or misunderstood in the existing literature.

Now we describe the content and the line of thought of the first part in more detail. 
 Section 2 recalls  translation-quivers, their
 universal coverings and fundamental groups and the fact that any  connected stable translation-quiver  is the quotient $\Z T/\Pi$ 
for some tree $T$ and some admissible group  $\Pi$ of automorphisms of $\Z T$. 

 In section 3 
  we prove  the existence of a covering functor relating the mesh category $k(\tilde{C})$ 
of the universal cover $\tilde{C}$ of an appropriate  
connected component $C$ of the Auslander-Reiten quiver of a locally bounded category $A$ to the category $ind\,C$ of indecomposables lying in $C$. This key result has
 two nice consequences 
for any  locally representation-finite category $A$ with Auslander-Reiten quiver $\Gamma_{A}$. First, any connected component of the stable part of $\Gamma_{A}$ is 
a quotient $\Z T/\Pi$ for some Dynkin-tree $T$ and some cyclic group $\Pi$.
 Second, any quotient $\tilde{\Gamma}_{A}/H$ by some admissible group $H$ of automorphisms of $\tilde{\Gamma}_{A}$ is the Auslander-Reiten quiver
 of its  standard-category  which is defined as the full subcategory of  $k(\tilde{\Gamma}_{A}/H)$ consisting of  the projective points. For the special case $H=\{e\}$ one obtains the universal
 cover $\tilde{A}$ of $A$ which has always a 'standard' presentation by zero- and commutativity relations. 

Some natural problems arise:

 1. Can we describe the standard-category to $k(\Gamma_{A})$ directly in terms of $A$?  

2. Does the Auslander-Reiten quiver determine the algebra?

3. Is it possible to classify all locally representation-finite categories with simply connected Auslander-Reiten quiver?

For the special case of selfinjective categories we give in sections 4 - 11 the solutions of these problems with complete proofs. 
For arbitrary locally  representation-finite categories 
the problems are also solved to a large extent as discussed
in the second part of this article.

First we answer the last question and  we study in section 4 a selfinjective locally representation-finite category $A\neq k$ with simply connected Auslander-Reiten quiver $\Gamma_{A}$.
 Its stable part is thus
 isomorphic to $\Z T$ for some Dynkin-diagram $T$, the tree class of $A$.
We define for  any subset $C$ of $\Z T$ a  
translation quiver
 $\Z T_{C}$ by adding for each $c \in C$ a new point $c^{\ast}$ which is projective and 
injective and arrows 
$c \rightarrow c^{\ast}$ , $c^{\ast} \rightarrow \tau^{-1}c$. 
 Here $C$ is called a configuration if $\Z T_{C}$ is the Auslander-Reiten
quiver of a locally representation-finite selfinjective category $A(C)$ which is then uniquely determined by $C$. Any 
configuration $C$ is stable under $\tau^{L(T)}$ where $L(T)-1$ is the Coxeter number of $T$. 
Thus the Nakayama-automorphism $\nu$ of $k(\Z T_{C})$ and $A(C)$ is induced by $\tau^{L(T)}$. 

For a better description of the configurations we fix an arbitrary  section  $S$   in $\Z T$. Given a configuration $C$ 
 we take the vector $d=d(S,C)$ of total dimensions of the indecomposable
 $A(C)$-modules lying on $S$. It follows by elementary Auslander-Reiten theory that the full subcategory $B=B(S,C)$ of $A(C)$ consisting of the supports of the modules in $S$ 
 is a simply connected algebra having $S$ as a
 section in its finite Auslander-
Reiten quiver. The pattern of projective points inside $\Z T$ is obtained from $d$ as usual by knitting backwards and the configuration 
by knitting and knotting forwards. 
Here knotting means inserting projective-injective points. 
This very simple and effective procedure induces  a bijection between the configurations and the patterns of projective points. Furthermore  one gets $A(C)$ from $B(S,C)$ easily by 
periodic repetitions and  by adding certain arrows and relations applying a simple recipe. 

Thus it remains to find the possible patterns for a suitable section or equivalently the algebras $B(S,C)$.
For the tree class $A_{n}$ we choose an equioriented section and we include in section 5 a detailed discussion of the emerging pedigree-algebras. Then one gets in section 6 the algebras
 for the tree class $D_{n}$ as
one-point extensions from those of $A_{n-1}$. The classes $E_{6},E_{7},E_{8}$ can be treated with the same method which is easy to implement on a computer.
 Our solutions of problem 3 is complete, but for the description of all  Auslander-Reiten quivers $\Z T_{C}/\Pi$   occurring for a
 representation-finite selfinjective algebra  we give in section 7 a table of the possible admissible groups $\Pi$.

To study  the first two problems  we need some general concepts most of which were introduced to  study selfinjective categories. 
First we consider multiplicative bases
in distributive categories. A locally bounded category $A$ is distributive resp. regular if  all morphism spaces $A(x,y)$ are uniserial as modules over $A(x,x)$ or over $A(y,y)$ 
resp. over both.  Locally representation-finite categories are always distributive. A  multiplicative basis in a distributive category $A$ is a set $M$ of morphisms 
containing all identities such that $M\cap \mathcal{R}A(x,y)$  is a basis of each radical $\mathcal{R}A(x,y)$ and such that the composition of two morphisms in $M$ is again in $M$ or $0$. 
The basis is called cancellative if 
- in addition - $fgh=fh\neq 0$
with $f,g,h \in M$ implies that $g$ is an identity. In general, a distributive category $A$ does not have a multiplicative basis but as substitutes there is the 
  stem-category $\hat{A}$ 
and the ray-category $\vec{A}$ of $A$. Both categories are defined in terms of the subbimodules of the various $A(x,y)$ and their products.
 The linearization  $k\hat{A}$  has a  multiplicative bases and $k\vec{A}$ even a cancellative one. Both categories determine each other and they coincide for a regular category $A$.
Fortunately,  any regular representation-finite selfinjective algebra $A$ is isomorphic to $k\hat{A}$.

In section 9 we look at Galois coverings. For any locally representation-finite category $A$ the fundamental group $\Pi$ 
acts on the universal cover $\tilde{A}$ by automorphisms and  on its quiver $\tilde{Q}$ by  permutations respecting  the standard 
presentation. The  resulting 'combinatorial quotient'  $\tilde{A}//\Pi$ has  a cancellative basis and the same Auslander-Reiten quiver as $A$ and two representation-finite 
algebras with the same Auslander-Reiten quiver have isomorphic stem-algebras. Therefore problems 1 and 2 are solved except for some non-regular algebras having tree class $D_{3m}$ and a configuration 
$C$ stable under $\tau^{2m-1}$. These algebras are studied in section 11 about large fundamental groups.

Section 10 deals with the trivial extension $T(F)$ of an algebra $F$. There is always the repetitive category $\hat{T}(F)$ 
as a Galois cover of $T(F)$ with group $\Z$.  For representation-finite  $T(F)$ one has $\hat{T}(F) \simeq A(C)$ with some configuration $C$. Given any fundamental algebra $F$
 of $A(C)$
- i.e. a convex connected subcategory hitting each $\nu$-orbit exactly once -
one gets $T(F)\simeq A(C)// \langle \nu \rangle$. Therefore the isomorphism classes of representation-finite trivial extensions are in bijection to those of locally representation-finite
 selfinjective simply connected categories. For the tree class $A_{n}$ the trivial extensions are the well-known Brauer-quiver algebras and for $D_{n}$ 
 one gets some modifications thereof.

I want to thank Dieter Vossieck and Henning Krause for their interest on part one and Andrew Hubery for proposing 'knitting and knotting' as a translation of 'Stricken und Sticken'.

\section{Translation quivers and their coverings}
The material resumed here stems from \cite[section 1]{Rie1}  and \cite[sections 1,2] {BG} where one finds also many illustrations.
 A translation quiver  $(\Gamma,\tau)$ is a pair 
consisting of a  quiver $\Gamma$ without loops and double-arrows and a bijection $\tau:X \rightarrow Y$ between two subsets of 
$\Gamma_{0}$. Furthermore, for any  $x \in \Gamma_{0}$ the set  $x^{+}= \{y \mid \exists x \rightarrow y\}$  and $x^{-}=\{z \mid \exists z \rightarrow x\}$  are finite with
  $(\tau x)^{+}=x^{-}$  for  each $x \in X$. In that case an arrow $\alpha:y \rightarrow x$ corresponds to an arrow  $\sigma \alpha:\tau x \rightarrow y$. We call the vertices 
in $\Gamma_{0}\setminus X$ projective, those in $\Gamma_{0} \setminus Y$ injective and $\tau$ the translation. Usually we only write $\Gamma$ for $(\Gamma,\tau)$.
 The most important examples of translation quivers are Auslander-Reiten quivers and some variants thereof.

A translation subquiver $\Gamma'$ is a translation quiver $(\Gamma',\tau' )$ such that $\Gamma'$ is a full subquiver of $\Gamma$. Furthermore a  point $x'$ of $\Gamma'$ is projective in $\Gamma'$
iff
 it is
 so in $\Gamma$ or if $\tau x'$ does not belong to $\Gamma'$ and for the non-projective points we have $\tau x'=\tau' x'$. Analogous statements  hold for $\tau^{-1}$. So a translation subquiver 
is uniquely determined 
by its point set $\Gamma'_{0}$ as a subset of $\Gamma_{0}$. We simply say that $\Gamma'$ is embedded in $\Gamma$.

A point $x$ is  stable if all powers of $\tau$ and $\tau^{-1}$ are defined on $x$. The stable part $\Gamma_{s}$ of $\Gamma$ is the translation subquiver consisting of the stable points.
For any  ( unoriented ) tree $T$  one has the stable translation quiver $\Z T$ with point set $\Z \times T_{0}$. An edge $x - y$ induces arrows
 $(i,x) \rightarrow (i,y) \rightarrow (i+1,x)$ for each integer $i$ and $\tau$ maps $(i,x)$ to $(i-1,x)$ for all $i$ and $x$. Quite often we consider a section $S$ in $\Z T$ i.e. a connected subquiver
 such that its points represent  the $\tau$-orbits and its arrows the $\sigma$-orbits.

A morphism of translation quivers is a map $f:\Gamma' \rightarrow \Gamma$  between the quivers with $f \tau' =\tau f$. Such a morphism is a covering if  $x'$ is projective resp.
 injective
iff $fx'$ is so and if $f$ induces bijections between  $x^{+}$ and  $fx^{+}$ resp. $x^{-}$ and $fx^{-}$. A group $G$ of automorphisms of $(\Gamma, \tau )$ is admissible if no orbit of $G$ on $\Gamma_{0}$ intersects
$\{x\} \cup x^{+}$ or $\{x\} \cup x^{-}$ in two points. Then one has  the quotient $\Gamma/G$ and a projection $\pi:\Gamma \rightarrow \Gamma/G$ which is a covering of translation quivers.
The points of $\Gamma /G$ are just the orbits $Gx$ of the points $x$ of $\Gamma$. There is one and only one way to define the arrows and the translation of $\Gamma /G$ such that 
$\pi(x)= Gx$ is a covering of translation quivers.

 For any connected translation quiver  $(\Gamma,\tau)$ the universal cover $(\tilde{\Gamma},\tilde{\tau})$  is constructed as follows.  One defines a new quiver
 $\hat{\Gamma}$ by adding a new arrow $\gamma_{x}:\tau x \rightarrow x$ 
for each $x$ in $X$. A walk in
 $\hat{\Gamma}$ is a  composition of old and new arrows and their formal inverses such that the occurring starts and ends fit together well and the composition of walks is defined in the
 obvious way. The homotopy is the smallest equivalence relation stable under left or right multiplication with the same walk, under 'inversion' and such
 that $\alpha\alpha^{-1}$ is equivalent to an identity for each
arrow $\alpha$ - old or new - or each formal inverse.  Furthermore $\gamma_{x}$ is equivalent to $\alpha \sigma(\alpha)$ for each arrow $\alpha$ ending in a non-projective vertex $x$.
Now  we fix an arbitrary  base point $z$. The following definitions are  independent of the choice of $z$ up to isomorphism.
The fundamental group $\Pi$ is defined as the set of all homotopy classes of walks starting and ending in $z$ endowed with the multiplication induced by the composition of walks.
The universal cover $\tilde{\Gamma}$ has the homotopy classes with start in $z$ as its points. The arrows, the translation and the covering $\pi:\tilde{\Gamma} \rightarrow \Gamma$ are all 
defined in the obvious way. Furthermore $\Pi$ acts as  an admissible group of automorphisms on $\tilde{\Gamma}$ and $\Gamma$ is the quotient.
A translation quiver is simply connected if it is connected and if its universal covering is the identity.
  
The fact that the homotopy relation is homogeneous provided one gives the new
 arrows $\gamma_{x}$ the degree 2 and the inverses of the arrows the corresponding negative degree implies that  all walks between two fixed points have the same 'length'.
As  an important though trivial consequence  there is  a morphism $\kappa$ of translation quivers from $\tilde{\Gamma}$ to ${\bf Z}A_{2}$   which makes it possible to argue 
 by induction. In particular  $\tilde{\Gamma}$  is interval-finite i.e. there are only finitely many paths between two points. In particular there is no oriented cycle.
 
A  basic  important result is the structure theorem for stable components from \cite[Struktursatz]{Rie1}:

\begin{theorem}
 A translation quiver  is   stable and  simply connected iff it  is isomorphic to $\Z T$ for some tree $T$ which is unique up to isomorphism.
\end{theorem}

\section{ Covering functors }
\subsection{Basic definitions and conventions}

We will always work over a fixed algebraically closed field $k$. A locally finite dimensional category $A$ is a $k$-category such that distinct objects are not isomorphic, all homomorphism spaces $A(x,y)$
 have finite dimension and the endomorphism algebras $A(x,x)$ are local. The radical $\mathcal{R}A$ is the two-sided ideal such that $\mathcal{R}A(x,y)$ is the set of 
non-invertible morphisms from $x$ to $y$. The quiver $Q_{A}$ of $A$ has the objects as points and $dim  (\mathcal{R}A(x,y)/\mathcal{R}^{2}A(x,y)$ arrows from $x$ to $y$. If there is at most one
 arrow beteen two points $A$ is called square-free. It is locally bounded if for any $x$ in $A$ there are only finitely many $y$ with $A(x,y)\neq 0 $ or $A(y,x) \neq 0$.
Any locally bounded category is isomorphic to $kQ/I$ where $Q$ is the locally finite quiver of $A$ and $I$ is a  twosided ideal in the path category. Basic finite dimensional algebras correspond to
locally bounded categories with finite quiver.

For a $k$-category $A$ a module $M$ is  a contravariant $k$-linear functor $M:A \rightarrow Mod \,k$ into the category of  vector spaces. We denote this category by $Mod \,A$, but
 we are mainly interested in its full subcategory $mod\,A$ of functors where the total dimension $dim\,M := \sum_{x} \,dim_{k}\,M(x)$ is finite. Such a functor is a finite direct sum of indecomposable 
functors with local endomorphism algebras whence the decomposition is essentially unique.. We write $Hom_{A}(M,N)$ for the set of morphisms in $mod\,A$ and $\mathcal{R}$ for the radical.  Yonedas lemma implies that in a 
locally bounded category $A$ an indecomposable 
is projective resp. injective iff it is isomorphic to $A(\,\,,x)$ resp. to $DA(x,\,\,)$ for some $x$ in $A$. As is well-known $mod\,A$ has almost split sequences.
A $k$-category is called locally representation-finite if for any point $x$ there are up to isomorphism only finitely many indecomposables $U$ in $mod\,A$ with $U(x)\neq 0$.
An equivalent condition is that the full subcategory $ind\,A$ of $mod\,A$ supported by representatives of the isomorphism classes of indecomposables is locally bounded.

An important example of a $k$-category is the mesh-category $k(\Gamma)$ to a translation quiver $(\Gamma,\tau )$.
 For a non-projective vertex $x$ in $\Gamma$ the mesh ending in $x$ is the full subquiver of $\Gamma$ supported by $x, \tau x$ and the predecessors of $x$. So let 
$\alpha_{i}:x_{i } \rightarrow x$, $1 \leq i \leq r$, be all arrows ending in $x$ and $\sigma(\alpha_{i}): \tau x \rightarrow x_{i}$ be all arrows starting in $\tau x$. Then the
 mesh-sum $\mu_{x}$ is defined as 
$\mu_{x}=\sum_{i} \alpha_{i}\sigma (\alpha_{i})$. The mesh-ideal $I_{\Gamma}$ is the ideal of $k\Gamma$ generated by all mesh-sums $\mu_{x}$ for non-projective points $x$ and $k(\Gamma)$ 
is the quotient $k\Gamma/I_{\Gamma}$.  If $k(\Gamma)$ is locally bounded we call $\Gamma$  locally bounded. 
A trivial but very useful fact is that a translation subquiver $\Delta$ of  a locally bounded translation quiver $\Gamma$ is again locally bounded. Namely $k(\Delta)$  is isomorphic
 to the quotient of $k(\Gamma)$ 
by the ideal generated by the paths leading through a point outside $\Delta_{0}$. 

A $k-$linear  functor $F:M \rightarrow N$  is a covering functor  if $F$ hits all objects in $N$ and if it induces for all $m \in M$ and $n \in N$ isomorphisms 
$$F_{(m,n)}:\bigoplus_{Fm'=n} M(m,m') \simeq N(Fm,n)  \, , \,F_{(n,m)}:\bigoplus_{Fm'=n} M(m',m) \simeq N(n,Fm).$$ The existence of a covering functor $F:M \rightarrow N$ implies that $M$ is locally bounded 
iff $N$ is so.

Let $F: M \rightarrow N$ be a covering functor and let $N'$ be a full subcategory of $N$ with 'preimage' $M'=f^{-1}N'$. Then there are induced covering functors $F':M' \rightarrow N'$ and 
$F'':M/M' \rightarrow N/N'$ where $M/M'$ resp. $N/N'$ are obtained by dividing out all morphism factoring through direct sums of objects in $M'$ resp. $N'$. 

Any covering $f:\Gamma' \rightarrow \Gamma$ of translation quivers 
induces a covering functor $k(f):k(\Gamma') \rightarrow k(\Gamma)$. A less obvious covering functor is constructed  in the next section.

\subsection{Riedtmanns  existence theorem} 

Let $A$ be a locally bounded category and $C$ a connected component of its Auslander-Reiten quiver with universal cover $\pi:\tilde{C} \rightarrow C$. For each point $y$ in $C$ 
we choose an indecomposable $U(y)$ and we denote by $ind\,C$ the full subcategory of $mod\,A$ with the $U(y)$'s as objects.
 Suppose that $ind\,C$ is square-free and that for all points $y,y'$ in $C$ there is a natural number $n(y,y')$ 
with  $\mathcal{R}^{n(y,y')}(U(y),U(y'))=0$. These conditions 
are satisfied if $C$ is the Auslander-Reiten quiver of a locally representation-finite category or the preprojective or preinjective component of the path algebra of a quiver without
oriented cycles or double arrows.

The next result is fundamental for all what follows. It is a slight generalization of sections 2.2 and 2.3 of   \cite{Rie1}. 

\begin{theorem} Using the above notations and assumptions we have: \begin{enumerate}

\item There is a $k$-linear functor $F:k(\tilde{C}) \rightarrow ind\,C$ mapping any point $\tilde{c}$ in $\tilde{C}$ to
 $U(\pi(\tilde{c}))$ and 
any arrow to an irreducible morphism.  

\item Each functor with these properties is a covering functor.

\end{enumerate}

\end{theorem}

Proof: For the first part  we have to construct  a functor $F'$ from the path category $k\tilde{C}$ to $ind \,C$ that vanishes on the mesh-sums. We map each object $\tilde{c}$ 
to $U(\pi\tilde{c})$.
 To define $F'$ on the arrows we use the morphism of translation quivers $\kappa:\tilde{C} \rightarrow \Z{\bf A}_{2}$ and we denote the points of $\Z{\bf A}_{2}$ by the integers as in  
\cite[section 1]{BG}.
We say that an arrow
 $\alpha:\tilde{x} \rightarrow \tilde{y}$ ends at level $i$ if $\kappa\tilde{y}=i$. An arrow ending at level 1 is sent to an arbitrary irreducible morphism between the images of the 
start point 
and the end point. Let $\alpha:\tilde{x} \rightarrow \tilde{y}$ be an arrow ending at level 2. If $\tilde{y}$ is projective we send $\alpha$ again to an arbitrary irreducible morphism. 
If $\tilde{y}=\tau^{-1}\tilde{x}$ let
 $\alpha_{i}:\tilde{x} \rightarrow \tilde{x}_{i}$, $1 \leq i \leq r$, be all the arrows starting at $\tilde{x}$. Then the morphism
 $ F'\tilde{x} \rightarrow \oplus_{i}\,F'\tilde{x}_{i}$ with the $F'\alpha_{i}$ as components is minimal left almost split and so there is by well-known arguments an almost split sequence
$$ 0 \rightarrow F'\tilde{x} \rightarrow \oplus_{i}\,F'\tilde{x}_{i} \rightarrow F'\tilde{y} \rightarrow 0.$$ We choose for $F'(\sigma^{-1}\alpha_{i})$ the component
 $F'\tilde{x}_{i} \rightarrow F'\tilde{y}$ in the almost split sequence. Then $F'$ vanishes on the mesh-sum $\mu_{\tilde{x}}$. We continue in the same way  on all arrows ending at 
level $3$ and so on. Dually we define $F'$ on all arrows starting at negative level and we get a functor $F'$ vanishing on all mesh-sums inducing the wanted functor 
$F:k(\tilde{C}) \rightarrow ind\,C$. In the sequel we write $\overline{w}$ for the image of a morphism $w$ of $k\tilde{C}$ in $k(\tilde{C})$.

  For the second part it is up to duality to show  that $F_{(\tilde{x},y)}$ is bijective for all $\tilde{x} \in \tilde{C}$ and $y \in C$. The surjectivity follows for $x=F\tilde{x}$ from 
 $\mathcal{R}^{n(x,y)}(U(x),U(y))=0$ by well-known arguments. The injectivity is the subtle point.

 We fix a $y$ in $C$ and we look at the set $\mathcal{M}$ of all triples $(\tilde{x},\tilde{z},(\phi_{\tilde{y}}))$ with $\tilde{x}$ and  $\tilde{z}$ in $\tilde{C}$, $\pi\tilde{z}=y$,
$(\phi_{\tilde{y}})$ in $\oplus_{\pi\tilde{y}=y} \,k(\tilde{C})(\tilde{x},\tilde{y})$ with $\sum_{\pi\tilde{y}=y} F(\phi_{\tilde{y}}) =0$ and $\phi_{\tilde{z}} \neq 0$.
We have to show that $\mathcal{M}$ is empty. If not we choose a triple $(\tilde{x},\tilde{z},(\phi_{\tilde{y}}))$ such that $d=\kappa(\tilde{z}) - \kappa(\tilde{x})$ 
is minimal. Because $\phi_{\tilde{z}} \neq 0$ we have $d\geq 0$. For $d=0$ we get $\tilde{z}=\tilde{x}$, $\phi_{\tilde{z}}= a\,id_{\tilde{z}}\neq 0$ and so 
$\sum_{\pi\tilde{y}=y} F(\phi_{\tilde{y}}) =0$ implies that $id_{y}$ belongs to the radical of $End_{A} U(y)$ which is impossible. The same contradiction follows 
if no arrow starts at $\tilde{x}$ in $\tilde{C}$.

Thus we have $d>0$ and there are $r>0$ arrows $\alpha_{i}:\tilde{x} \rightarrow \tilde{x}_{i}$, $1 \leq i \leq r$, starting at $\tilde{x}$. Then we have 
$\phi_{\tilde{y}}= \sum_{i} \psi_{\tilde{y},i}\overline\alpha_{i}$ for each $\tilde{y}$ with some appropriate $\psi_{\tilde{y},i}$ whence 
$$0=\sum_{\tilde{y}}F\phi_{\tilde{y}}= \sum_{\tilde{y}}(\sum_{i} F\psi_{\tilde{y},i}F\overline\alpha_{i})=\sum_{i}(\sum_{\tilde{y}} F\psi_{\tilde{y},i})F\overline\alpha_{i}.$$ 
Denoting by 
$f_{\tilde{x}}: F\tilde{x} \rightarrow \oplus_{i}\,F\tilde{x}_{i}$ the  minimal left almost split morphism with components $F\overline{\alpha_{i}}$ and by $\psi: \oplus_{i}\,F'\tilde{x}_{i} \rightarrow U(y)$ the 
homomorphism with components
 $\sum_{\tilde{y}} F\psi_{\tilde{y},i}$, $1 \leq i \leq r$ one obtains $\psi \,f_{\tilde{x}}=0$. 

If $\tilde{x}$ is injective then $f_{\tilde{x}}$ is epi whence $\psi =0$. Now we have $0\neq \phi_{\tilde{z}}= \sum_{i} \psi_{\tilde{z},i}\overline\alpha_{i}$ and so there is an index $j$ with
 $\psi_{\tilde{z},j}\neq 0$. But from $\psi=0$ we get $\sum_{\tilde{y}} F\psi_{\tilde{y},j}=0$. Thus $(\tilde{x}_{j},\tilde{z},(\psi_{\tilde{y},j}))$ belongs to $\mathcal{M}$ and we get
 the contradiction $\kappa(\tilde{z})-\kappa (\tilde{x}_{j}) < d$.

If $\tilde{x}$ is not injective we consider  the minimal right almost split morphism $g:\oplus F\tilde{x}_{i} \rightarrow  F\tau^{-1}\tilde{x}$ with the $F\overline{\sigma^{-1}\alpha_{i}}$ as components.
Because of $\psi f_{\tilde{x}}=0$  we  have $\psi = l \,g$ for some homomorphsm $l$. Using the surjectivity of $F_{(\tau^{-1}\tilde{x},y)}$ there is an element $(\chi_{\tilde{y}})$ in 
$\oplus_{\pi\tilde{y}=y} \,k(\tilde{C})(\tau^{-1}\tilde{x},\tilde{y})$ with $\sum_{\tilde{y}}F\chi_{\tilde{y}}= l$.Then we obtain for each $i$ the equality
$\sum_{\tilde{y}}F\psi_{\tilde{y},i} = \sum_{\tilde{y}}F\chi_{\tilde{y}} F\sigma^{-1}\alpha_{i}$. Now  by $$0 \neq \phi_{\tilde{z}}=\sum_{i}\psi_{\tilde{z},i}\overline{\alpha_{i}}=
\sum_{i}\psi_{\tilde{z},i}\overline{\alpha_{i}} - \chi_{z}\sum_{i}\overline{(\sigma^{-1}\alpha_{i})(\alpha_{i})} =
\sum_{i}(\psi_{\tilde{z},i}- \chi_{\tilde{z}}\overline{\sigma^{-1}\alpha_{i}})\overline{\alpha_{i}}$$ there is an index $j$ with 
$\psi_{\tilde{z},j}- \chi_{\tilde{z}}\overline{\sigma^{-1}\alpha_{j}} \neq 0$. Then $(\tilde{x}_{j},\tilde{z},(\psi_{\tilde{y},j}- \chi_{\tilde{y}}\overline{\sigma^{-1}\alpha_{j}}))$ belongs to
$\mathcal{M}$. The contradiction is again $\kappa(\tilde{z})-\kappa (\tilde{x}_{j}) < d$.
  \hspace{0.1cm}q.e.d.

Observe that the construction of $F$ involves many choices.

\subsection{Stable components of locally representation finite categories}

The first application of theorem 2 is the 'Hauptsatz' from  \cite{Rie1}.

\begin{theorem}
 Let $C$ be a connected component of the stable part $\Gamma^{s}_{A}$ of the Auslander-Reiten quiver of a locally representation-finite connected category $A$. 
Then $C$  is isomorphic to the quotient $\Z T/G$   where $T$ is a Dynkin-quiver and $G$ an admissible automorphism group which is cyclic.
\end{theorem}

Proof:  Consider  the universal covering $\pi:\tilde{\Gamma}_{A} \rightarrow \Gamma_{A}$ of $\Gamma_{A}$,  the induced covering $\pi':C':=\pi^{-1}C \rightarrow C$ and the universal covering 
$\rho:\tilde{C'} \rightarrow C'$ 
of $C'$. Here $C, C'$ and $\tilde{C'}$ are all connected stable translation-quivers and $\tilde{C'}$ is simply connected whence isomorphic to $\Z T$ for some tree $T$ by theorem 1. Because 
 $\pi ' \circ \rho$ is a covering it is the universal covering of $C$ which is therefore isomorphic to $\Z T/G$ where $G$ is the fundamental group of $C$.

By theorem 2 we  have a covering functor $F:k(\tilde{\Gamma}_{A}) \rightarrow ind\,A$. Since $ind\,A$ is locally bounded the same holds for $\tilde{\Gamma}_{A}$, its stable sub translationquiver $C'$
 and finally for $\tilde{C'} \simeq \Z T$. It remains to show that $T$ is a Dynkin quiver.

 If $T$ contains a Euclidean subtree $T'$ then $\Z T'$ and also $\N T'$ are full 
translation subquivers of $\Z T$ whence locally bounded. By  theorem 1 we have an equivalence $\N T' \rightarrow ind\,P$ where $P$ is the preprojective component of 
the path category to the opposite quiver of $T'$. Here $P$ is infinite 
and so 
$\N T'$ is not locally bounded. Thus $T$ contains no Euclidean subtree and so it has at most one branching point.

If $T$ is infinite then $\Z T$  contain an infinite line $L$ where all arrows point away from the starting point. This is a full translation subquiver which is not locally bounded. Thus $T$ is
 a finite tree not containing a Euclidean subtree whence  a  Dynkin-quiver. 

In   \cite[section 4]{Rie1} it is shown that only cyclic groups occur as admissible automorphism groups.  
\hspace{0.1cm}q.e.d.

Note that the case $G=1$ is allowed in the theorem above. This does not follow by the method of Todorov in \cite{T} resp. Happel/Preiser/Ringel in \cite{HPR} which requires a periodic module
 in the component, but
 it follows also from the results of Bautista/Brenner in \cite{BB}.

\subsection{Auslander-categories and standard-categories}

A $k$-category $D$ is an Auslander-category if it is equivalent to the full subcategory $ind\, A$ of the indecomposables over some locally representation-finite category $A$. An nice result 
of Auslander gives a homological characterization of this property.

\begin{theorem} A $k$ -category $D$ is an Auslander-category iff  $D$ is locally bounded, has global dimension at most 2 and any projective $p$ in $D$ admits an exact sequence
 $0 \rightarrow p \rightarrow i_{0} \rightarrow i_{1}$ where the $i_{k}$ are 
projective and injective. 
 
\end{theorem}

A proof of this and of the next result is in section 3 of \cite{BG}. 
\begin{theorem}
 Let $F:M \rightarrow N$ be a covering functor between square-free locally finite dimensional k-ctagories. Then $M$ is an Auslander-category iff $N$ is so.
\end{theorem} 

We apply this to a locally representation-finite category $A$ with Auslander-Reiten quiver $\Gamma_{A}$ and universal cover $\pi:\tilde{\Gamma}_{A} \rightarrow \Gamma_{A}$.
 By a result of Bautista mentioned in \cite{BB} $ind\,A$, $k(\Gamma_{A})$ and $k(\tilde{\Gamma}_{A})$ are all square-free and we have covering functors $F:k(\tilde{\Gamma}_{A}) \rightarrow ind\,A$
 and $k(\pi):k(\tilde{\Gamma}_{A}) \rightarrow k(\Gamma_{A})$.

Thus there are equivalences $k(\tilde{\Gamma}_{A}) \simeq ind\,\tilde{A}$ and $k(\Gamma_{A}) \simeq ind\,A^{s}$ for some appropriate locally representation-finite categories $\tilde{A}$ and $A^{s}$.
The next result shows that both are uniquely determined up to equivalence by $A$, that $A$ and $A^{s}$ have the same Auslander-Reiten quiver and that $\tilde{\Gamma}_{A}$
 is the Auslander-Reiten quiver 
of $\tilde{A}$. 
We call $A^{s}$ the standard-category of $A$, $\tilde{A}$ the universal cover of $A$ and $A$ standard if $A\simeq A^{s}$. 

\begin{proposition} Let $B$  be locally representation-finite and let $\Gamma$ be a translation quiver. Then the following are equivalent:
\begin{enumerate}
 \item There is an equivalence $F:k(\Gamma) \rightarrow ind\,B$. 
\item $B$ is equivalent to the full subcategory
$Pro(k(\Gamma))$ of $k(\Gamma)$ consisting of the projective points and the two translation quivers $\Gamma$ and  $\Gamma_{B}$ are isomorphic. 
\end{enumerate}

\end{proposition}
Proof: We show that the first assertion implies the second. Of course $F$ induces an equivalence between the categories $fun (k(\Gamma))$ and $fun(ind\,B)$ of finitely
 generated contravariant functors with values in $mod\,k$ which  preserves  the simple functors and their minimal projective resolutions. Now the simple functors correspond to the points $x$ 
of $\Gamma$ resp. the ( isomorphism classes of ) indecomposables $U$ in $mod\,B$. The easy  lemma 2.6 in \cite{BG} describes the minimal projective resolutions of the simple functors $k_{x}$ in
 $fun (k(\Gamma))$. The projective dimension of $k_{x}$ is at most 1 iff $x$ is a projective point  and otherwise the minimal projective resolution of $k_{x}$ is induced by the mesh stopping at $x$.
In $fun(ind\,B)$ the situation is completely analogous by basic properties of almost split sequences. The simple functor $S_{U}$ corresponding to the indecomposable $U$ has projective dimension
 at most 1 iff $U$ is projective and otherwise the minimal projective resolution 
is induced by the almost split sequence ending at $U$. Thus $F$  induces an equivalence $Pro(k(\Gamma)) \simeq Pro(ind\,B)$ with $B \simeq Pro(ind\,B)$ by Yonedas lemma and $F$ induces also
 an isomorphism $\Gamma \simeq \Gamma_{B}$.

Reversely, using theorems 2 and 5 we see that $k(\Gamma_{B})$ is equivalent to $ind\,C$ for some locally representation-finite $C$. By the direction just shown we know that $C$ is
 equivalent  to $Pro(k(\Gamma_{B}))$. From $k(\Gamma) \simeq k(\Gamma_{B})$ we obtain $C\simeq Pro(k(\Gamma_{B})) \simeq Pro(k(\Gamma)) \simeq B$ whence $k(\Gamma)\simeq ind\,B$.\hspace{0.1cm}  
    q.e.d.

The definitions of $\tilde{A}$ and $^{s}A$ strongly depend on the existence of a covering functor $F:k(\tilde{\Gamma}_{A}) \rightarrow  ind\,A$ for the locally representation finite category $A$.
 Of course the nice presentation of $ind\,A^{s} \simeq k(\Gamma_{A})$ by mesh relations 
should correspond to a nice presentation of $A^{s}$ by quiver and relations and the relation between $A$ and $A^{s}$ should be analyzed. This will be discussed in section 9.

 For the universal cover $\tilde{A}$ one gets a 'standard presentation'  using  two observations of mine \cite{B1,B2}. The easy proof is given in \cite[section 5.3]{BLR}. 

\begin{theorem} The universal cover $\tilde{A}$ of a locally representation finite category has a presentation $\phi:k\tilde{Q} \rightarrow \tilde{A}$ where all paths 
between two points have the same image.

\end{theorem}

\section{The classification in the simply connected case}
\subsection{Configurations}
Let $A\neq k$ be a  connected selfinjective locally representation-finite category with simply connected Auslander-Reiten quiver $\Gamma_{A}$. 
It is well-known that the almost split
 sequences having a projective-injective $P$ as a direct summand of the
 middle term look like $0 \rightarrow rad\,P
\rightarrow P \oplus rad\, P/soc\,P \rightarrow P/soc\,P \rightarrow 0$.  If $rad\,P/soc\,P =0$ for one $P$ one gets $\Gamma_{A} \simeq \Z A_{2}$ and $\Gamma_{A}^{s} \simeq \Z A_{1}$. 
Thus  the stable translaton quiver $\Gamma_{A}^{s}$ is always connected ( as a translation quiver  ) 
and even simply connected. By theorem 2 we know that $\Gamma_{A}^{s}$ is isomorphic to $\Z T$ for a uniquely determined Dynkin diagram $T$ that we call the tree class of $A$.

We define for  any subset $C$ of $\Z T$ a  
translation quiver
 $\Z T_{C}$ by adding for each $c \in C$ a new point $c^{\ast}$ which is projective and 
injective and arrows 
$c \rightarrow c^{\ast}$ , $c^{\ast} \rightarrow \tau^{-1}c$. 
 Here $C$ is called a configuration if $\Z T_{C}$ is the Auslander-Reiten
quiver of a locally representation-finite selfinjective category $A(C)$.

To analyze these categories we fix some definitions and notations most of which are familiar.
 For any point $x$ in $T_{0}$
 one has permutations 
 $s^{+}_{x}$
 resp. $s^{-}_{x}$ of $\Z T$ mapping any $(j,x)$ to $(j+1,x)$ resp. to $(j-1,x)$ and leaving all points $(i,y)$ with $ y \neq x$ fixed.
Let $S$ be a section in  $\Z T$  i.e. 
a full connected subquiver such that for each $x \in T$ there is exactly one $l$ with $x(S)=(l,x) \in S$. If $x(S)$ is a source in  $S$ we say that $x$ is a source in $S$. 
Then the full subquiver 
$s^{+}_{x}S$ supported by $s^{+}_{x}S_{0}$ is again a section. For sinks the situation is analogous and any two sections can be transformed into each other through several appropriate $s_{x}^{+}$ resp.
$s_{y}^{-}$. Also a $(+)$-admissible enumeration $x_{1},x_{2},\ldots ,x_{r}$ of $T_{0}$  exists for any $S$ i.e. a sequence such that $x_{1}$ is a source in
$S$, $x_{2}$ a source in $s^{+}_{x_{_{1}}}S$ and so on. Then  
$\tau^{-1}S=s^{+}_{x_{r}} s^{+}_{x_{r-1}} \ldots s^{+}_{x_{1}S}$ is obvious and corresponding  results hold for sinks.

Given any subset $C$ of $\Z T$ the set $(\Z T_{C})_{0}$ is partially ordered by  $x \preceq y$
 provided there is a path from $x$ to $y$ in $\Z T_{C}$. For a section $S$ in $\Z T$ and a point $x$ in $\Z T_{C}$ we write $x \preceq S$  if $x \preceq y$ holds for some $y \in S$. Similarily
 $S \preceq
x $ is defined. If we have  a second section $S'$ the interval $[S',S]$ resp. $[S',S[$ is the set of all $x$ with $S' \preceq x \preceq S$ resp. $x \prec S$. Here  $x \prec S$ means 
$x \preceq S$ but $x \not\in S$.

\vspace{1cm}
\setlength{\unitlength}{0.4cm}
\begin{picture}(30,6)
\put(0,6){\line(1,0){30}}

\put(0,2){\line(1,0){30}}

\put(8,5){\line(1,1){1}}
\put(8,5){\line(1,-1){2}}
\put(10,3){\line(-1,-1){1}}

\put(12,5){\line(1,1){1}}
\put(12,5){\line(1,-1){2}}
\put(14,3){\line(-1,-1){1}}

\put(16,5){\line(1,1){1}}
\put(16,5){\line(1,-1){2}}
\put(18,3){\line(-1,-1){1}}

\put(8,6.3){$\tau^{i}S$}
\put(12,6.3){$S$}
\put(16,6.3){$\tau^{-i}S$}

\put(13,-0.5){figure 1}
\end{picture}
\vspace{0.5cm}

Several times we will work in the  situation indicated in figure 1 above:
We have a locally representation-finite category $A$ with simply connected Auslander Reiten quiver
 $\Gamma_{A}$ which is 
 embedded into some $\Z T_{C}$ and contains the interval  $[\tau^{i}S, \tau^{-i}S]$ for some section $S$ and some $i> 0$.  To the left of $\tau^{i}S$ and to the right of $\tau^{-i}S$
there might be  infinitely many points of $\Gamma_{A}$ resp. $C$ or not. In particular, $C$ does not have to be a
 configuration but it is so 'locally' 
between $\tau^{i}S$ and $\tau^{-i}S$.

By theorem  2 there is  a $k$-linear equivalence $M:k(\Gamma_{A}) \rightarrow ind\,A$ which is 'the identity' on the objects and we define the vector $d(S)=d(S,A) \in \N ^{T_{0}}$ by $d(S)(x)= dim\,M(x(S))$. Similarly 
we define $d(S')$ for any section of $\Z T$ that belongs to $\Gamma_{A}$.   
Finally for any subset $X$ of 
$\Gamma_{A}$ the module $M(X)$ is the direct sum of the modules $M(x)$ with $x \in X$. Given points $x,y$ in $\Gamma_{A}$ with $ x \preceq S \preceq y$ any homomorphism $f:M(x) \rightarrow M(y)$
 factors through some power $M(S)^{n}$ because any path $p$ from $x$ to $y$ hits $S$ and $f$ is a linear combination of such $M(p)$.

The next easy proposition is important.

\begin{proposition}
 Keeping all the above assumptions and notations we have:
\begin{enumerate}
 \item $d(S)$ determines $C\cap [S,\tau^{-i}S]$ and $d(\tau^{-i}S)$. Similarly $d(S)$ determines also $C\cap [\tau^{i}S,S]$ and $d(\tau^{i}S)$.
\item $d(\tau^{i}S)= d(S)$ implies $C\cap [\tau^{i}S,S]= \tau^{i}(C \cap[S,\tau^{-i}S])$ and $d(S)=d(\tau^{-i}S)$. The dual also holds.
\item $C\cap [\tau^{i}S,S]= \tau^{i}(C \cap[S,\tau^{-i}S]$ and $Hom_{A}(M(\tau^{i}S),M(S))=0$ implies $d(S)=d(\tau^{-i}S)$.
\end{enumerate}

\end{proposition}

Proof:  Recall that  $x_{1}(S)$ is a
 source in $S$. We set
$$s:= -d(S)(x_{1}) + \sum _{x_{1}(S) \rightarrow y(S)}d(S)(y).$$ Two cases can occur. For $s > 0$ the point $x_{1}(S)$ is not in $C$ and we have $d(s^{+}_{x_{1}}S)(x_{1})=s$. For $s < 0$ the 
point $x_{1}(S)$ belongs to $C$ and we have $d(s^{+}_{x_{1}}S)(x_{1})=d(S)(x_{1})$. For $y \neq x_{1}$ the entries of $d(S)$ and $d(s_{x_{1}}^{+}S)$ coincide in both cases.  
Thus we always get $d(s^{+}_{x_{1}}S)$ and we see whether $x_{1}(S)$ belongs to $C$ or not.
We apply the analogous argument to $x_{2}$ and $d(s^{+}_{x_{1}}S)$ and so on  until we reach $\tau^{-i}S$  thereby proving the first half of part i).
To go in the other direction we argue with sinks. 

Part ii) is obvious because to go from $\tau^{i}S$ to $S$ or to go from $S$ to $\tau^{-i}S$ affords the same calculations up to a shift in the indices.

The assumption $Hom(M(\tau^{i}S),M(S))=0$ implies that $d(S)(x)$ is for any $x$ the sum of the dimensions of the $Hom(M(c^{\ast}),M(x(S)))$ with $c$ is in $C\cap [\tau^{i}S,S]$. 
By the other assumption  the two mesh-categories consisting in the intervals $[\tau^{i}S,S]$ and $[S,\tau^{-i}S]$ are isomorphic. So we obtain also $Hom(M(S),M(\tau^{-i}S))=0$ and $d(S)=d(\tau^{-i}S)$
 follows.  \hspace{0.1cm}q.e.d.

\subsection{Patterns}

Again we choose a section $S$ in $\Z T$. An $S$-pattern is a subset $P$ of $\Z T$  such that there is a simply connected algebra $B=B(P)$ and an embedding of $\Gamma_{B}$ into $\Z T$ with
 $S \subseteq \Gamma_{B}  \subseteq \Z T$ such that $P$ is the set of projective points. Then the dimension vector $d(P) \in \N ^{T_{0}}$ is defined by $d(P)(x) = dim\,M(x(S))$. Here $M(x(S))$ is the
indecomposable $B$-module corresponding to $x(S)$. Of course $P$ determines $d(P)$  and then $B(P)$  by the usual knitting forward and by theorem 6. Reversely,  $d(P)$ determines 
 $P$ by knitting backward.

\begin{theorem}
 Let $T$ be a Dynkin-tree and let $S$ be  a section in $\Z T$. For a vector $d \in \N^{S}$ there is a configuration $C$ with $d=d(S,A(C))$ iff there is a pattern $P$ with $d=d(S,P)$.

\end{theorem}

Thus we get the bijection $P \mapsto d(P)=d(S,C) \mapsto C$ between the finite set of $S$-patterns  and the set of configurations.
 The proof provides a slightly cumbersome algorithm to obtain the configuration from the pattern by knitting and knotting an infinite periodic carpet starting with the pattern.
Fortunately there is a direct way to get the category $A(C)$ from $B(P)$ which we describe in subsection 4.6.

 The proof uses only basic  Auslander-Reiten theory, one-point extensions and
 properties of the starting functions $s_{x}$ in the mesh category $k(\Z T)$ of a Dynkin-quiver as explained in \cite{G4}. Here $s_{x}:(\Z T)_{0} \rightarrow \N$ is defined by
 $s_{x}(y)= dim\,k(\Z T)(x,y)$.

\subsection{Some calculations in the mesh-categories of $\Z T $ and $\Z T_{C} $}

 A morphism $f:X \rightarrow Y$ in an additve   category is  right-complete if $f$ is not zero but the composition $fg$ with any non-retraction $g$ is.
Similarly one defines left-
complete and complete morphisms.  
 In the category $mod\,A$ of finite dimensional modules over a locally bounded category $A$ the projective cover $p_{x}:P_{x} \rightarrow k_{x}$ of each simple $k_{x}$ is left complete and 
 and the injective hull  $i_{x}:k_{x} \rightarrow I_{x}$ is right complete. A morphism is complete iff it is a non-zero element in the two-sided socle of $mod\,A$. 

Any non-zero morphism $f:X \rightarrow Y$ between indecomposables can be prolongated to a complete 
morphism $hfg=i_{x}p_{x}$ whenever $k_{x}$ occurs in the top of the image of $f$. Namely
 $f=f''f'$ with $f':X \rightarrow Im\,f$ and $f'':Im \,f
\rightarrow Y$ the canonical morphisms. Let $p:Im\,f \rightarrow k_{x}$ be a non-zero map. Then we find morphisms $g$ and $h$ with $p_{x}=pf'g$ and $i_{x}p=hf''$. The  dual argument applies to any
 simple submodule $k_{y}$ of $Im \,f$ so that $f$ can eventually be  'completed' in many different ways, but $i_{x} p_{x}$ is the only complete morphism starting in $P_{x}$ 
up to multiplication by a non-zero scalar.

For any  Dynkin diagram $T$ we define $L(T)=c(T)-1$ where $c(T)$ is the Coxeter number. Thus we have  $L(A_{n})=n$, 
$L(D_{n})=2n-3$, $L(E_{6})=11$, $L(E_{7})=17$ and $L(E_{8})=29$.

The next well-known  facts from \cite{G4} are very important for us. 

\begin{lemma} Let $T$ be a Dynkin diagram. Then there is an automorphism $\nu_{T}$  of $\me$ with the following properties:
\begin{enumerate}
\item For any $f \neq 0$ in $\me(x,y)$ one has $\nu_{T}x \preceq x \preceq y \preceq \nu_{T}^{-1}x$.
\item $0 \neq f:x \rightarrow y$ is right-complete if and only if $x= \nu_{T}y$ and left-complete if and only if $y= \nu_{T}^{-1}x$.
\item For any $x$ any path from $x$ to $\nu_{T}^{-1}x$ has length $L(T)-1$.
\item $\nu_{T}^{2}=\tau^{L(T)-1}$.
\end{enumerate}

\end{lemma}

Now we consider the same situation as in subsection  1.1.
We have a locally representation-finite category $A$ with simply connected Auslander Reiten quiver
 $\Gamma_{A}$ which is 
 embedded into some $\Z T_{C}$ and contains the interval  $[\tau^{3L}S, \tau^{-3L}S]$ for some section $S$ and $L=L(T)$. ( see figure 2 ).
\vspace{0.5cm}

\setlength{\unitlength}{0.4cm}
\begin{picture}(30,7)
\put(0,6){\line(1,0){30}}

\put(0,2){\line(1,0){30}}

\put(4,5){\line(1,1){1}}
\put(4,5){\line(1,-1){2}}
\put(6,3){\line(-1,-1){1}}

\put(8,5){\line(1,1){1}}
\put(8,5){\line(1,-1){2}}
\put(10,3){\line(-1,-1){1}}

\put(12,5){\line(1,1){1}}
\put(12,5){\line(1,-1){2}}
\put(14,3){\line(-1,-1){1}}

\put(12,5){\line(1,1){1}}
\put(12,5){\line(1,-1){2}}
\put(14,3){\line(-1,-1){1}}

\put(16,5){\line(1,1){1}}
\put(16,5){\line(1,-1){2}}
\put(18,3){\line(-1,-1){1}}

\put(20,5){\line(1,1){1}}
\put(20,5){\line(1,-1){2}}
\put(22,3){\line(-1,-1){1}}

\put(24,5){\line(1,1){1}}
\put(24,5){\line(1,-1){2}}
\put(26,3){\line(-1,-1){1}}

\put(28,5){\line(1,1){1}}
\put(28,5){\line(1,-1){2}}
\put(30,3){\line(-1,-1){1}}

\put(4,-0.3){$\tau^{3L}S$}\put(8,-0.3){$\tau^{2L}S$}
\put(12,-0.3){$\tau^{L}S$}
\put(16,-0.3){$S$}
\put(20,-0.3){$\tau^{-L}S$}
\put(24,-0.3){$\tau^{-2L}S$}
\put(28,-0.3){$\tau^{-3L}S$}
\put(13,7){figure 2}

\end{picture}\vspace{0.5cm}

\begin{proposition} 
We keep the notations and assumptions. Then we have:
\begin{enumerate}

 \item For $c$ in $C \cap [\tau^{L}S,S]$ the point $d= \tau^{-L}c$ belongs to $C$ and $M(d^{\ast})$ is the injective hull of the top of $M(c^{\ast})$. Dually for $d \in [S,\tau^{-L}S]$ 
the point $\tau^{L}d$ belongs to $C$ and $M(c^{\ast})$ is the projective cover of the socle of $M(d^{\ast})$. In particular $ [\tau^{L}S,S]$ and  $[S,\tau^{-L}S]$ are just translates of each other.
 \item We have $Hom_{A}(M(S),M(\tau^{-L}S))=0$.
\end{enumerate}
\end{proposition}
 Proof: For $x,y$ in the full translation subquiver $[\tau^{3L}S,\tau^{-3L}S]$ we have 
$$Hom_{A}(M(x),M(y))\simeq k(\Z T_{C})(x,y) \,\,\mbox{and} \,\,\underline{Hom}_{A}(M(x),M(y))\simeq k(\Z T)(x,y).$$
 Here  we consider the quotient annihilating maps factoring through  arbitrary  direct sums of some $M(c^{\ast})$ with $c \in [\tau^{3L}S,\tau^{-3L}S]$.  
Furthermore $$Hom(M(x),M(s))= \underline{Hom}(M(x),M(s))$$ and its dual hold always for simple $M(s)$ and indecomposable non-projective $M(x)$. 
The statements follow easily from this and the last lemma.

 Namely for $c \in C \cap [\tau^{L}S,S]$ and $s$ with $M(s)=top M(c^{\ast})$ we find  $0 \neq f \in Hom(M(\tau^{-1}c),M(s))=\underline{Hom}(M(\tau^{-1}c),M(s))$ which implies
 $s \preceq \nu_{T}(\tau^{-1}c)\preceq \tau^{-L}S$. Since $f$ is right-complete in $\underline{mod}\,A$ we get $\tau^{-1}c = \nu_{T}^{-1}s$. Similarly let $M(i')$ be the radical of the injective hull
of $M(s)$. Then there is $0 \neq g \in Hom( (M(s),M(i'))= \underline{Hom}(M(s),M(i'))$ which implies $i' \preceq \tau^{-2L}S$. Since $g$ is left-complete we obtain $s=\nu_{T}^{-1}i'$. The rest
follows from $\nu_{T}^{2}=\tau^{-L+1}$ and duality.

  Suppose $0\neq f:M(x) \rightarrow M(y)$ with $x \in S,y \in \tau^{-L}S$. Then we find as in the beginning of section 1.3 a simple $M(s)$ with projective cover $P=M(p)$ and injective hull $I=M(i)$
such that $Hom(M(x),M(s))\neq 0$, $Hom(P,M(x)) \neq 0$ and $Hom(M(y),I) \neq 0.$ Thus we get $S \preceq s$. From $Hom(P/soc\,P,M(s))\neq 0$ we infer $\tau^{L}S \preceq p$ and so  $p=c^{\ast}$
 for some $c \in C\cap [\tau^{L}S,S]$. This implies  $I=M(d^{\ast})$ with $d \prec \tau^{-L}S$ contradicting $Hom(M(y),I) \neq 0$.\hspace{0.1cm}q.e.d. 
\vspace{0.2cm}

The assumptions of the proposition are obviously satisfied for $A(C)$ coming from a configuration $C$. Thus $C$  is stable under multiplication 
with $\tau^{L(T)}$ and $\tau^{-L(T)}$. Therefore the automorphism $\tau^{L(T)}$ of $\Z T$ can be extended to an automorphism $\nu=\nu_{C}$ of  $\Z T_{C}$ inducing 'Nakayama'-automorphisms 
of $k(\Z T_{C}) $  and 
 $A(C)$  still denoted by  $\nu$.

The configuration  $C$ is known once we know  for any section $S$ the fundamental domain $P(S,C):=C \cap [\nu S,S[$. Any 
section
 $S'$ is obtained from a fixed section $S$ by applying several appropriate $s^{+}_{x}$ or $s^{-}_{y}$. For a source we have 
$P(s^{+}_{x}S,C)=P(S,C)$ or $P(s^{+}_{x}S,C)=(P(S,C)\setminus \{x\}) \cup \{\nu ^{-1} x\}$  depending on $x \not\in C$ or $x \in C$. The analogous result holds for sinks and so we get that 
$p(C)=|P(S,C)|$ is independent of $S$. In fact we always have $p(C)=|T_{0}|$ as we are going to show.

\subsection{The pattern to a configuration}

In this subsection $C$ is a configuration of $\Z T$, $S$ a section, $A=A(C)$ and $d=d(S,C)$. We consider the full subcategory $B=B(S,C)$ of $A(C)$ supported by the points $c^{\ast}$ with
 $c \in P(S,C)$. 
 Here $B(S,C)$ strongly depends on $S$ but it has always $p(C)$ points. 

 $B$ is a convex subcategory and we identify as usual the $B$-modules with the $A$-modules vanishing outside of
$B$. Thus for each indecomposable $B$-module $U$ there is a uniquely determined point $u$ in $\Gamma_{A}$ with $U\simeq M(u)$ but the map $U \mapsto M(u)$ is 
not compatible with the translations $\tau_{B}$ and $\tau_{A}$. Since $A$ is locally representation-finite and $\Gamma_{A}$ is simply connected the Auslander-Reiten quiver $\Gamma_{B}$ is a 
finite preprojective
 component which is even simply connected as follows from our next result:

\begin{proposition}Keeping the assumptions and notations $S$ is a section in $\Gamma_{B}$. Furthermore there is an $S$-pattern $P$ with $B(P)=B(S,C)$ and $d(P)=d(S,C)$.

\end{proposition}
 
Proof: From $Hom(M(\tau^{L}S),M(S))=0$ we get that the only $c \in C$ with $Hom(M(c^{\ast}),M(S)) \neq 0$ are those in $P(S,C)$. For any such  $c$ there is  an $f \neq 0$ in
 $Hom(M(c^{\ast}),M((\tau^{-L}c))^{\ast})$. Since $S \preceq \tau^{-L}c$ the morphism $f$ factors through some power of $M(S)$. We conclude that  $M(S)$ is a sincere $B$-module
and so the projective indecomposable $B$-modules $M(p)$ satisfy $p \preceq S$ and the injectives $M(i)$ satisfy $S \preceq i$. 
We have to show that the points of $S$ are representatives of the $\tau_{B}$-orbits and the arrows $\sigma_{B}$- orbits of the irreducible maps in $\Gamma_{B}$.

First assume that for some $x\neq y$ in $S$ the $\tau_{B}$-orbits intersect. Then we get  up to symmetry $M(y)=\tau_{B}^{k} M(x)$ with 
$k \geq 1$.  The almost split sequence $0 \rightarrow \tau_{A}M(x) \rightarrow X \rightarrow M(x) \rightarrow 0$ is a push out of the almost
 split sequence
$0 \rightarrow \tau_{B}M(x) \rightarrow Y \rightarrow M(x) \rightarrow 0$ by a non-zero morphism. Thus we find a chain of irreducible maps from $M(y)$ through $\tau_{B}M(x)$ and 
 $\tau_{A}M(x)$ to  $M(x)$. 
Therefore  $\tau_{A}M(x)$ and $M(x)$ lie on $S$ which is a contradiction. 

It follows that $|T_{0}| \leq p(C)$ and that $|T_{0}|=p(C)$ if and only if each $\tau_{B}$-orbit hits $S$. Thus it suffices to prove this if $S$
 contains only sinks and sources. If it does not hold not there are  indecomposables $U=M(u)$ and $V=M(v)= \tau_{B}^{-1}U$  with $u \prec S \prec v$. Let 
$$ 0 \rightarrow U \rightarrow \oplus M(x_{i}) \rightarrow V \rightarrow 0$$
 be the almost split sequence in $mod\,B$. For $x_{i} \prec S$ the irreducible map $M(x_{i}) \rightarrow V$ factors properly through $M(S)$ in $mod\,B$. Thus we have $S \preceq x_{i}$. Dually we 
find $x_{i} \preceq S$. So all $x_{i}$ belong to $S$. Suppose $x_{1}$ is a source in $S$. Let $0 \rightarrow M(x_{1}) \rightarrow \oplus M(y_{j}) \rightarrow \tau_{A}^{-1}M(x_{1}) \rightarrow 0$
be the almost split sequence. Then all $M(y_{j})$ are $B$-modules except possibly $M((x_{1})^{\ast})$ which maps only trivially to $M(v)$. Thus the map $M(x_{1}) \rightarrow M(v)$ is not irreducible 
in $mod \,B$, a contradiction. Dually one finds that $M(u) \rightarrow M(x_{1})$ is not irreducible if $x_{1}$ is a sink in $S$.

The irreducible morphisms occurring in $S$ are a fortiori irreducible in $m0d\,B$ and  we prove that any
 irreducible morphism $f:U=M(u) \rightarrow V=M(v)$ is in the $\sigma_{B}$-orbit of such a morphism. First assume that $u$ and $v$ are in $S$.
 Then we have $f=\sum g_{i}f_{i}$ where the $f_{i}:U \rightarrow X_{i}$ are the components
 of the minimal left almost split morphism $U \rightarrow \oplus M(x_{i})$ in $mod \,A$. For 
$g_{i}f_{i} \neq 0$ the module $M(x_{i}$ is in the convex section $S$ whence a $B$-module. So we have factored $f$ inside $mod\,B$ and so $f$ is isomorphic to some $f_{i}$.
In the general case  
 $u \prec S \prec v$ is not possible because 
then $f$ factors through some power of the $B$-module $M(S)$. Since the orbits of $U$ and $V$ hit $S$ we can apply an appropriate power of $\tau_{B}$ to $f$ until one of the modules is in $S$. 
If the other is not 
yet in $S$ we apply $\tau_{B}$ or  $\tau_{B}^{-1}$ to it. The resulting module belongs to the convex set $S$ because the orbit hits $S$.

Finally $M(S)$ is a $B$-module so that $d=d(S,C)$ is the vector of total dimensions of $B$-modules  on a section. Knitting backwards we get the wanted pattern $P$ of projectives and  $\Gamma_{B}$ embeds 
into $\Z T$. \hspace{0.1cm}q.e.d.

\subsection{ The configuration to a pattern}

 Starting with a  with a vector $d=d(P)$ coming from a pattern $P$ we will construct by knitting and knotting  a configuration $C$ with $d=d(C)$.
 Here  knotting means taking an appropriate one-point extension with the help of the next simple lemma:

\begin{lemma} 
 Let $A$ be an algebra and $A'$ a one-point extension by a non-zero module $U$.  Let $P$ be the projective cover of the new simple $S$. 
Then $P$ is also injective iff $U$ is an injective $A$-module with endomorphism algebra $k$.
\end{lemma}

Proof: Because $P$ is projective and $S$ is injective the exact sequence $0 \rightarrow U \rightarrow P \rightarrow S \rightarrow 0$ induces for any simple $A'$-module $T$ the exact sequence 
$ Hom_{A'}(T,S) \rightarrow Ext^{1}_{A'}(T,U) \rightarrow Ext^{1}_{A'}(T,P) \rightarrow 0$ and also 
$0=Hom_{A'}(S,P) \rightarrow k\simeq  Hom_{A'}(P,P) \rightarrow Hom_{A'}(U,P)\simeq  Hom_{A'}(U,U) \rightarrow  Ext^{1}_{A'}(S,P) \rightarrow 0.$
Using this and  $Ext^{1}_{A'}(T,U)=Ext^{1}_{A}(T,U)$ for all simples $T$ not isomorphic to $S$ the lemma follows. \hspace{0.1cm}q.e.d.\vspace{0.2cm}

Now suppose we have  a simply connected representation-finite algebra $A$ such that its Auslander-Reiten quiver $\Gamma_{A}$ is embedded 
into some $\Z T_{C}$ where $T$ is a Dynkin-tree and $C$ a finite set. Assume furthermore that there is a section $S$ in $\Z T$ which is a final section in $\Gamma_{A}$. This means that the points 
in $S$ 
are representatives of the $\tau_{A}$-orbits of the indecomposable injectives $I$ with $S\preceq I$ and the arrows of $S$ are representatives of the $\sigma_{A}$-orbits between these orbits.
Moreover
we need to know the vector $d$ giving the total dimensions of the $A$-modules in $S$. Then we take a source $x$ is $S$ and we move one step further. Whether we knit or knot depends on the number
$$s:=  - d(x) + \sum _{x \rightarrow y}d(y).$$

For $s > 0$ we knit as usual. So we set $A'=A$, $C'=C$, $S'=s^{+}_{x}S$ and we define $d'$ by $d'(x)=s$, $d'(y)=d(y)$ for $y\neq x$ .

For $s \leq 0$ we knot. Namely the $A$-module $U$ corresponding to $X$ is an injective with endomorphism algebra $k$. Thus in the one point extension $A'$ of $A$ by $U$ the new projective is also
 injective and $S'=s^{+}_{x}S$ is a final section of $\Gamma_{A'}$. Then $d'=d$ is the vector of the total dimensions of the $A'$- modules lying on $S'$ and one can go on to the right to find 
the remaining part of $\Gamma_{A'}$ by the usual 
knitting with 
the total dimensions of the modules. Since $S'$ is
 a Dynkin-quiver this stops after finitely many steps whence $A'$ is again a representation-finite algebra. It is even simply connected since any closed walk in $\Gamma_{A'}$ is null-homotopic.
 Finally
we set $C'=C \cup \{x\}$ and then the dashed data satisfy in both cases the original assumptions so that we can go on. 

\begin{proposition}
 Let $T$ be a Dynkin-tree, $S$ a section in $\Z T$ and $P$ an $S$-pattern.  Then there is a configuration $\tilde{C}$ with $d(S,\tilde{C})=d(P)$ and $B(S,\tilde{C})=B(P)$.
\end{proposition}

Proof:  By definition there is  an embedding of $\Gamma_{B(P)}$ into $\Z T$ such that
 $S \subseteq \Gamma_{B(P)}$. We take a $(+)$-admissible sequence $x_{1},x_{2}, \ldots ,x_{r}$ of $S$ and abbreviate $L(T)$ by $L$. 

We will construct the wanted configuration inductively by the procedure explained before the proposition. To start the induction we take $A=B(P)$, $C= \emptyset$, 
$S=S$ and $d=d(P)$. Then all assumptions are satisfied and we can choose $x=x_{1}$ as the first source. The procedure above ends in data $A',C',S',d'$ still satisfying the assumptions.
 Moreover $x_{2}$
is a source for $S'$. We can go on like this as often as we want to. 
Having done this $t\cdot L \cdot r$ times we denote by $A(t), C(t), S(t)=\tau^{-tL(T)}S, d(t)$ the final result of our knittings and knottings.

For $t= 6$  proposition 4 shows that  the sets $C(6)\cap[\tau^{-2L}S,\tau^{-3L}S]$ and $C(6)\cap[\tau^{-3L}S,\tau^{-4L}S]$ are translates of 
each other and  $Hom(M(\tau^{-2L}S),M(\tau^{-3L}S)=0$. By part 3 of proposition 3  we find $d(\tau^{-2L}S) = d(\tau^{-3L}S)$.
Then the remaining parts of proposition 1 imply that for all $j=0,1, \ldots 5$ we have $d(S) = d(\tau^{-(j+1)L}S)$ and   $Hom(M(\tau^{-jL}S),M(\tau^{-(j+1)L}S)=0.$
Again by the frst two parts of proposition 3 we find that 
$C(t)\cap[\tau^{-(j+1)L}S,\tau^{-(j+2)L}S]$ ia a  translates of $C(6)\cap[S,\tau^{-L}S]$ for all $t\geq 6$ and all $j \leq t-1$. It follows that any indecomposable module over $A(t)$ is a translate
of an indecomposable over 
$A(1)$. Thus there is natural number $N$ such that the dimension of any indecomposable over $A(t)$ is  bounded by $N$ for any $t$.

Now we define $\tilde{C}$ as $$\tilde{C}:=\bigcup_{j \in \Z} \tau^{jL}(C(6) \cap[S,\tau^{-L}S[).$$

To see that this is a configuration we have to show that the full subcategory $A(\tilde{C})$ of $k(\Z T)_{\tilde{C}}$ consisting of projective points $c^{\ast}$ with $c \in \tilde{C}$ is 
locally representation-finite with 
$(\Z T)_{\tilde{C}}$ as its Auslander-Reiten quiver. Since $\tilde{C}$ is by definition periodic we have $Hom(M(\tau^{iL}S),M(\tau^{(i-1)L}S)=0$ for all $i$ and so $A(\tilde{C})$ is locally bounded.
 The modules occurring in any almost split sequence of $A(\tilde{C})$ 
 live on a finite piece  which can be translated into some appropriate $A(t)$ so that it becomes an almost split sequence in
 $mod \,A(t)$ and so the Auslander-Reiten quiver of $A(\tilde{C})$ is $(\Z T)_{\tilde{C}}$.

The remaining assertions of the proposition hold by construction. \hspace{0.1cm}q.e.d.\vspace{0.2cm}

We illustrate the reasoning by an example of tree class $A_{7}$ namely by the algebra $B$ given by the full subquiver with points
 $a^{\ast},b^{\ast},c^{\ast},d^{\ast},e^{\ast},f^{\ast},g^{\ast}$ shown in figure 3 subject to the  relations 
 $c^{\ast} \rightarrow d^{\ast} \rightarrow f^{\ast}=0$ and $ e^{\ast}\rightarrow f^{\ast} \rightarrow g^{\ast}=0$.
In  figure 4 we start to calculate the Auslander-Reiten quiver of $B$ until we reach the
 equioriented section $S$ leading upwards containing the bold-faced numbers that define the vector $d=(1,4,3,2,4,3,4)$. To go further to the right we  knit or knot. 
The points where we have to knot are encircled and denoted by 
$a,b,c,d,e,f,g$.  After $7$ steps the calculations become periodic since $\nu$ is induced by $\tau^{7}$. Applying $\nu$ and $\nu^{-1}$ again and again to the encircled points 
we get the configuration $C$ to the vector $d$.

\setlength{\unitlength}{0.4cm}
\hspace{1.5cm}\begin{picture}(32,8)

\put(3,2){$a^{\ast}$}\put(3.6,2.3){\vector(1,0){1.5}}
 
\put(0,0){$b^{\ast}$}\put(0.6,0.3){\vector(1,0){1.5}}

\put(2,0){$c^{\ast}$}\put(2.6,0.5){\vector(3,2){2.3}}

\put(5,2){$d^{\ast}$}\put(5.6,2.3){\vector(1,0){1.5}}
\put(5,0){$e^{\ast}$}\put(5.6,0.5){\vector(1,1){1.5}}
\put(7,2){$f^{\ast}$}\put(7.6,2.3){\vector(1,0){1.5}}
\put(9,2){$g^{\ast}$}\thicklines\put(9.6,2.3){\vector(1,0){1.5}}\thinlines

\put(11,2){$\nu^{-1} a^{\ast}$}
\put(8,0){$\nu^{-1} b^{\ast}$}
\put(13,0){$\nu^{-1} e^{\ast}$}

\put(7.5,-1){figure 3}
\thicklines
\put(5.7,1.9){\vector(3,-2){2}}
\put(7.8,1.9){\vector(3,-1){4.5}}\thinlines
\end{picture}
\vspace{0.5cm}

\setlength{\unitlength}{0.4cm}
\hspace{1cm}\begin{picture}(32,8)
\put(4,5){1}
\put(1,2){2}
\put(0,3){1}
\put(2,3){1}

\put(1,0){\bf 1}
\put(2,1){\bf 4}
\put(3,2){\bf 3}
\put(4,3){\bf 2}
\put(5,4){\bf 4}
\put(6,5){\bf 3}
\put(7,6){\bf 4}\put(7.3,6.3){\circle{1}}

\put(3,0){3}\put(3.3,0.3){\circle{1}}
\put(4,1){2}
\put(5,2){1}
\put(6,3){3}
\put(7,4){2}
\put(8,5){3}
\put(9,6){4}

\put(5,0){3}\put(5.3,0.3){\circle{1}}
\put(6,1){2}
\put(7,2){4}
\put(8,3){3}
\put(9,4){4}
\put(10,5){5}
\put(11,6){1}

\put(7,0){3}
\put(8,1){5}
\put(9,2){4}
\put(10,3){5}
\put(11,4){6}\put(11.3,4.3){\circle{1}}
\put(12,5){2}
\put(13,6){1}

\put(9,0){2}\put(9.3,0.3){\circle{1}}
\put(10,1){1}
\put(11,2){2}
\put(12,3){3}
\put(13,4){6}
\put(14,5){5}
\put(15,6){4}

\put(11,0){2}
\put(12,1){3}
\put(13,2){4}
\put(14,3){7}
\put(15,4){6}
\put(16,5){5}\put(16.3,5.3){\circle{1}}
\put(17,6){1}

\put(13,0){1}
\put(14,1){2}
\put(15,2){5}
\put(16,3){4}
\put(17,4){3}
\put(18,5){5}
\put(19,6){4}\put(19.3,6.3){\circle{1}}
 
\put(15,0){1}
\put(16,1){4}
\put(17,2){3}
\put(18,3){2}
\put(19,4){4}
\put(20,5){3}
\put(21,6){4}\put(21.3,6.3){\circle{1}}

\put(7,7){$a$}

\put(3,-1){$b$}

\put(5,-1){$c$}

\put(11,5){$d$}
\put(9,-1){$e$}
\put(16,6){$f$}
\put(19,7){$g$}

\put(21,7){$\nu^{-1} a$}

\put(10,-2){figure 4}
\end{picture}
\vspace{0.5cm}

\subsection{ Fundamental algebras}

A pattern algebra $B$ inside $A(C)$ is a convex connected full subcategory that contains exactly one point from each $\nu$-orbit. Any algebra with these properties is called a fundamental algebra.
For the classification of the selfinjective algebras we need only the pattern algebras, but the fundamental algebras are exactly those having
 a representation-finite trivial extension as we will see in secton 10. A fundamental algebra comes from an $S$-pattern iff its points lie between $S$ and $\nu S$ after
 translation with some power of $\nu$.

Given any fundamental algebra $F$   the quiver of $A(C)$ contains the disjoint union 
of the quivers of all $\nu^{i}F$ as a subquiver and the only missing arrows are  $\nu$-translates of the  arrows from an $x \in F$ to an $y \in \nu^{-1} F$. 
The next lemma says when there is such an arrow.

\begin{lemma}
 In the above situation there is an irreducible morphism $\alpha:x \rightarrow y$ in $A(C)$ iff there is a complete morphism $f:\nu y\ \rightarrow x$ in $F$.
\end{lemma}

Proof: Let $\alpha$ be irreducible. It can be prolongated by an $ f:\nu y \rightarrow x$ in $F$ to a complete morphism $\alpha f$ in $A(C)$. If $f$ is not complete in $F$
there are 
two possibilities. Assume first that there is some non-invertible $g:x \rightarrow z$ in $F$ with $gf \neq 0$. Then this can be prolongated to a complete morphism $hgf$ in $A(C)$, 
whence $\alpha$ has a proper factorization. Next, let $fh$ be a prolongation in $F$. Then there is a complete morphism $fhi:\nu  x \rightarrow x$ in $A(C)$ showing
 that $\nu (\alpha)$ 
is not irreducible.

Reversely the morphism $f$ can be prolongated by some $\alpha:x \rightarrow y$ to a complete morphism in $A(C)$. Suppose $\alpha$ admits a proper factorization $hg$ with
$g:x \rightarrow z$ and  $h:z \rightarrow y$. Then $ z \in \nu^{-1} F$ because $f$ is complete in $F$ and $f \nu (h)$ is not zero which is also impossible 
because  $f$ is complete in $F$. \hspace{0.1cm}q.e.d.\vspace{0.2cm}

 For the algebra $B$ considered in 4.5 the complete morphisms are $b^{\ast} \rightarrow c^{\ast} \rightarrow d^{\ast}$,
$a^{\ast} \rightarrow d^{\ast} \rightarrow f^{\ast} \rightarrow g^{\ast}$ and 
$e^{\ast} \rightarrow f^{\ast} \rightarrow g^{\ast}$. These give rise to the additional arrows $d^{\ast} \rightarrow \nu^{-1} b^{\ast}$, $g^{\ast} \rightarrow \nu^{-1} a^{\ast}$ and $f^{\ast}
 \rightarrow \nu^{-1} e^{\ast}$. The fundamental algebra
$$ e^{\ast} \rightarrow f^{\ast} \rightarrow g^{\ast} \rightarrow \nu^{-1} a^{\ast} \rightarrow \nu^{-1} d^{\ast} \rightarrow \nu^{-2}b^{\ast} \rightarrow \nu^{-2}c^{\ast}$$
with zero-relations  $e^{\ast} \rightarrow f^{\ast} \rightarrow g^{\ast}$ and $\nu^{-1} a^{\ast} \rightarrow \nu^{-1} d^{\ast} \rightarrow \nu^{-2}b^{\ast}$ does not come from any $S$-pattern.

 Recall that theorem 6 gives  a presentation of  $A(C)$ where all paths with the same start and end have the same image and that the only complete morphisms go from any  $x$ to $\nu x$. Thus $A(C)$ 
is uniquely determined by its quiver which  in turn is given by the quiver of a fundamental algebra $F$ and the complete morphisms of $F$.

\subsection{Configurations and tilting modules}

Again $T$ is a Dynkin-diagram and and $S$ a section in $\Z\,T$. As a consequence of our results we get a bijection between $\mathcal{C}$ and the set $\mathcal{T}(S,T)$ of 
isomorphism classes of basic tilting ( left ) modules over the path algebra $kS$. Namely given any configuration $C$ the $A(C)$-module 
$M(C)= \bigoplus_{s \in S} k(\Z\,T_{C})(\,\,,s)$ has the simply connected algebra $B=B(S,C)$ as its support and so it is a finite-dimensional left-module over its endomorphism algebra
 $kS$. Since $S$ is also a section in $\Gamma_{B}$ the module $M(C)$ satisfies the criteria of \cite[section 2]{BT} for a tilting module over $B$ and so $M(C)$ is also a basic tilting module over
$kS$ by the main result of classical tilting theory \cite{}. 

\begin{theorem}
 Using the above assumptions and notations the map  sending a configuration $C$ to the isomorphism class of $M(C)$ is a bijection $\mathcal{C} \simeq \mathcal{T}(S,T)$.
\end{theorem}

Proof: If $M(C)$ is isomorphic to $M(C')$ then both modules have the same dimension-vector $d=d(C)$ which implies $C=C'$ by proposition 2.
To show that the map is also surjective let $K$ be a basic tilting ( left ) module over $kS$ i.e. a rtilting richt-module over $A=kS^{op}$. The indecomposable injectives lie on a section  isomorphic to $S$ in 
$\Gamma_{A}$ which gives a section of the same shape in $\Gamma_{B}$ under the functor $Hom(K,\,\,)$ from $mod\, A$ to $mod\, B$ with $B=End\,K$. Thus the projectives of $B$ form an $S$-pattern with
$M(C) \simeq X$ for the corresponding configuration. \hspace{0.1cm}q.e.d.\vspace{0.2cm}

As  a corollary one gets in \cite{BLR} that a dimension vector $d$ occurs as some $d=d(S,C)$ for a configuration $C$ iff the generic representation of $S$ with dimension vector $d$ has $\mid T \mid$ pairwise 
non-isomorphic indecomposable summands. 

The reader should observe that the sets of these vectors  depend on the chosen section $S$ as do the sets of isomorphism  classes of basic tilted 
algebras $End\,M(C)$ which are all of the same type $T$. For instance in  section 6.4 on $D_{4}$ the eight fundamental algebras listed there are all tilted algebras. Up to symmetries and duality there are only two possible sections and for one of them one gets only 5 isomorphism classes and for the other 7. 
 
 However, in section 10 we will see that for each section $S$ in $\Z,T$ the algebras $End\, M(C)$ hit all reflection equivalence classes of iterated tilted algebras of type $T$.

\section{Tree class $A_{n}$}
\subsection{Pedigrees}
Remember that we are looking at contravariant representations of quivers. In the following we consider an equioriented quiver $$S= 1 \rightarrow 2 \rightarrow 3 \ldots n-1 \rightarrow n$$
which we embed  into $\Z A_{n}$  by taking the points with first coordinate $0$.

Pedigrees are certain oriented trees. A point of a tree is a leaf if it has only one neighbour. The  pedigree $T(m)$ with $m$ generations 
is defined recursively as follows. $T(0)$ is just one point $\omega$. To go from $T(m)$ to $T(m+1)$ one adds to each leaf $x$ of $T(m)$ two points $x',x''$ and 
two arrows $\beta_{x}: x' \rightarrow
 x$ and $\alpha_{x} :x \rightarrow x''$. Thus $T(m)$  has $2^{m}$ leaves and $2^{m+1}-1$ points. Usually one  draws $T(m)$ in the plane as indicated in figure 5  that shows $T(4)$. Here the $\beta$'s point downwards
and the $\alpha$'s upwards and we forget about the indices.

 We write the walks in $T(m)$ with start in $\omega$ as words in $\beta^{-1}$, $\alpha$ and the empty word $\epsilon$ starting and ending in $\omega$ and we order them lexikographically 
( reading from the right to the left ) 
by  $\beta^{-1} < \epsilon < \alpha$. 
The points of $T(m)$ are just the end points
 of these walks. In figure 5 the walk from $\omega$ to $x$ is 
$\alpha \beta^{-2} \alpha$. 

For any walk $w$ as above there is a thin indecomposable representation $U(w)$ not vanishing on $x$ iff the walk passes through $x$. This indecomposable $U(w)$ is uniquely determined 
up to isomorphism by $w$ and it satisfies all relations $\alpha \beta =0$.  For any words $v ,w$ one has an epimorphism $U(v\beta^{-1}w) \rightarrow U(w)$ and a monomorphism $U(w) \rightarrow U(v \alpha w)$ 
and so there is a non-zero morphism $U(w_{1}) \rightarrow U(w_{2})$ whenever $w_{1} \leq w_{2}$.

\setlength{\unitlength}{0.4cm}
\hspace{2cm}\begin{picture}(15,7)
\put(6,5){figure 5}
\put(8,0){\vector(4,1){4}}
\put(4,1){\vector(4,-1){4}}
\put(12,1){\vector(2,1){2}}
\put(10,2){\vector(2,-1){2}}

\put(2,2){\vector(2,-1){2}}
\put(4,1){\vector(2,1){2}}

\put(6,2){\vector(1,1){1}}
\put(1,3){\vector(1,-1){1}}
\put(5,3){\vector(1,-1){1}}
\put(2,2){\vector(1,1){1}}

\put(10,2){\vector(1,1){1}}
\put(9,3){\vector(1,-1){1}}
\put(13,3){\vector(1,-1){1}}
\put(14,2){\vector(1,1){1}}

\put(3,3){\vector(1,2){0.5}}
\put(0.5,4){\vector(1,-2){0.5}}
\put(5,3){\vector(1,2){0.5}}
\put(4.5,4){\vector(1,-2){0.5}}

\put(1,3){\vector(1,2){0.5}}
\put(2.5,4){\vector(1,-2){0.5}}
\put(6.5,4){\vector(1,-2){0.5}}
\put(7,3){\vector(1,2){0.5}}

\put(11,3){\vector(1,2){0.5}}
\put(8.5,4){\vector(1,-2){0.5}}
\put(9.6,4){x}
\put(13,3){\vector(1,2){0.5}}
\put(12.5,4){\vector(1,-2){0.5}}

\put(9,3){\vector(1,2){0.5}}
\put(10.5,4){\vector(1,-2){0.5}}
\put(14.5,4){\vector(1,-2){0.5}}
\put(15,3){\vector(1,2){0.5}}
\put(7.9,-0.5){$\omega$}
\end{picture}

\vspace{1cm}

A pedigree $L$ is a subtree of some $T(m)$ containing the root $\omega$. The corresponding pedigree 
algebra $A_{L}$ is given by the quiver $L$ and all zero-relations $\alpha \beta$ occurring in $L$. 
\begin{proposition}
  For any pedigree $L$ with $n$ points the Auslander-Reiten quiver of $A_{L}$ has a unique section $U_{1} \rightarrow U_{2} \rightarrow \ldots \rightarrow U_{n-1} 
\rightarrow U_{n}$ starting with  the projective indecomposable $A_{L}$-module $P_{\omega}$. This section consists of the indecomposables $U(w)$ defined for the walks starting in $\omega$. The map 
$$L \mapsto d(L)=(dim\,U_{1},dim\,U_{2}, \ldots \,,dim\,U_{n})$$ is a bijection between the set of all pedigrees with $n$ elements and the set of all dimension vectors of
 $S$-patterns in $\Z \,T$.

\end{proposition}

Proof:  $A_{L}$ has a simply connected component in its Auslander-Reiten quiver as any tree-algebra. We show first that a section as claimed starting in $P_{\omega}$ is unique and consists of the 
$U(w)$.
Namely such a section induces an embedding of $\Gamma_{A_{L}}$ into $\Z \,A_{n}$
mapping 
$U_{i}$ to $(0,i)$. The support of the starting function $s_{(0,1)}$ shows that the $U_{i}$ are the indecomposables not vanishing on $\omega$ and they are linearly ordered by the existence
 of a non-zero homomorphism whence $U_{n}$ is the injective to the point $\omega$. Moreover for any walk $w$ as above $U(w)$ does not vanish at $\omega$ and so it belongs to the section. Since there are $n$ walks starting 
in $\omega$ the indecomposables in the section are just the $U(w)$.

Next we  show by induction on $n$ that the section exists. 
The case $n=1$ is trivial and in the inductive step  we assume first that $L$ does not contain $\omega''$. 
Then $L'=L \setminus \{\omega\}$ is a pedigree with root $\omega'$ and so $A_{L'}$ has by induction a section $V_{1} \rightarrow V_{2} \ldots \rightarrow V_{n-1}$ starting with $P_{\omega'}$.
Now $A_{L}$ is the one-point extension of $A_{L'}$ by $P_{\omega'}$ and we obtain by knitting the section 
$$U_{1}=P_{\omega} \rightarrow \tau^{-1}V_{1} \rightarrow \tau^{-1}V_{2} \ldots \tau^{-1}V_{n-2}\rightarrow \tau^{-1}V_{n-1}$$ which ends in the simple $S_{\omega}$ which is injective.

The case where $L$ does not contain $\omega'$ is dual and now we find a section starting with the projective $S_{\omega}$.
 In the general case we look at the connected components $L_{1}$ resp. $L_{2}$ of $\omega$ in $L\setminus \{\omega'\}$ resp.
 in $L\setminus \{\omega''\}$. Given a representation $V$ of $A_{L}$ we choose in $V(\omega)$ a vector-space
supplement $H$ of the image $I$ of $V(\alpha)$. That way we obtain a 
decomposition  $V=V' \oplus V''$  where $V'$ vanishes on $L''$ and $V''$ on $L'$ and 
with $V''(\omega)=I$, $V'(\omega)=H$. Thus any indecomposable except $S_{\omega}$ has support either in  $L_{1}$ or in  $L_{2}$. 
Then  $\Gamma_{A_{L}}$ is obtained by glueing $\Gamma_{A_{L_{1}}}$ and $\Gamma_{A_{L_{2}}}$ at the common point $S_{\omega}$. In particular we obtain the wanted 
section by glueing the section of $\Gamma_{A_{L_{1}}}$ with start $P_{\omega}$ and end $S_{\omega}$ with the section 
of
$\Gamma_{A_{L_{2}}}$ with start $S_{\omega}$ and end $I(\omega)$. Thus the map $L\mapsto d(L)$ is well-defined and d(L) is the dimension-vector of an $S$-pattern.

Furthermore $d(L)$ determines $L$ as one sees again by  the above
 inductive arguments. Namely there is exactly one index $i$ with $d(i)=1$. The first part $d_{1}=(d(1),d(2),\ldots ,d(i))$ of $d(L)$ consists in the dimensions of the indecomposables
 with support in $L_{1}$ and $d(j+1)=dim V(j)$ for $j=1,2,\ldots ,i-1$ are the dimensions of the indecomposables in the section of $L_{1}\setminus \{\omega\}$ considered above. By induction 
$L_{1}\setminus \{\omega\}$ is uniquely determined and so also  $L_{1}$. Similarly one gets  $L_{2}$ and finally $L$.

The map $L \mapsto d(L)$ is also surjective: Given an $S$-pattern with dimension vector   $d=(d_{1},d_{2}, \ldots ,d_{n})$ there is 
a simply connected algebra $A$ with quiver $Q$, an  embedding of 
translation-quivers
 $\gamma:\Gamma_{A} \rightarrow \Z A_{n}$, a section $U_{1} \rightarrow U_{2} \rightarrow U_{n}$ in $\Gamma_{A}$ such that $dim\,U_{i} = d_{i}$ and $\gamma \,U(i)= (0,i)$ holds for all $i$.
 Here $U_{1}$ is projective, say $U_{1}=P_{x}$, because otherwise $\tau\,U_{1}$ has no non-zero map to an injective. Let $i$ be the greatest index such that $U_{i}$ is projective. 
Then the radical of $U_{i}$
 is local or the direct sum of two local modules and so $Q$  is a tree.

Since $\Gamma_{A}$ is a sub translation quiver of $\Z A_{n}$ the homomorphism spaces 
between indecomposables have dimension $1$ at most and so all indecomposables are thin and $Q$  does not contain a subquiver of type $D_{4}$ without relation. 
Thus the support of any indecomposable $U$
is a linear 
subquiver without relation. Any point $y$ of $Q$ is the end of a unique walk
 $w(y)$ 
starting in $x$ and any arrow $\delta$ occurs inside such a walk. We choose an orientation by calling $\delta$ an $\alpha$-arrow if it points away from $x$ and a $\beta$-arrow otherwise.

The section  contains all the indecomposables not vanishing on $x$ and these are linearly ordered by the existence of a non-zero homomorphism. Two $\alpha$-arrows from $x$ to $x''_{1}$ and 
$x''_{2}$  would give two incomparable indecomposables with supports $x,x''_{i}$. Dually at most one $\beta$-arrow ends at $x$. For $\alpha \beta \neq 0$ there is an indecomposable with support
$x',x,x''$ which is incomparable to the simple $S_{x}$. Thus $\alpha \beta =0$ holds at $x$ and so any indecomposable $U$ vanishes at $x'$ or at $x''$ by the argument given above. It follows that 
the indeomposables in the section are just the thin modules $U(w(y))$ and  no zero relation occurs on such a path $w(y)$. On the other hand the argument given above for the point $x$ applies
 to any point $y$ 
and it shows that at most one $\alpha$ starts at $y$, one $\beta$ ends at $y$ and the relation $\alpha \beta =0$ holds.
Thus there is a unique isomorphism of $Q$ with a pedigree $L$ mapping $x$ to $\omega$ and $A$ is isomorphic to $A_{L}$.\hspace{0.1cm}q.e.d.
\vspace{0.2cm}

As seen in the proof above there is at most one equioriented section going 'upstairs', but there might be also another one going 'downstairs'. 
 This happens iff $A$ is hereditary iff $A$ has a sincere indecomposable which is then the intersection of the two sections.

\subsection{The simply connected categories of tree class$ A_{n}$}
Any  $A(C)$ comes from a pattern algebra $A_{L}$ and so it only remains to  find the complete morphisms in $A_{L}$ which is easy. An $\alpha$-
path 
 of length $r\geq 0$ is a path $p$ in $L$ of length $r$ containing only $\alpha$-arrows. Such a path $p$ is an $\alpha$-string if $r\geq 1$ and $p$ is not a proper subpath of another
 $\alpha$-path in $L$. Similarly $\beta$-strings are defined and it is clear that these strings give the complete morphisms.

Later on we   have to know the last two modules in the section $$U_{1} \rightarrow U_{2} \ldots U_{n-1} \rightarrow U_{n}$$ of a pedigree algebra with $n\geq 3$ points
and so  we fix some notation. Let $\omega= a_{0} \rightarrow a_{1} \rightarrow \ldots \rightarrow a_{r}=a$ be the longest $\alpha$-path starting in $\omega$. If there is a $\beta$-arrow
$b \rightarrow a$ let 
$b_{s}  \rightarrow b_{s-1} \ldots  b_{0}=b$ be the longest $\beta$-path ending in $b$  and let $b=c_{0} \rightarrow c_{1}\rightarrow \ldots c_{t}$ be the longest
$\alpha$-path starting in $b$. Here $r,s$ or $t$  might be $0$. The pedigree $L(r,s,t)$ shown in figure 6 consists only of the $r+s+t+2$ points lying on these three paths.  The next lemma is an immediate consequence of the lexikographic order defined on the walks in $L$ starting
in $\omega$.

\begin{lemma} In the above situation we always have $U_{n}=U(\alpha^{r})$. If there is a $\beta$-arrow ending in $a$  we have 
 $U_{n-1}=U(\alpha^{t}\beta^{-1}\alpha^{r})$  and $U_{n}$ is not projective. Otherwise $U_{n}$ is projective and we have 
$r\geq 1$ and $U_{n-1}=U(\alpha^{r-1})$.
 
\end{lemma}

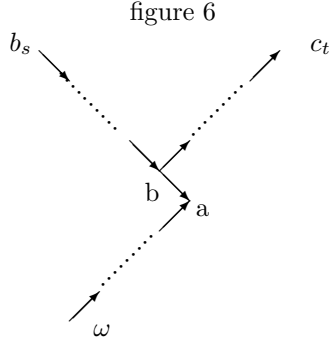
\setlength{\unitlength}{0.4cm}
\hspace{2cm}\begin{picture}(14,12)

\put(4,1){\vector(1,1){1}}
\multiput(5,2)(0.2,0.2){9}{$\cdot$}
\multiput(8,7)(0.2,0.2){9}{$\cdot$}
\multiput(3.6,9)(0.2,-0.2){10}{$\cdot$}
\put(7,4){\vector(1,1){1}}
\put(7,6){\vector(1,1){1}}

\put(10,9){\vector(1,1){1}}
\put(7,6){\vector(1,-1){1}}
\put(6,7){\vector(1,-1){1}}
\put(3,10){\vector(1,-1){1}}
\put(6,11){figure 6}
\put(4.8,0.5){$\omega$}
\put(8.2,4.5){a}
\put(6.5,5){b}
\put(2,10){$b_{s}$}
\put(12,10){$c_{t}$}
\end{picture}

\section{Tree class $D_{n}$ }
\subsection{One point extensions again}

We have already used one point extensions to get the configuration to a pattern and to characterize pedigree algebras. Now they help us to  construct the patterns resp.configurations 
for a Dynkin
 tree $T$ out of those for one or several smaller trees.

So let $S$ be a section inside $\Z\,T$  with  a sink $x$ and let $d \in \N^{S_{0}}$ be a vector.   Denote by $x_{1},x_{2},\ldots x_{r}$ the neighbors of $x$ in $S$,
  by $S_{i}$ the connected components of $x_{i}$ in $S \setminus \{x\}$ and by 
$d_{i}$ the restriction
of $d$ to $S_{i}$. Theorem 7 says that $d=d(S,P)$ for an $S$-pattern $P$ with pattern algebra $B$ iff $d=d(S,A(C))$ for a configuration $C$. The key observation is that then 
the indecomposable $B$-module $U(x)$ is projective iff $d(x)= 1 + \sum_{i=1}^{r} d(x_{i})$ iff $\tau x$ belongs to $C$. In  that case we have for each $i$ that $d_{i}=d(S_{i},P_{i})$ where 
$P_{i}$  is some
$S_{i}$-pattern defining the algebra $B_{i}$. Furthermore $B$ is then the one-point extension of the product of these $B_{i}$ by the direct sum of the $B_{i}$- modules $U(x_{i})$.
Reversely, given for each $i$ a vector $d_{i}=d(S_{i},P_{i})$ belonging to an $S_{i}$-pattern we add to  the disjoint union of the $S_{i}$ a point $x$ and arrows $x_{i}\rightarrow x$. Then 
$x$ is a sink in $S$. If  $S$ is a Dynkin-quiver
the vector $d$ defined by $ d(x)= 1 + \sum_{i=1}^{r} d(x_{i})$ and by $d(y)=d_{i}(y)$ for
 $y \in S_{i}$ belongs to an $S$-pattern. Details to all this are in   \cite[section 6]{BG}.

Thus we have an effective way to construct all configurations in $\Z T$ containing the point $x$. Obviously  $A(C)$ and $A(\tau C)$ are isomorphic for each configuration and so it suffices to 
consider only points of the form $x=(0,t)$ for some $t \in T$. For the tree class $D_{n}$ we see in the next section that one gets all 
configurations by taking only one point $t$ at the end of a short branch. This reduces the classification described in  6.3 to the case $A_{n-1}$.

The algorithm described above can be implemented easily on a computer to treat the trees $E_{6},E_{7},E_{8}$ and this was done in    \cite{BLR}.  One shows first that
 each  configuration contains a point at the end of a longest branch and so one gets even $r=1$ in the above algorithm.

\subsection{Combinatorial configurations}

 A set $E$ of  points in $\Z T$ with a Dynkin tree $T$ is called  a combinatorial configuration if it has the following two properties: 

C1: For any $x \in \Z T$ there is an $e \in E$ with $k(\Z T)(x,e) \neq 0$

C2: $k(\Z T)(e,f)=0$ holds for all $e\neq f$,$e,f \in E$. 

By 4.3 the dual of C1 holds too. A very useful fact is:
\begin{proposition}
 A configuration $C$ in $k(\Z\, T)$ is also a combinatorial configuration.
\end{proposition}

Proof: Recall the equivalence $M:k(\Z\, T_{C}) \rightarrow ind\,A(C)$ which induces an equivalence between the stable categories $k(\Z\,T) \simeq \underline{ind}\,A(C)$. Let $E$ 
be the set of points in $\Z\,T$ corresponding to simple modules.
In $\underline{ind}\,A(C)$ we have  $\underline{Hom}(M(e),M(f))=0$ for $e\neq f$ and for any $x$ there is some $e$ with $0 \neq Hom(M(x),M(e)) = \underline{Hom}(M(x),M(e))$. Thus $E$ is a
 combinatorial configuration. The configuration consists in  the first syzygies of the simples. Since taking first syzygies induces  an auto-equivalence of the
 stable category $\underline{mod} \,A(C)$ the claim follows.
 \hspace{0.1cm}q.e.d.\vspace{0.2cm}

 Riedtmann considered combinatorial configurations in \cite{Rie2} and classified them for type  $A$ in \cite{Rie2} and for $D$ in \cite{Rie3}. 
They coincide with the configurations. The same holds for type $E$ as verified by computer in \cite{BLR}. A theoretical argument argument for this seems still to be missing.

The tree $D_{4}$ is treated in  6.5. Here we consider $D_{n}$ with $n\geq 5$. We take a section $S$ in $\Z \,D_{n}$ such that the endpoint of the  long branch is the only source
and we embed  it as  $(0,1),(0,2), \ldots (0,n-2),(0,n-1),(0,n)$ with $(0,n-2)$ as the branching point.
Then  $\Z\,D_{n}$ has the shape shown in \cite{BLR}. The points $(i,x)$ of $\Z D_{n}$ with $x\geq n-1$ are called high, the others low. 
\begin{lemma}
In the situation above  take $F=\{(i,x) \mid 0 \leq i < 2n-3 \}$ as a  fundamental domain for the action of $\langle  \tau^{2n-3} \rangle$.
 For any configuration $C$ the number $h$ of high points in $C \cap F$ is $2$ or $3$.
If $h=2$ both high points have the same first coordinate. Otherwise the first coordinates are pairwise different. 
\end{lemma}
Proof: Let $A_{i}$ be the support of the point $(i,1)$ in $k(\Z\,D_{n})$. Then $A_{i}$ consists in the points $(i,j)$, $ 1 \leq j \leq n$ and $(i+1,n-2),(i+2,n-3), \ldots ,(i+n-2,1)$ as
 shown in figure  7.
By C1 any $A_{i}$ contains at least one point in $C$ and by C2 at most one low point. Moreover a high point belongs to exactly one $A_{i}$ a low point to two, 
once 'upstairs' and once 'downstairs'.
Now we are looking at the quotient map $q: \Z T \rightarrow \Z T/\langle \tau^{2n-3} \rangle$, $(i,x) \mapsto (\overline{i},x)$. Let $a$ resp. b resp c be the number of 
$i$ with $0 \leq i <2n-3$  such that 
$q(A_{i})$ contains two high points resp. one high point resp. no high point of $q(C)$. We get $a+b+c=2n-3$ and $2a + b + \frac{c}{2}=|q(C)|=n$. So we  find $3a+b=3$ 
which leads to $a=1,b=0$ or to $a=0,b=3$. \hspace{0.1cm}q.e.d.

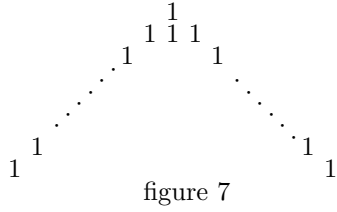
\setlength{\unitlength}{0.3cm}
\hspace{2cm}\begin{picture}(14,8)
\put(6,-1){figure 7}
\put(0,0){1}
\put(1,1){1}

\multiput(2,2)(0.5,0.5){6}{$\cdot$}
 
\multiput(10,4)(0.5,-0.5){6}{$\cdot$}
\put(5,5){1}
\put(6,6){1}
\put(7,7){1}
\put(7,6){1}

\put(8,6){1}
\put(9,5){1}

\put(13,1){1}
\put(14,0){1}

\end{picture}
\vspace{1cm}

\subsection{Two- and three-cornered categories}
The configuration $C$ and the category $A(C)$ are called two-cornered resp. three-cornered if $h=2$  resp. $h=3$. 

By 4.6 any category $A(C)$ is uniquely determined by a fundamental algebra $B$ and its complete morphisms. Such a $B$ is obtained as a one-point extension from an arbitrary pedigree $L$
with algebra $A_{L}$. Thus the quiver of $B$ contains one more point and some more arrows and also the complete morphisms change, but all this concerns only the part of $L$ 
shown in figure 6 of 5.2. We use the notation introduced there in the following.

For $h=2$ we get $B$ by iterating two one-point extensions which we combine to a 'two-point' extension
of a pedigree algebra $A_{L}$ with $n-2$ points by $U(\alpha^{r})$.
Thus the quiver of $B$ consists  of $L$, two additional points $x_{1},x_{2}$ and two additional arrows $\gamma_{1}:a \rightarrow x_{1}$ and 
$\gamma_{2}:a \rightarrow x_{2}$. The complete morphisms correspond to the $\alpha$- and $\beta$-strings as before, but $\gamma_{1}\alpha^{r}$ and $\gamma_{2}\alpha^{r}$ are two
 new complete morphisms
that replace for $r\geq 1$ the string $\alpha^{r}$. 

Remember that the pedigree algebra $A_{L}$ is uniquely determined by the dimensions of the modules on the 'subsection' $S':=(0,1) \rightarrow (0,2) \rightarrow \ldots (0,n-2)$ of $S$.
Since we have $d(n-1)=d(n)=d(n-2)+1$ we obtain:

\begin{proposition}
 The map $d \mapsto d|_{S'}$ induces a bijection between the isomorphism classes of two-cornered algebras of tree class $D_{n}$ and the set of pedigrees with $n-2$ points.
\end{proposition}

The three-cornered configurations are more complicated since one has three choices for a high point up to $\nu$-translates and also several cases for 
$r,s,t$ have to be distinguished. 

The algebra $B$  is now a one-point extension of 
a pedigree algebra $A_{L}$ with $n-1$ points by the module $U(n-2)$ occurring in the section $P(\omega)=U(1) \rightarrow U(2) \ldots \rightarrow U(n-2) \rightarrow U(n-1)$ of $\Gamma_{A_{L}}$. 
A point $p$ in $A(C)$ is a  $D_{4}$-point if it has   three neighbors $q_{1},q_{2},q_{3}$ such that the full subcategory supported by these four points
 is the path algebra of a $D_{4}$-quiver.  All statements in the following result are easy verifications. Observe that $r+s+t >0$ since $n \geq 5$ and so there is always a $D_{4}$-point.

\begin{proposition}
In the situation above we have:
\begin{enumerate} 
 \item  In $L$ there is by 5.2 a $\beta$-arrow $b \rightarrow a$ that we denote now by $\gamma_{2}$. Thus we have   $U(n-2)=U(\alpha^{t}\beta^{-1}\alpha^{r})$.
\item In the quiver of $B$ there is a new point $x$ and an arrow $a \rightarrow x$ denoted by $\gamma_{3}$.
\item  The morphism defined by $\gamma_{3}\gamma_{2}$ is complete in $B$ whence there is an arrow $\gamma_{1}:x \rightarrow \nu b$ in $A(C)$.
\item We have  $r >0$ iff there is in $A(C)$ an arrow ending in $a=a_{r}$ different from $\gamma_{2}$ iff there is in $A(C)$ an arrow different from $\gamma_{1}$ starting in $x$. 
Then $a$ and $x$ are $D_{4}$-points and the complete morphism $\alpha^{r}$ of $A_{L}$ is replaced by $\gamma_{3}\alpha^{r}$.
\item We have  $s >0$ iff there is in $A(C)$ an arrow ending in $b=b_{0}$ different from $\gamma_{1}$ iff there is in $A(C)$ an arrow different from $\gamma_{3}$ starting in $a$. Then $x$ and $b$ 
are $D_{4}$-points. For $s=0$ the arrow $\beta$ defines  no longer a complete morphism in $A(C)$.
\item We have  $t >0$ iff there is in $A(C)$ an arrow $\delta:c_{s} \rightarrow x$ ending in $x$ different from $\gamma_{3}$ iff there is in $A(C)$ an arrow different from $\gamma_{2}$
 starting in $b=c_{0}$. Then $a$ and $b$ are $D_{4}$-points.
\item Only $a,b,x$ can be $D_{4}$-points. There are two such points iff two of $r,s,t$ are $0$ and three otherwise.
\item Apart from the few cases mentioned before the complete morphisms in $B$ are just the $\alpha$- and $\beta$-strings in $L$.

\end{enumerate}

\end{proposition}

\subsection{Tree class $D_{4}$}
We discuss  this easy case with some details.
We denote the point of order $3$ by  $a$ and the others by $b,c,d$. 
All path algebras are pattern algebras for all sections of $\Z\,D_{4}$. So we apply knitting and knotting  to the quivers first with $a$ as the only source and then with
 $c$ as
 the only sink. We stop the procedure as soon as it becomes periodic.  Figure 8 shows with the encircled points the two non-isomorphic configurations obtained that way. A short calculation 
with starting functions shows that 
there are no further combinatorial configurations up to isomorphism.

\vspace{0.5cm}
\setlength{\unitlength}{0.4cm}
\begin{picture}(32,5)

\put(2,1){1}
\put(3,2){2}
\put(3,1){2}
\put(3,0){2}

\put(4,1){5}

\put(5,2){3}
\put(5,1){3}
\put(5,0){3}

\put(6,1){4}\put(6.3,1.3){\circle{1}}
\put(7,2){1}
\put(7,1){1}
\put(7,0){1}

\put(8,1){4}
\put(9,2){3}
\put(9,1){3}
\put(9,0){3}

\put(10,1){5}
\put(11,2){2}\put(11.3,2.3){\circle{1}}
\put(11,1){2}\put(11.3,1.3){\circle{1}}
\put(11,0){2}\put(11.3,0.3){\circle{1}}

\put(12,1){1}
\put(13,2){2}
\put(13,1){2}
\put(13,0){2}

\multiput(14,0)(0,0.3){10}{\line(0,1){0.2}}

\put(15,2){1}

\put(15,0){1}

\put(16,1){3}
\put(17,2){2}
\put(17,1){4}
\put(17,0){2}

\put(18,1){5}
\put(19,2){3}
\put(19,1){1}
\put(19,0){3}

\put(20,1){2}
\put(21,2){3}\put(21.3,2.3){\circle{1}}
\put(21,1){1}
\put(21,0){3}\put(21.3,0.3){\circle{1}}

\put(22,1){5}
\put(23,2){2}
\put(23,1){4}\put(23.3,1.3){\circle{1}}
\put(23,0){2}

\put(24,1){3}
\put(25,2){1}
\put(25,1){4}\put(25.3,1.3){\circle{1}}
\put(25,0){1}

\put(26,1){3}

\put(27,1){4}

\put(12,-1.5){figure 8}

\end{picture}

\vspace{1.5cm}
Figure  9  shows on the left the part of the infinite categories  $A(C)$ that is obtained by two neighbored copies of the fundamental algebras we started with and on the right the four fundamental
algebras occurring in $A(C)$ up to automorphisms of $A(C)$. The first four algebras have $rad^{2}=0$ the second four $rad^{3}=0$ and the commutativity relation.

\vspace{1cm}
\setlength{\unitlength}{0.4cm}
\hspace{1cm}\begin{picture}(20,6)
\put(0,5){\vector(1,1){1}}
\put(0,5){\vector(1,-1){1}}
\put(0,5){\vector(1,0){1}} 
\put(1,6){\vector(1,-1){1}}
\put(1,5){\vector(1,0){1}}
\put(1,4){\vector(1,1){1}}

\put(2,5){\vector(1,1){1}}
\put(2,5){\vector(1,-1){1}}
\put(2,5){\vector(1,0){1}}

\put(6,5){\vector(1,1){1}}
\put(6,5){\vector(1,-1){1}}
\put(6,5){\vector(1,0){1}}
\put(9,6){\vector(1,-1){1}}
\put(9,5){\vector(1,0){1}}
\put(9,4){\vector(1,1){1}}

\put(13,6){\vector(1,-1){1}}
\put(14,5){\vector(1,0){1}}
\put(13,4){\vector(1,1){1}}

\put(18,5){\vector(1,1){1}}
\put(18,5){\vector(1,0){1}}
\put(17,4){\vector(1,1){1}}

\end{picture}

\setlength{\unitlength}{0.4cm}
\begin{picture}(20,5)

\put(1,2){\vector(1,-1){1}}
\put(1,0){\vector(1,1){1}}
\put(2,1){\vector(1,0){1}}

\put(3,1){\vector(1,1){1}}
\put(3,1){\vector(1,-1){1}}
\put(4,2){\vector(1,-1){1}}
\put(4,0){\vector(1,1){1}}
\put(5,1){\vector(1,0){1}}

\put(7,1){\vector(1,1){1}}
\put(7,1){\vector(1,-1){1}}
\put(8,2){\vector(1,-1){1}}
\put(8,0){\vector(1,1){1}}

\put(10,2){\vector(1,-1){1}}
\put(10,0){\vector(1,1){1}}
\put(11,1){\vector(1,0){1}}

\put(14,1){\vector(1,1){1}}
\put(13,1){\vector(1,0){1}}
\put(12,0){\vector(1,1){1}}

\put(17,1){\vector(1,-1){1}}
\put(17,1){\vector(1,1){1}}
\put(16,1){\vector(1,0){1}}
\put(11,-1.5){figure 9}

\end{picture}

\vspace{2cm}

Taking into account the action of the automorphism group we obtain $5$ different configurations in the first case and $15$ in the second case. Thus here are $20$ configurations in $\Z\,D_{4}$.

\section{Fundamental groups}

To go on from the simply connected algebras $A(C)$ and their Auslander-Reiten quivers $\Z\,T_{C}$ we have to know the admissible groups $\Pi$ we can divide out i.e. 
the possible fundamental groups.

In the following tabular we have listed in the first two columns the tree  $T$ and the automorphism group $Aut (\Z \,T)$. The third column contains 
the possible groups  of automorphisms $f$ of $\Z\,T$ with $fC=C$. This  
 depends of course on $C$ and we have listed  the largest possible group which is not always an admissible group. The fourth column shows all occurring fundamental groups $\Pi$
and the last $L=L(T)$. The period $e$ of $C$ is the smallest natural number with $\tau^{e}C=C$ and  so $e$ divides $L(T)$. 
The fundamental group induces an admissible automorphism group of $\Z \,T$ whence it
 is generated by $\tau^{r}\alpha$ for some possibly trivial $\alpha$. For $r\geq L(T)$ the fundamental group is called small and large otherwise. 
As indicated by the double line large fundamental groups are only possible  for the first three trees where they really occur as will be seen in section  11. Furthermore we divide 
the tree class $D_{n}$ into the cases $h=2$ or $h=3$.
 The letters $\phi, \psi ,\chi $ denote the obvious automorphisms of $T$ of order 
$2$, $\rho$ is a glide reflection generating  $Aut (\Z A_{2n})$, $\sigma$ the element of order $3$ in $S_{3}$ and $s \geq 1$ a natural number.
\vspace{0.5cm}
 
\begin{tabular}{c|c|c|c|c}
Tree-class&$Aut(\Z T)$&$Aut(C) $&$\Pi $&$L$\\ \hline \hline
$A_{2n},n\geq 1$&$\langle \rho \rangle$&  $\langle \tau^{e} \rangle$  &$\langle \tau^{se} \rangle$  &$2n$\\ \hline
 $A_{2n+1},n\geq 1$&$\langle \tau \rangle \times \langle \phi \rangle$ &$\langle \tau^{e} \rangle,\, \langle \tau^{L} \rangle \times \langle \phi \rangle$&$\langle \tau^{se} \rangle, 
\langle \tau^{sL} \phi\rangle$&$2n+1$\\ \hline 

$D_{3m},m\geq 2$&$\langle \tau \rangle \times \langle \phi \rangle $&$h=3$,$\langle \tau^{2m-1} \rangle$&$\langle \tau^{s(2m-1)} \rangle$&$6m-3$\\ \hline\hline
$D_{4}$&$\langle \tau \rangle \times S_{3}$&$\langle \tau^{5} \rangle \times S_{3}$&$\langle \tau^{5s}\psi \rangle,\, \langle \tau^{5s}\sigma \,\rangle \,\langle \tau^{5s}\rangle$&$5$\\ \hline
$D_{n},n\geq 5$&$\langle \tau \rangle \times \langle \psi \rangle$&$h=2$,$\langle \tau^{L} \rangle \times \langle \psi \rangle $&$\langle \tau^{sL} \rangle, \, \langle \tau^{sL}\psi \rangle$ &$2n-3$\\ \hline
$D_{n},n \geq 5, 3\not| n, $&$\langle \tau \rangle \times \langle \psi \rangle$&$h=3,\langle \tau^{L} \rangle$&$\langle \tau^{sL} \rangle$&$2n-3$\\ \hline
$E_{6}$&$\langle \tau \rangle \times \langle \chi \rangle$&$\langle \tau^{L} \rangle \times \langle \chi \rangle$&$\langle \tau^{sL} \rangle,\,\langle \tau^{sL}\chi \rangle$&$11$\\ \hline
$E_{7}$&$\langle \tau \rangle$&$\langle \tau^{L} \rangle$&$\langle \tau^{sL} \rangle$&$17$\\ \hline
$E_{8}$&$\langle \tau \rangle$&$\langle \tau^{L} \rangle$&$\langle \tau^{sL} \rangle$&$29$\\ \hline
\end{tabular}

\begin{proposition}
 The data given in the table are correct.
\end{proposition}

Proof: Any automorphism $f$ commutes with $\tau$ whence $\langle \tau \rangle$ is a central subgroup of $Aut \Z T$. Moreover $f$ induces an automorphism of the graph $\overline{T}$ underlying $T$ 
whose automorphism group $Aut(\overline{T}$) is trivially isomorphic to $\Z / 2\Z$ or to the symmetric group $S_{3}$. In all cases except $A_{2n}$ the group $Aut(\overline{T})$ fixes a point $x$ in $T$ and  it
 can be lifted to a subgroup of $Aut(\Z T)$  fixing all  points $(i,x)$ which we denote again by $Aut(\overline{T})$.  Therefore  $Aut(\Z T) \simeq \langle \tau \rangle \times Aut(\overline{T})$ holds. The automorphism group of  $ A_{2n}$ is generated by the obvious 
glide reflection $\rho$ 
with $\rho^{2}=\tau$. 

Thus the first column is correct and we consider the first row. Suppose $AutC$ contains some $\rho'=\rho^{2r+1}$. Then $e=2r+1$ since $\rho'^{2}=\tau^{2r+1}$. We have $ef=2n$ and this is the 
cardinality of $C/\langle \tau^{L(T)}  \rangle$ . Any orbit has $2f$ elements and so with $m$ denoting the number of orbits $m$ we find $ef=2n=2fm$ and the contradiction  $e=2m$. Thus $AutC$
 contains only translations for type $A_{2n}$.

In all remaining cases $L$ and $e$ are odd.
Suppose now there is an element $\tau^{r}\alpha$ in $Aut(C)$  with $\alpha$ of order $2$ and  $0 \leq r <e$. Then $e$ divides $2r$ and for $r>0$ we
 obtain the contradiction $e=2r$. 
Thus $r=0$ and  $\alpha$ stabilizes $C$ and it has a
 fixed point $x$ in $C$.  

For type $A_{2n+1}$ the fixed point  $x$ is $(i,n+1)$ for some $i$. The union of the supports of 
$k(\Z A_{2n+1})(x, \,)$ and  $k(\Z A_{2n+1})( \,,x)$ cannot contain another point $(j,n+1)$ since $C$ is a combinatorial configuration.  Thus we get $e=2n+1$.
The type $D_{4}$ is clear by section 6.4. In the case $D_{n}$ with $n\geq 5$ we get $e=L$ for $h=2$ and $3e=L$ for $h=3$. Finally in the exceptional cases $E_{6},E_{7},E_{8}$ the period $e$ divides 
the prime $L$ which is stricly greater than the cardinality of $T$. Thus for $e=1$ the configuration  contains too many points. \hspace{0.1cm}q.e.d.
\vspace{0.2cm}

The  Auslander-Reiten quivers of all connected locally representation-finite selfinjective categories are just the quotients $\Z\,T_{C}/\Pi$
 and  
$\Z\,T'_{C'}/\Pi '$ is an equivalent category iff $T=T'$ and there is an $f \in Aut(\Z\,T)$ with $fC=C'$ and $f \Pi f^{-1}=\Pi '$. 

Similarly the ordinary quivers are just the quotients of the quiver of $A(C)$ by $\Pi$ but  we do not yet know the relations of the standard category and whether there are non-isomorphic
 algebras
 with the same Auslander-Reiten quiver. These questions are completely answered in \cite{Rie2,BLR,Rie3} following the first approach but we concentrate now more on the algebras, 
their quivers with relations and their coverings. This leads to simpler proofs and it applies also to certain representation-infinite algebras.

\section{Multiplicative bases in distributive categories}
\subsection{Distributive categories and their stem and ray categories}

First we give some definitions from \cite{BGRS}.
A base category $B$ is a category with uniquely determined zero morphisms in each $B(x,y)$ such that different objects are not isomorphic, all endomorphisms different 
from the identity are nilpotent and for each $x$ there are only finitely many morphisms different from zero starting or ending in $x$. For any field $k$ the linearization $kB$ is the locally bounded
 $k$-linear category with the
 same objects as $B$ and $kB(x,y)$ is the quotient of the $k$-vector space with basis $B(x,y)$ divided by the subspace spanned by the zero morphism. The composition is 
induced by the composition in $B$.

Any base category $B$ is given by a quiver $Q_{B}$ and relations. 
 The points of $Q_{B}$ are the objects of $B$ and the arrows the irreducible morphisms in $B$, i.e. the non-zero 
non-invertible morphisms that are not a composition of non-invertible morphisms. The  ( non-linear ) path category $PQ_{B}$ has as morphisms the paths and 
additional zero-morphisms. The composition is concatenation and there is  a 
full functor $p:PQ_{B} \rightarrow B$ which is the identity on the objects and the arrows. Thus each base category can be
described by a quiver and a stable equivalence relation on its paths i.e. an equivalence relation which is invariant under left and right multiplication. Of course $Q_{B}$ is the quiver
 of $kB$ for any field and the kernel of $kp:kQ_{B} \rightarrow kB$ is generated by paths and differences of paths.

A $k$-category $A$ is distributive if it is locally bounded and if the lattice of two-sided ideals is distributive or -equivalently by Kupischs observation in \cite{Kupisch}- if  $A(x,x)\simeq k[T]/(T^{n_{x}}$ is uniserial for any $x$
 and 
 $A(x,y)$ is always cyclic  as a right $A(x,x)$-module  or a left $A(y,y)$-module. In case  $A(x,y)$ is always cyclic from both sides $A$ is called regular. By an old result of Jans in \cite{Jans}
 any locally representation-finite category is distributive.

 For subsets $J$ of
$A(x,y)$ and $J'$ of $A(y,z)$ we write $J'J$ for
the set of all compositions and  $J'\cdot
J$ for the subbimodule generated by $J'J$. Furthermore we abbreviate by $N(x)$ the radical of the algebra $A(x,x)$.

To any distributive category $A$ one  associates two base categories: the stem category $\hat{A}$ and the ray category $\vec{A}$.   Both categories have the same objects as $A$ 
and the  subbimodules of $A(x,y)$ form the morphism set $\hat{A}(x,y)= \vec{A}(x,y)$ but the compositions differ slightly. In $\hat{A}$ the composition of two morphisms $J'$ and $J$ 
 is just the product $J'\cdot J$ which is clearly associative.
In $\vec{A}$ we denote the composition by $J'\ast J$. Here again 
$J'\ast J = J' \cdot J$ unless  $J'\cdot N(y) \cdot J = J'\cdot J \neq 0$ in which case one defines   $J' \ast J =0$. Note that the exceptional case never occurs in a regular category and so
 we have $\hat{A}=\vec{A}$ in that  case.

 Somewhat surprisingly the $\ast$-product is  associative whence $\vec{A}$ is a category indeed:

\begin{lemma}
 let $A$ be distributive with objects $x,y,z,t$ and subbimodules $J_{1} \subseteq A(x,y),J_{2} \subseteq A(y,z),J_{3} \subseteq A(z,t)$. Then we have 
$J_{3}\ast (J_{2} \ast J_{1})= (J_{3}\ast J_{2}) \ast J_{1}.$ Furthermore $J_{3}\ast J_{2} \ast J_{1}= J_{3} \ast J_{1} \neq 0$ implies $J_{2}=A(y,y)$.
\end{lemma}
Proof: For $J' \ast J \neq 0$ we have $J' \ast J = J' \cdot J$. Thus $\ast$ is associative if both sides are not $0$. 

Up to duality we have to exclude 
 $J_{3}\ast (J_{2} \ast J_{1})= 0 \neq (J_{3}\ast J_{2}) \ast J_{1} = J_{3}\cdot J_{2}\cdot J_{1}$. Suppose first $J_{2}\ast J_{1}=0$ i.e. $J_{2}\cdot J_{1}= J_{2} \cdot N(y) \cdot J_{1}$.
Then we get  
$(J_{3}\ast J_{2}) \cdot J_{1}= J_{3}\cdot J_{2} \cdot J_{1}= J_{3}\cdot J_{2} \cdot N(y) \cdot J_{1} = (J_{3}\ast J_{2}) \cdot N(y)\cdot J_{1}$
whence the contradiction $ (J_{3}\ast J_{2}) \ast J_{1} = 0$. 

Thus we have $J_{2}\ast J_{1}\neq 0$ but $J_{3}\cdot N(z)\cdot  J_{2}\cdot J_{1}= J_{3}\cdot J_{2}\cdot J_{1}$. From
$J_{3}\ast J_{2}\neq 0$ we see $J_{3}\cdot N(z)\cdot J_{2} \neq J_{3}\cdot J_{2}$.

 If $A(y,t)$ is cyclic over $A(t,t)$  we have 
$J_{3}\cdot N(z) \cdot J_{2} \subseteq N(t) \cdot J_{3}\cdot J_{2}$ leading to the  contradiction $N(t) \cdot J_{3}\cdot J_{2} \cdot J_{1}= J_{3}\cdot J_{2}\cdot J_{1}\neq 0$.

If $A(y,t)$ is cyclic over  $A(y,y)$ we have 
$J_{3}\cdot N(z) \cdot J_{2} \subseteq  J_{3}\cdot J_{2} \cdot N(y)$ leading to $ (J_{3}\ast J_{2}) \cdot J_{1}= (J_{3}\ast J_{2}) \cdot N(y) \cdot J_{1}$
whence to $(J_{3}\ast J_{2}) \ast J_{1}=0$. We have shown that $\ast$ is associative.

The last statement of the lemma is obvious. \hspace{0.1cm}q.e.d.\vspace{0.2cm}

\subsection{Multiplicative bases}

 For a set $M$ of morphisms in a distributive category $A$ we consider the following axioms \cite[section 1.1]{BGRS}:

M1: $M$ contains all identities and for each pair $x,y$ of points  $M \cap \mathcal{R}A(x,y)$ is a basis of $\mathcal{R}A(x,y)$.

M2: For $m \in M\cap A(x,y)$ and $n \in M\cap A(y,z)$ we have $mn=0$ or $mn \in M$.

M3: For $l,m,n$ in $M$ the relation  $ ln = lmn \neq 0$ implies that $m$ is an identity.

A subset $M$ satisfying $M1$ and $M2$ is called a multiplicative basis and a cancellative  basis if it satifies in additon $M3$. 
If $M$ is  a  multiplicative basis $M \cap \mathcal{R}^{n}A(x,y)$ is a basis of $\mathcal{R}^{n}A(x,y)$ whenever this is not $0$. Therefore each subbimodule $J\neq 0$ of some $A(x,y)$ 
is generated as a bimodule by exactly one $m=m_{J} \in M$ and so the set $\mathcal{M}$ of non-zero bimodules parametrizes any 
multiplicative basis. Finally one should observe that a  multiplicative basis is uniquely determined by the $m_{J}$ corresponding to the bimodules 
$J=\mathcal{R}A(x,y) \neq \mathcal{R}^{2}A(x,y)$ that define the arrows of the quiver.

The next lemma narrows down the possibilities for the product of two elements in a multiplicative basis.

\begin{lemma}
 Let $M$ be a multiplicative basis in the distributive category  $A$ and let $0\neq J \subseteq A(x,y)$,$0\neq J' \subseteq A(y,z)$ be two subbimodules. Then we have:
\begin{enumerate}
 \item For  $J'\cdot N(y) \cdot J \neq J' \cdot J $ the product $yx$ of any two generators $x$ of $J$ and $y$ of $J'$ is a generator of $J'\cdot J$. In particular
$m_{J'}m_{J}=m_{J'\cdot J}$ holds..
\item $M$ is  cancellative iff $m_{J'}m_{J}=0$ whenever  $J'\cdot N(y)\cdot J = J' \cdot J$. 

\end{enumerate}

\end{lemma}

Proof: $J'\cdot J$ is generated as a vector space by $yx$ and by products $rysxt$ where  at least one of $r,s,t$ lies in the radical of $A$.
 By the assumption made in part i) all these products are in the radical of the bimodule $J' \cdot J$. Thus $yx$ generates $J'\cdot J$.

For part two suppose  first that $M$ is cancellative. For $J'\cdot J =0$ one has obviously $m_{J'}m_{J}=0$. So we can assume $0 \neq J'\cdot J=J'\cdot N(y)\cdot J$. Then $A(x,y)$ and $A(y,z)$ are both 
cyclic over $A(y,y)$. Thus we have 
$m_{J'}m_{J}=m_{J'} \rho m_{J}$ for some $\rho= \sum a_{i}m_{J_{i}}$ in $N(y)$. If now $m_{J'}m_{J}\neq 0$ one uses the 
linear independence and the 
multiplicativity to get $m_{J'}m_{J}=m_{J'}m_{J_{i}}m_{J}$ for some $J_{i}$ contained in $N(y)$. The contradiction $id_{y} \in N(y)$ follows.

Reversely $0 \neq m_{J'}m_{J}=m_{J'}m_{J''}m_{J}$ implies $J'\cdot J \neq J' \cdot N(y) \cdot J$ whence $m_{J'}m_{J}=m_{J' \cdot J}$ by the first assertion of the lemma.
 Then $m_{J''}=id_{y}$ follows. \hspace{0.1cm}q.e.d.\vspace{0.2cm}

Any  multiplicative basis $M$ in a distributive category defines a  category $B(M)$. The objects are the same as in $A$, the non-zero morphisms are the elements 
in $M \cap A(x,y)$ to which one adds zero-morphisms $0_{xy}$
 and the composition is the restriction of the composition in $A$.

The next result is a consequence of the last lemma. 
\begin{lemma}
 Let $M,M'$ be two multiplicative bases in a distributive category $A$. Then $B(M)$ and $B(M')$ are equivalent provided $M$ and $M'$ are also cancellative.

\end{lemma}

  If $M$ is a multiplicative basis of the distributive category $A$ 
then $A$ and $kB(M)$ are equivalent. Furthermore $k\hat{A}$ has $\hat{A}$ as a multiplicative basis, $k\vec{A}$ has $\vec{A}$ as a cancellative basis and $\hat{A}$  is not equivalent to $\vec{A}$ 
if the exceptional case $J'\cdot N(y) \cdot J = J'\cdot J \neq 0$ occurs in $A$. Nevertheless the linearizations $k\hat{A}$ and $k\vec{A}$ might be isomorphic and this might even depend on
 the field as we will see in 9.5. 
Finally a distributive category $A$ has a cancellative basis iff $A\simeq k\vec{A}$.

Recall that in general an algebra has no filtered multiplicative basis at all and if it has it can have many different such bases even with many non-equivalent base categories.
 For instance let $A$ be given by the 
 quiver with three points $a,b,c$, three arrows $\mu:a \rightarrow b$, $\rho:b \rightarrow b$ and $\nu:b \rightarrow c$ and by the relations $\rho^{3},\nu\mu, \nu\rho\mu$.
Then we  have three different filtered multiplicative bases by taking as generators always $\rho$ and $\nu$ but as the third generator either $\mu$ or $\mu + \rho\mu$ or $\mu + \rho^{2}\mu$. 
Then the corresponding base-categories 
are pairwise not equivalent and in the last two cases there occur also commutativity relations, but the linearizations are of course all isomorphic to $A$. 

Finally note that for a distributive category $\hat{A}$ and $\vec{A}$ always determine each other.\newpage

\subsection{Regular representation-finite selfinjective algebras have a multiplicative basis}

The following  nice result of  \cite{Kupisch5} applies to representation-finite algebras. 

\begin{theorem}
 Let $A$ be a basic connected regular selfinjective algebra such that for no idempotent $e$ the quiver of $eAe$  contains a bipartite quiver of type $\tilde{A}_{2n}$ with $n\geq 1$.
Then $A$ has a  multiplicative basis.
\end{theorem}

Proof: The assumption implies that up to duality there is a point $x$ in the quiver of $A$ where only one arrow $\alpha:x \rightarrow y$ starts.  Let $\nu$ be the Nakayama automorphism  of $A$. 
Then $\nu$ acts on the quiver $Q$ of $A$ and we denote by $N$ the orbit of $x$ and by  
$x_{1}=x,x_{2}, \ldots x_{r}$ its elements.

First assume that $y \in N$. Then there is also only one arrow starting in $y$ and so on whence the points of $N$ lie on on oriented cycle $C$. An arrow $\beta:z \rightarrow x_{i}$ with $z \not\in N$ 
induces a non-zero morphism in $A(z,x_{i})$ that can be prolongated to a non-zero morphism in $A(z,\nu z)$ with $\nu z \not\in  N$. This  is impossible since there is no path leaving $C$.
Thus $Q$ is an oriented cycle and $A$ is a generalized Nakayama algebra and so it has a multiplicative basis.

We can assume that there is no arrow in $Q$ between two points of $N$. We define the idempotent $e$ as the sum of all idempotents corresponding to points outside $N$.
By induction on the dimension the algebra $eAe$ has a multiplicative basis $B'$. 

 For any $i$ there is an arrow $\alpha_{i}:x_{i} \rightarrow y_{i}$ and we choose an element 
$a_{i} \in  \mathcal{R}A(x_{i},y_{i}) \setminus \mathcal{R}^{2}A(x_{i},y_{i})$. Since there is no simple projective there is at least one arrow stopping at each point. We denote by  
$\beta_{i,1}:z_{i,1} \rightarrow x_{i},\beta_{i,2}:z_{i,2} \rightarrow x_{i},\ldots \beta_{i,s}:z_{i,s} \rightarrow x_{i}$  the arrows stopping at $x_{i}$. 
For each $i$ and $j$ we choose elements  $b_{i,j} \in \mathcal{R}A(z_{i,j},x_{i}) \setminus \mathcal{R}^{2}A(z_{i,j},x_{i})$.
The non-zero morphism $b_{i,j}$ can be prolongated to a non-zero morphism from $z_{i,j}$ to $\nu z_{i,j}$ whence  $0\neq a_{i} b_{i,j}$ holds for all $i$ and $j$. Therefore  we get
 $a_{i} b_{i,j}= b'_{i,j}u_{i.j}$. Here $b'_{i,j}$ is the unique element in $B'$ generating the same subbimodule as  
$a_{i} b_{i,j}$ and $u_{i,j}$ is a unit. Replacing each $b_{i,j}$ chosen before by $b_{i,j}u_{i,j}^{-1}$ we see that all products $a_{i} b_{i,j}$ belong to $B'$.

Let $M$ be the set consisting of all $a_{i}$, all $b_{i,j}$ and all elements in $B'\cap rad\,A$. 
We claim that the set $B$ consisting of all non-zero products of elements in $M$ and of all idempotents $e_{p}$, $p \in Q$, is a  multiplicative basis in $A$. By construction $B$
generates $A$ as a vector space and $B \cup \{0\}$ is stable under multiplication whence contained in the union of the $A(q,p)$ , $p,q \in Q$. To prove that  $B$ has $dim\,A$ elements we
show that for all $p,q \in Q$ and for all $v,w \in  A(q,p)\cap B$ always $v=w$ follows from $A_{v}=A_{w}$. Here $A_{u}$ denotes the subbimodule of $ A(q,p)$ generated by some $u$.

Four cases are possible. For $p$ and $q$ not in $N$ both elements $v$ and $w$ lie in $B'$ and so $A_{v}=A_{w}$ implies $v=w$. Next take $p=x_{i} \in N$ and $q \not\in N$. Then we have
 $v=b_{i,j}v'$ and $w=b_{i,j'}w'$ and $v=wu$ for some unit $u$ in $A(q,q)$. The non-zero morphisms $v$ and $w$ in $A(q,p)$ can be prolongated to non-zero morphisms in $A(q,\nu q)$.
Thus we have  $a_{i}v \neq 0 \neq a_{i}w$ and both elements lie in $B'$ and generate the same subbimodule of
 $A(q,y_{i})$. From the first case we get $a_{i}v=a_{i}w$ or $a_{i}(v-w)=0$. This   implies $v=w$ because otherwise the non-zero morphism $v-w$ could be prolongated non trivially to $\nu q$.

 Finally assume $q=x_{j}$ and $p$ arbitrary. Then we have $v=v'a_{j}$ and  $w=w'a_{j}$.  From $v'= \rho w'$ with some nilpotent $\rho$ we get $v=w\rho$ contradicting $A_{v}=A_{w}$. Similarily 
$w'=\rho v'$ with some nilpotent $\rho$ is impossible. 
 Thus we have $A_{v'}=A_{w'}$ and $v'=w'$ follows from the cases already considered. \hspace{0.1cm}q.e.d.

\section{Galois coverings}
\subsection{Quotients by group actions}Covering functors arise often as quotients under group actions: 
Let $C$ be a  locally finite dimensional category on which a group $G$ of automorphisms acts freely  - i.e. no group element  except $1 $ fixes an object  -
 and locally bounded - i.e. for given 
$a,b$ in $C$ there are only finitely many $g$ with $C(a,gb)\neq 0$. Then there is a
 quotient $F:C \rightarrow C/G$ called a Galois covering and this is a covering functor. Here  any  functor
$ E:C \rightarrow D$ with $Eg=g$ for all $g$ in $G$ factors uniquely as $E= E'\circ  F$ and  $E'$ is an equivalence iff E is a covering functor with the orbits as fibres.

Of course, the action of $G$ on $C$ induces an action $m \mapsto m^{g}$ on $C$-mod by scalar extension. This action is free if $m\not\simeq m^{g}$ holds for all $g\neq 1$ and 
$m\neq 0$. If $G$ is torsion-free this follows.

\begin{theorem}\label{Galois}( \cite{Galoiscovering} ) Let $C$ be a locally bounded connected category and $G$ a group of automorphisms acting freely  on $C$ and $C-mod$.
 Denote by $F$ the quotient and by
$F_{\bullet}:C-mod \rightarrow C/G-mod$ the left adjoint to the restriction.
\begin{enumerate}
 \item $F_{\bullet}$ is exact. It preserves dimensions, indecomposibility and infinite families of pairwise non-isomorphic indecomposables of the same dimension.
\item $C$ is locally representation-finite iff $C/G$ is so. In that case $F_{\bullet}$  induces a covering
$ind C \rightarrow ind (C/G)$ and a covering $\Gamma_{C} \rightarrow \Gamma_{C/G}$ between the Auslander-Reiten quivers.

\end{enumerate}
 
\end{theorem}

The situation is even better for 'combinatorial actions' i.e. if $C$ has a presentation $C=kQ/I$ by quiver and relations and  $G$ is given by a group of permutations of $Q$ 
such that the induced automorphisms of $kQ$ let $I$ invariant. Then we write $C//G$ for the quotient.
 The quiver of $C//G$ is  given by the orbits of points and arrows  and the defining ideal is the image of $I$ under the obvious functor $kQ \rightarrow k(Q/G)$. 
Furthermore $F_{\bullet}$  is on the level of representations simply defined  by 'taking direct sums of vector spaces and linear maps'.

Such a combinatorial action occurs twice for any locally representation-finite category $A$. Namely the 
 fundamental group $\Pi$ of the universal cover $\pi:\tilde{\Gamma}_A \rightarrow \Gamma_{A}$ acts on the mesh category $k(\tilde{\Gamma}_{A})$ letting the mesh-ideal invariant 
and the combinatorial quotient is  $k(\pi): k(\tilde{\Gamma}_{A}) \rightarrow k(\Gamma_{A})$.  Now $\Pi$ acts also on the the full subcategory $\tilde{A}$ of projectives in  
$k(\tilde{\Gamma}_{A})$ inducing the Galois covering $G:\tilde{A} \rightarrow \tilde{A}/\Pi =A^{s}$. Moreover $\Pi$ permutes the points of the quiver $\tilde{Q}$ of $\tilde{A}$ and also the arrows
  because $\tilde{A}$ is square-free. The kernel of any standard presentation $\phi$ as in theorem      is  invariant under these permutations and we get a combinatorial quotient
$H:\tilde{A} \rightarrow \tilde{A}//\Pi$. The $\Pi$-invariant cancellative basis of $\tilde{A}$ induces a cancellative basis on  the quotient. Observe that for any arrow $\alpha$ in $\tilde{Q}$ 
and any $g \in \Pi$ the elements $\phi (g\alpha)$ and $g\phi(\alpha)$ are proportional but in general not equal i.e. $\phi$ is not $\Pi$-invariant.

\begin{theorem}
 In the above situation there is a $\Pi$-invariant standard presentation whence the quotients $A^{s}$ and $\tilde{A}//\Pi$ are isomorphic.
\end{theorem}

 The independent proofs of this remarkable result in \cite{BrG,delapena2} are based on  the argument of Gabriel in \cite[section 5]{BG} that some second cohomology group vanishes because
 $\Pi$ 
is a free group as shown by Gabriel in  \cite[section 4]{BG}.  Theorem 11 is not used in our version of the classification. 

We illustrate the situation by an example. Let $A(x)$ be the category with point set $\Z \times \Z$ shown in figure  10. In the rows the arrows fly
from the left to the right, in the columns they fly downwards. 
The points are periodically numbered by $1,2,3$ as indicated. The composition of two arrow in a row or in a column is $0$ and all squares are commutative with the exception of those starting 
in a point numbered one. Then we  have $1 \rightarrow 2 \rightarrow 1 = x(1 \rightarrow 3 \rightarrow 1)$ for some non-zero scalar $x$. The translation by three points to th right and those by 
three 
points downwards induce a combinatorial action of $\Z \times \Z$ on $A(x)$ and one obtains a one-parameter-family of tame self-injective algebras as combinatorial quotients. Nevertheless all $A(x)$
 are isomorphic to $A(1)$, but there is no $\Z\times \Z$-equivariant equivariant isomorphism.

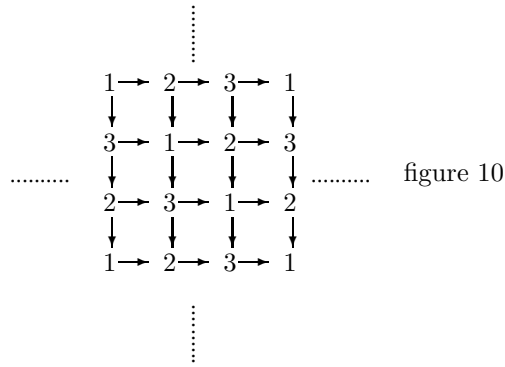
\setlength{\unitlength}{0.4cm}
\hspace{1cm}\begin{picture}(14,14)
\put(4,10){1}\put(6,10){2}\put(8,10){3}\put(10,10){1}
\put(4,8){3}\put(6,8){1}\put(8,8){2}\put(10,8){3}
\put(4,6){2}\put(6,6){3}\put(8,6){1}\put(10,6){2}
\put(4,4){1}\put(6,4){2}\put(8,4){3}\put(10,4){1}

\put(4.5,10.3){\vector(1,0){1}}\put(6.5,10.3){\vector(1,0){1}}\put(8.5,10.3){\vector(1,0){1}}

\put(4.5,8.3){\vector(1,0){1}}\put(6.5,8.3){\vector(1,0){1}}\put(8.5,8.3){\vector(1,0){1}}
\put(4.5,6.3){\vector(1,0){1}}\put(6.5,6.3){\vector(1,0){1}}\put(8.5,6.3){\vector(1,0){1}}
\put(4.5,4.3){\vector(1,0){1}}\put(6.5,4.3){\vector(1,0){1}}\put(8.5,4.3){\vector(1,0){1}}

\put(4.3,9.8){\vector(0,-1){1}}\put(6.3,9.8){\vector(0,-1){1}}\put(8.3,9.8){\vector(0,-1){1}}\put(10.3,9.8){\vector(0,-1){1}}
\put(4.3,7.8){\vector(0,-1){1}}\put(6.3,7.8){\vector(0,-1){1}}\put(8.3,7.8){\vector(0,-1){1}}\put(10.3,7.8){\vector(0,-1){1}}
\put(4.3,5.8){\vector(0,-1){1}}\put(6.3,5.8){\vector(0,-1){1}}\put(8.3,5.8){\vector(0,-1){1}}\put(10.3,5.8){\vector(0,-1){1}}

\multiput(1,7)(0.2,0){10}{\circle*{0.1}}
\multiput(11,7)(0.2,0){10}{\circle*{0.1}}
\multiput(7,11)(0,0.2){10}{\circle*{0.1}}
\multiput(7,1)(0,0.2){10}{\circle*{0.1}}
\put(14,7){figure 10}
\end{picture}
\vspace{1cm}

\subsection{Covering functors with the same fibres}
The following elementary result is apparently new. It does not depend of any assumptions on the representation type:

\begin{theorem} Let $F:\tilde{A} \rightarrow A$ and $F':\tilde{A} \rightarrow A'$ be two covering functors between  locally bounded categories. Suppose that  $\tilde{A}$ is directed and that $F$ and $F'$
 have the same fibres on the objects. 
 Then  $A$ is distributive resp. regular iff $A'$ is so.  In that case both stem categories are  equivalent and therefore also the ray categories.

\end{theorem}
 
In the proof we denote by $rad^{i}A(x,y)$ the $i^{th}$ radical of $A(x,y)$ as a bimodule.

\begin{lemma}Let $F:\tilde{A} \rightarrow A$  be a covering functor between locally bounded categories where $\tilde{A}$ is directed and $A$ is distributive. Take points $x,y \in A$ with
$dim\,A(x,y) = t+1 \geq 1$.
Then the following holds:
\begin{enumerate}

\item Suppose  $A(x,y)$ is cyclic over $A(y,y)$. For any $\tilde{x}$ with $F\tilde{x}=x$  there are points $\tilde{y}_{0},\tilde{y}_{1}, \ldots \tilde{y}_{t}$ and  morphisms
 $v_{0}:\tilde{x} \rightarrow \tilde{y}_{0},v_{i}:\tilde{y}_{i-1} \rightarrow \tilde{y}_{i}$,$1 \leq i \leq t$, such that $rad^{i}A(x,y)$ is 
generated by $\mu_{i}:=F(v_{i}v_{i-1  \ldots }v_{0})$  for each $i$. Here the points are uniquely determined by $\tilde{x}$ and they are the only points $\tilde{y}$ 
in the fibre of $y$ with $\tilde{A}(\tilde{x},\tilde{y})\neq 0$.

\item The dual also holds.
\item $\tilde{A}$ is Schurian.
\end{enumerate}

\end{lemma}

Proof: a) First we define $v_{0}$. For $x=y$ we take $\tilde{y}_{0}=\tilde{x}$ and $v_{0}= id_{\tilde{x}}$. For $x\neq y$ the isomorphism 
$\bigoplus_{F\tilde{y}=y} \tilde{A}(\tilde{x},\tilde{y}) \simeq A(x,y)$ implies the existence of a $v_{0}:\tilde{x} \rightarrow  \tilde{y}_{0}$ such that $Fv_{0}$ generates $A(x,y)$. 

Any covering functor respects the radicals and all objects of $\tilde{A}$ have one-dimensional endomorphism algebras because $\tilde{A}$ is directed.
Therefore $F$ induces for each lifting $\tilde{y}$ of $y$ an isomorphism between $\bigoplus_{F\tilde{z}=y,\tilde{z}\neq \tilde{y}} \tilde{A}(\tilde{y},\tilde{z})$
  and $R:= \,rad\,A(y,y)$. So there are morphism $v_{i}:\tilde{y}_{i-1} \rightarrow \tilde{y}_{i}$ with $\tilde{y}_{i-1} \neq \tilde{y}_{i}$ such that $Fv_{i}$ generates $R$ for $i=1,2, \ldots t$.
By construction $F(v_{t}v_{t-1} \ldots v_{0}) \neq 0$ and so $\mu_{i}$ generates  $rad^{i}A(x,y)$ and $\tilde{A}(\tilde{x},\tilde{y}_{i}) \neq 0$ for each $i$. 
Since $\tilde{A}$ is directed the $\tilde{y}_{i}$'s are pairwise different. From $dim \,A(x,y) =t+1$ we infer  that $dim \,\tilde{A}(\tilde{x},\tilde{y}_{i}) =1$ for all i and 
$\tilde{A}(\tilde{x},\tilde{y}) = 0$ if $\tilde{y}$ is not one of the $\tilde{y}_{i}$.

Part b) holds because the definition of a covering functor is self-dual and $ dim\, A(x,y) \leq 1$ is shown for $A(x,y)$ cyclic over $A(y,y)$ in the proof of part a). Thus part c) is true
 because $A(x,y)$ is always cyclic from one side.\hspace{0.1cm}q.e.d. 
 \vspace{0.2cm}

Now we prove the theorem. For any object $x \in A$ there is an object $\tilde{x} \in \tilde{A} $ with $F\tilde{x}=x$ and $x'= F'\tilde{x}$ is by assumption independent of the choice of
the lifting $\tilde{x}$. Obviously
$x \mapsto x'$ is a bijection between the sets of objects  and for any $x,y$ the vector spaces $A(x,y)$ and $A'(x',y')$ are isomorphic to $\oplus 
_{F\tilde{y}=y} \tilde{A}(\tilde{x},\tilde{y}) =\oplus _{F'\tilde{y}=y'} \tilde{A}(\tilde{x},\tilde{y})$.

Suppose  that $A$ is distributive and that $A(x,y)$ is a uniserial  $A(y,y)$-module. We show that $A'(x',y')$ is uniserial as a  $A'(y',y')$-module. This is clear for $dim\,A(x,y)=t+1 =1$. For
$t \geq 1$ we choose a point $\tilde{x}$ with $F\tilde{x}=x$ and we take a sequence $v_{0},v_{1}, \ldots v_{t}$ of morphisms as in part a) of the last lemma. Then 
$F(v_{t}v_{t-1} \ldots v_{0}) \neq 0$ implies $v_{t}v_{t-1} \ldots v_{0}\neq 0$ and so $F'(v_{t}v_{t-1} \ldots v_{0}) \neq 0$. Here $F'v_{i}$ is in $R':=rad\,A'(y',y')$ for all $i\geq 1$ 
because of $\tilde{y}_{i-1}\neq \tilde{y}_{i}$. Thus $(R')^{t}\neq 0$ and  $\mu'_{i}:=F'(v_{i}v_{i-1} \ldots v_{0})$ generates $rad^{i}A'(x',y')$ for each  $i$. 

Dually $A'(x',y')$ is uniserial from the right if  $A(x,y)$ is so. In particular $A'$ is distributive resp. regular if $A$ is. The other implication follows  symmetrically.

We compare the stem categories of $A$ and $A'$. We have a bijection $\Phi$ given by $x \mapsto x'$ on the objects and by $rad^{i}\,A(x,y) \mapsto rad^{i}\,A'(x',y')$ on the morphisms.
Given three not necessarily distinct points $x,y,z$ in $A$ we set $dim\,A(x,y)=t+1, dim\,A(y,z)=s+1,dim\,A(x,z)=r+1$. Up to duality we can assume that $A(x,z)$ is cyclic over $A(z,z)$.
To prove that $\Phi$ respects the composition  we distinguish three cases. In the first two cases the product of generators of the considered bimodules is a generator of the product
by lemma 7 in 8.2. Furthermore $rad^{j}A(y,z) \cdot rad^{i}A(x,y) \mapsto rad^{j}A'(y',z') \cdot rad^{i}A'(x',y')$ is obvious for $i=0$ or $j=0$.

First let $A(x,y)$ and $A(y,z)$ be also cyclic from the left. Then we take an  $\tilde{x}$ with $F\tilde{x}=x$, a sequence of points  $\tilde{y}_{0},\tilde{y}_{1}, \ldots \tilde{y}_{t}$ and
  morphisms
 $v_{0},v_{1}, \ldots v_{t}$. In the same vein we take for $y,z$ the lifting $\tilde{y}_{i}$, points  $\tilde{z}_{0},\tilde{z}_{1}, \ldots \tilde{z}_{s}$ and 
 morphisms
 $w_{0},w_{1}\ldots w_{s}$ with $w_{0}:\tilde{y}_{i} \rightarrow \tilde{z}_{0},w_{p}:\tilde{z}_{p-1} \rightarrow \tilde{z}_{p}.$ We set $f=w_{j}w_{j-1}\ldots w_{0}v_{i}\ldots v_{0}$.
 Then 
$\mu_{i}:=F(v_{i}v_{i-1  \ldots }v_{0})$ generates $rad^{i}A(x,y)$, $\mu_{i}':=F'(v_{i}v_{i-1  \ldots }v_{0})$ generates $rad^{i}A'(x',y')$,
$\nu_{j}:=F(w_{j}w_{j-1  \ldots }w_{0})$ generates $rad^{j}A(y,z)$, $\nu_{j}':=F'(w_{j}w_{j-1  \ldots }w_{0})$ generates $rad^{j}A'(y',z')$ and 
 $\nu_{j}\mu_{i}=Ff$ resp.  $\nu_{j}'\mu_{i}'=F'f$ generate the products.  For $f=0$ the products are  $0$. For $f\neq 0$
 we apply again the first part of the lemma to $x,z$ and the lifting $\tilde{x}$. We obtain
a chain
 of morphisms $\tilde{x} \rightarrow \tilde{z}'_{0} \rightarrow \tilde{z}'_{1} \rightarrow \ldots \tilde{z}'_{r}$ and the non-zero morphism $f$ starting in $\tilde{x}$ ends in $\tilde{z}'_{k}$ 
for a uniquely determined $k$.
Then we have  $rad^{j}A(y,z) \cdot rad^{i}A(x,y)=rad^{k}A(x,z)$ and $rad^{j}A'(y',z') \cdot rad^{i}A'(x',y')=rad^{k}A'(x',z')$.

In the second case  $A(y,z)$ is cyclic from the left but $A(x,y)$ is not. This time we choose a $\tilde{y}$ with $F\tilde{y}=y$ and we apply part a) of the 
lemma to obtain  points $\tilde{z}_{0},\tilde{z}_{1}, \ldots \tilde{z}_{s}$ and 
 morphisms
 $w_{0},w_{1}\ldots w_{s}$ with $w_{0}:\tilde{y}_{i} \rightarrow \tilde{z}_{0},w_{p}:\tilde{z}_{p-1} \rightarrow \tilde{z}_{p}$ such that
 $\nu_{j}:=F(w_{j}w_{j-1  \ldots }w_{0})$ generates $rad^{j}A(y,z)$ and $\nu_{j}':=F'(w_{j}w_{j-1  \ldots }w_{0})$ generates $rad^{j}A'(y',z')$. Dually we find points $\tilde{x}_{0},\tilde{x}_{1}, \ldots \tilde{x}_{t}$ and 
 morphisms
 $v_{0},v_{1}\ldots v_{s}$ with $v_{0}:\tilde{x}_{0} \rightarrow \tilde{y},v_{p}:\tilde{x}_{p} \rightarrow \tilde{x}_{p-1}$ such that $\mu_{i}=F(w_{0}w_{1} \ldots w_{i})$
 generates $rad^{i}A(x,y)$ and $\mu_{i}'=F'(w_{0}w_{1} \ldots w_{i})$
 generates $rad^{i}A'(x',y')$. We set $f= w_{j}\ldots w_{0}v_{0} \ldots v_{i}$. Then the product $rad^{j}A(y,z)\cdot rad^{i}A(x,y)$ is generated
 by $\nu_{j}\mu_{i}=Ff$ and $rad^{j}A'(y',z')\cdot rad^{i}A'(x',y')$ is generated
 by $\nu_{j}'\mu_{i}'=F'f$. For $f=0$ both products are $0$. For $f\neq 0$
 we apply again the first part of the lemma to $x,z$ and the lifting $\tilde{x_{j}}$. We obtain
a chain
 of morphisms $\tilde{x_{j}} \rightarrow \tilde{z}'_{0} \rightarrow \tilde{z}'_{1} \rightarrow \ldots \tilde{z}'_{r}$ and the non-zero morphism $f$ starting in $\tilde{x_{j}}$ 
ends in $\tilde{z}'_{k}$ 
for a uniquely determined $k$.
Then we have  again $rad^{j}A(y,z) \cdot rad^{i}A(x,y)=rad^{k}A(x,z)$ and $rad^{j}A'(y',z') \cdot rad^{i}A'(x',y')=rad^{k}A'(x',z')$.

In the last case $A(x,y)$ and $A(y,z)$ are both cyclic over $A(y,y)$. We choose $\tilde{x}$ with $F\tilde{x}=x$ and we get 
 a sequence of points  $\tilde{y}_{0},\tilde{y}_{1}, \ldots \tilde{y}_{t}$ and
  morphisms
 $v_{0},v_{1}, \ldots v_{t}$ and also points  $\tilde{z}_{0},\tilde{z}_{1}, \ldots \tilde{z}_{r}$ and
  morphisms
 $u_{0},u_{1}, \ldots u_{r}$ such that $F(v_{n}v_{n-1}\ldots v_{0})$ generates $rad^{n}A(x,y)$ for all $n\leq t$ and $F(u_{m}u_{m-1}\ldots u_{0})$ generates $rad^{m}A(x,z)$ for all $m\leq r$.
Analogous assertions hold in the dashed situation.
 We have $rad^{j}A(y,z) \cdot rad^{i}A(x,y) = A(y,z)\cdot N^{i+j}(y)\cdot A(x,y)$ and so it is generated as a vector space by all products $\psi F(v_{p}\ldots v_{0})$ with $i+j \leq p \leq t$
and $\psi \in A(y,z)$. For any $p$ as above there is the isomorphism $\bigoplus_{F\tilde{z}=z} \tilde{A}(\tilde{y}_{p},\tilde{z}) \simeq A(y,z)$. Therefore $rad^{j}A(y,z) \cdot rad^{i}A(x,y)$ is 
generated as a vector space by all $F(\tilde{\psi}v_{p} \ldots v_{0})$ with $\tilde{\psi}:\tilde{y}_{p} \rightarrow \tilde{z}$ for some $p \geq i+j$ and some $\tilde{z}$ in the fibre of $z$. 
 From $F(\tilde{\psi}v_{p} \ldots v_{0}) \neq 0$ we get $\tilde{z}=\tilde{z}_{i}$ for some $i$ between $0$ and $r$. We let $k$ be the smallest index occurring that way and we fix a 
$\tilde{\psi}:\tilde{y}_{p} \rightarrow \tilde{z}$ with $F(\tilde{\psi}v_{p} \ldots v_{0}) \neq 0$. Then $F(\tilde{\psi}v_{p} \ldots v_{0}) \neq 0$ generates $rad^{k}A(x,z)$ and 
 $rad^{j}A(y,z) \cdot rad^{i}A(x,y) = rad^{k}A(x,z)$ follows. Again the dashed situation is analogous.\hspace{0.1cm}q.e.d.

\subsection{Most representation-finite selfinjective algebras are regular and standard}

Fortunately almost all representation-finite selfinjective  algebras are regular. To see that $A$ is regular one shows that for all irreducible morphisms $\overline{\alpha} \in A(x,y)$ one has
 $N(y)\overline{\alpha}=\overline{
\alpha} N(x)$.

\begin{theorem}
 Let $A$ be a locally representation-finite selfinjective connected category. Then $A$ is regular except possibly for the tree class $D_{3m}$ and a large fundamental group.
\end{theorem}

Proof: For a small fundamental group we have $dim\, A(x,y)\leq 1$ for all points 
$x\neq y$
 and so $A$ is regular. By the table in section 7 the fundamental group is always small for the tree classes  $E_{6},E_{7},E_{8}$  and for type $D_{n}$ only the mentioned exceptions exist.

So let $A$ be an algebra of tree  class $A_{n}$  having a large fundamental group $\langle  \tau^{se} \rangle$ with $se< n$.
We work with the universal cover $\tilde{A}$ of $A$. For any $t$ we denote by $B_{t}$ the combinatorial quotient $\tilde{A}//\langle  \tau^{te} \rangle$. 
Theorem  says that $A$ is regular if $B_{s}$ is regular. We prove first that $B=B_{1}$ is regular.
 Let $\alpha:x \rightarrow y$ be an arrow in the quiver $Q$ of $A$. We prove $N(y)\overline{\alpha}\subseteq \overline{\alpha}N(x)$. The other inclusion is shown symmetrically. We can assume 
$x\neq y$ and $0\neq N(y)\overline{\alpha}$. Then we choose a path $p:y_{1} \rightarrow y_{2} \dots \rightarrow y_{r}$ in $Q$ such that $\overline{p}$ generates the ideal $N(y)$. 
We lift $p\alpha$ to a path $\tilde{p}\tilde{\alpha}$ in the quiver $\tilde{Q}$ of $\tilde{A}$. Then this path can be prolongated to a path 
$\tilde{q}:\nu^{-1}\tilde{y}_{r}\rightarrow   \ldots \rightarrow \tilde{x} \rightarrow \tilde{y}_{1} \rightarrow  \ldots \tilde{y}_{r}$ with $\overline{\tilde{q}} \neq 0$ in $\tilde{A}$.
Since $e$ divides $n$ the point $\nu^{-1}\tilde{y}_{r}$ is also mapped to $y$ in $Q$ and so $N(y)^{2}\neq 0$. Now we lift $p$ again this time with $\tilde{y}_{1}$ as thz ending point. 
Since $\tilde{y}_{1}$ is not a $D_{4}$-point in $\tilde{A}$ we find that 
$\tilde{y}_{r-1}\rightarrow \tilde{y}_{r}$ equals $\tilde{\alpha}$  and the wanted inclusion $N(y)\overline{\alpha}\subseteq \overline{\alpha}N(x)$ follows.

For general $s$ we have a Galois cover $B_{s} \rightarrow B_{1}$ where $B_{1}$ is regular and selfinjective and representation-finite.. This implies that $B_{s}$ is also regular. 
Namely one can reduce to the case where $B_{1}$ has only 
two points and then it follows trivially from 9.5 by inspection. \hspace{0.1cm}q.e.d.\vspace{0.2cm}

\begin{theorem}
 Let $A$ be a representation-finite regular selfinjective algebra with universal cover $\tilde{A}$ and fundamental group $\Pi$. Then $A$, its standard form $A^{s}$ and the 
combinatorial quotient $\tilde{A}//\Pi$ are isomorphic to each other.
\end{theorem}

Proof: Riedtmanns theorem induces a covering functor $F:\tilde{A} \rightarrow A$. Furthermore we have the two Galois covers $G:\tilde{A} \rightarrow A^{s}$ 
and $H:\tilde{A} \rightarrow \tilde{A}//\Pi$ introduced in 9.1. All three functors have the same fibres, whence they are all regular with the same stem category $\hat{A}$ by theorem 11.
Thus $A$,$A^{s}$,$\tilde{A}//\Pi$ and $k\hat{A}$ all satisfy the assumption of theorem 8 whence they are all isomorphic. \hspace{0.1cm}q.e.d.

\subsection{Coverings of base cateories and  quivers with relations}
Another combinatorial quotient occurs for coverings of base categories.
A functor $F:B' \rightarrow B$ between base categories is a covering if it satisfies the following conditions a), b) and the dual of b). Condition a) says that $F\mu=0$ is
 equivalent to $\mu=0$. Condition b) means that any point $x$ in $B$ can be lifted to a pont $x'$ and any $\mu:x \rightarrow y$ can be lifted to a unique $\mu':x' \rightarrow y'$.
 It follows that $\mu$ is irreducible iff $F\mu$ is, whence a covering induces a covering between the quivers. The covering $\pi:\tilde{B} \rightarrow B$ is 'the' universal covering if 
for any covering $F:B' \rightarrow B$ any $x$ in $B$ with liftings $x'$ in $B'$ and $\tilde{x}$ in $\tilde{B}$ there is exactly one functor $G:\tilde{B} \rightarrow B'$ with $\pi= F G$
 and $G\tilde{x}=x'$. Then $G$ is again a covering and  even an automorphism for $F=\pi$. The group of all these automorphisms is called the fundamental group $\Pi_{B}$ of $B$ and $B$
 is simply connected if this group is trivial.

Assume now that $B$ is connected i.e. its quiver $Q_{B}$ is connected. Then 
one obtains the universal covering by the following construction. A walk $w=\alpha_{n}\ldots \alpha_{1}$ of length $n$ from $x$ to $y$  is a formal composition of arrows $\beta$ in $Q_{B}$ and formal
inverses $\beta^{-1}$ such that the domains and codomains fit together well and $x$ is the domain of $\alpha_{1}$, $y$ the codomain of $\alpha_{n}$. Two walks $v$, $w$ can be composed to the 
walk
$wv$ if the end of $v$ is the start of $w$. The homotopy is the smallest equivalence relation on the set of all walks such that:1) $\alpha\alpha^{-1} \sim id_{y}$ and
$ \alpha^{-1}\alpha \sim id_{x}$  for all $\alpha:x \rightarrow y$, 2) $v \simeq w$ and $v^{-1} \simeq w^{-1}$ for all paths $v,w$ mapped to the same non-zero morphism under the canonical
 presentation $PQ_{B} \rightarrow B$ and 3) $v \sim w$ implies $uv \sim uw$ resp. $vu \sim wu$ whenever these compositions are defined. Now the points of the universal covering are the homotopy 
classes of walks with a fixed start $x$ and the fundamental group consists of the homotopy classes with start and end in $x$. The multiplication is induced by the composition of walks. Since 
$Q_{B}$ is connected this construction is essentially independent of the chosen base point $x$.

For any field $k$ the universal covering $p:\tilde{B} \rightarrow B$ induces a Galois- cover $kp:k\tilde{B} \rightarrow kB$ with group $\Pi$ to which theorem 1 applies in case $\Pi$ is torsion-free.
The problem is that $p$ might be the identity and so the theorem is useless. This situation occurs already in the famous example of Riedtmann that we discuss in the next section.

\subsection{Representation-finite selfinjective algebras with two simples}
To illustrate coverings we classify the  basic connected representation-finite selfinjective algebras $A=kQ/I$ with two simples. There are   four different types. 
That the listed algebras are representation-finite can be checked with the finiteness criterion of section 13.4, but we explain also in section 11
where they occur in the general classification.

 So the quiver $Q$ has two arrows 
$\alpha:x \rightarrow y$ and $\beta:y \rightarrow x$.
If there is no loop in $Q$ then $A$ is uniserial and defined by the vanishing of all paths of length $p$ for some $p\geq 2$.This is type 1.
 
 If $\gamma:y \rightarrow y$ is a loop there is no loop at $x$ and so $A$ is weakly symmetric.
Thus  we obtain a relation
$$\alpha\beta= \gamma^{n} + a_{n+1}\gamma^{n+1} + \ldots + a_{m}\gamma^{m}$$ for some $m\geq n\geq 2$ and $a_{j} \in k$.

In the sequel we consider algebras defined by zero-relations and the relation above with $n=m$. All these algebras have an obvious combinatorial Galois cover  $\tilde{A}_{n}$ with group $\Z$.
The quiver $\tilde{Q}$ is independent of $n$ and it has arrows $\gamma_{i}:y_{i} \rightarrow y_{i-1}$, $\beta_{i}:y_{i} \rightarrow x_{i}$, $\alpha_{i}:x_{i} \rightarrow y_{i-n}$ for $i \in \Z$.
 The relations are 
$$\alpha_{i}\beta_{i}=\gamma_{i+1-n}\ldots \gamma_{i-1}\gamma_{i}\,,\,i \in \Z $$ and the lifted zero-relations and the covering functor maps the $x_{i}$ to $x$, the $y_{i}$ to $y$. 

Suppose first that $\beta\gamma \in I$. Then $\gamma\alpha$ and $\gamma^{n+1}$ are both in $I$ which is generated by $\beta\gamma$, $\gamma\alpha$ and $\gamma^{n}- \alpha\beta$. This is type 2.

For $\beta\gamma \neq 0$ we suppose first that $\beta\alpha$ is not in the ideal $J$ generated by  all $\beta\gamma^{j}\alpha$, $j\geq 1$ and $I$. Then we add $\gamma^{n+1}$ to the relations and 
we find an algebra of the type described above, where the full subcategory of $\tilde{A}$ supported by $x_{i+n},x_{i},y_{i+1},y_{i},y_{i-1}$ is representation-infinite. 

In the remaining case we 
have  $\beta\alpha= \sum _{j=1}^{p}\,b_{j}\beta\gamma^{j}\alpha$. Replacing $\beta$ by $\beta - \sum _{j=1}^{p}\,b_{j}\beta\gamma^{j}$
we get$\beta\alpha=0$. For $n\geq 3$ we divide by $\gamma^{n+1}$ so that $\tilde{A}$ contains the representation-infinite algebra with support
 $x_{i+n},x_{i},y_{i+2},y_{i+1},y_{i},y_{i-1}$. Thus we have  $n=2$ and we end up with algebras $A(a)$ defined by the  relations $\beta\alpha=0$, $\alpha\beta=\gamma^{2} + a\gamma^{3}$ and
 $\gamma^{4}=0$. For $a=0$ the last relation follows and for $a\neq 0$ one has $A(a^{-1})\simeq A(1)$ as one sees by muliplying all three arrows by $a$.
 
 All $A(a)$ have the same ray category and the same stem category. The ray category $P$  is defined by $Q$ and the relations of $A(0)$ whereas the stem category $B$ is defined by $Q$ and the relations $\alpha\beta \sim \gamma^{2}, \beta\alpha \sim \beta\gamma\alpha, \gamma^{4 } \sim 0$. 
 Moreover  $A(0)\simeq kP$ is obvious and $A(1)\simeq kB$ follows if one replaces $\alpha$ by $\alpha - \gamma\alpha$, $ \gamma$ by $-\gamma$ and $\beta$  by
 $\beta$. 
 \begin{lemma}
 The linearizations $ kP$ and $kB$ are isomorphic iff the characteristic is not 2.
 
 \end{lemma}

Proof: Up to an inner automorphism any isomorphism $\phi:kP \rightarrow kB$ is induced by sending $\alpha$ to $x_{1}\alpha + x_{2}\gamma\alpha$,
$\beta$ to $y_{1}\beta + y_{2}\beta\gamma$ and $\gamma$ to $z_{1}\gamma + z_{2}\gamma^{2} + z_{3}\gamma^{3}$ for some scalars with
$x_{1}y_{1}z_{1}\neq 0$. The relations of $kP$ have to be mapped to $0$ in $kB$ and one obtains the three equations
$$x_{1}y_{1}=z_{1}^{2},x_{2}y_{1}+x_{1}y_{2}=2z_{1}z_{2}, x_{1}y_{1}+x_{1}y_{2}+x_{2}y_{1}=0.$$

There is no solution for $2=0$ and otherwise the solution $1=x_{1}=y_{1}=z_{1}=y_{2}=-x_{2}=-2z_{2}$, $y_{2}=z_{3}=0$.\hspace{0.2cm}q.e.d.

\section{Representation-finite trivial extensions}

\subsection{Trivial extensions and repetetive categories }
 First of all we define trivial extensions and repretitive categories in the language of locally bounded categories. Given such a 
 category $B$ and a functor $M:B^{op} \times B \rightarrow mod\, k$ the trivial extension $T(B,M)$ has the same objects as
$B$ and morphism spaces $T(B,M)(x,y)=B(x,y) \times M(x,y)$. The composition is given by $(g,n)(f,m)= (gf,M(id,g)m\, +\,M(f,id)n )$. 

The objects of the  repetetive category $\hat{T}(B,M)$ are the pairs $(x,i)$
with $x \in B$ and  $i \in \bf Z$. The morphism-sets are zero except for $\hat{T}(B,M)((x,i),(y,i)) = B(x,y) \times \{i\}$ and  
$\hat{T}(B,M)((x,i),(y,i+1)) = M(x,y) \times \{i\}$. 
The composition is defined by $(g,i)(f,i)=(gf,i)$, $(m,i)(f,i)=(M(f,id)m,i)$ and $(g,i+1)(m,i)=(M(id,g)m,i)$. 
 The functor $\hat{T}(B,M) \rightarrow T(B,M)$ forgetting the grading is the quotient for the obvious action of $\Z$ by shifting on $\hat{T}(B,M)$.
The repetitive category was introduced in the special case  $M(x,y)=DB(y,x)$ by Hughes and Waschbüsch in \cite{HW} to study  representation-finite trivial extensions. Then $\hat{T}(B,M)$ coincides with the universal cover of $T(B,M)$ ( see corollary 1 ), but later many more applications and some generalizations    were found.

\begin{proposition} Let $B$ be a connected locally bounded category and $M \neq 0$ such  that $T(B,M)$ is locally representation-finite. Then:
\begin{enumerate}
 \item $\hat{T}= \hat{T}(B,M)$ is locally representation-finite with connected Auslander-Reiten quiver $\hat{\Gamma}$ and there is a covering $\hat{\pi}:\hat{\Gamma} \rightarrow \Gamma$ with 
automorphism
 group $\hat{G} = \bf Z$.
\item  $\hat{\pi}$ is isomorphic to the universal cover $\tilde{\pi}:\tilde{\Gamma} \rightarrow \Gamma$ iff the automorphism group $\tilde{G}$ of $\tilde{\pi}$ is  $\bf Z$. 
In that case $B$ is directed and simply connected. 
\end{enumerate}

\end{proposition}

Proof: a)  For any $i$ the full subcategory $B(i)$ of objects $(x,i)$ is isomorphic to
$B$ under $\hat{\pi}$ and therefore a convex connected subcategory. Since $M(x,y)\neq 0$ for some $x,y$ each $B(i)$ is connected to $B(i+1)$
By theorem 10  $\hat{T}$ is locally representation-finite and so $\hat{\Gamma}$ is connected. Furthermore $\hat{\pi}$ exists with automorphism group  $\bf Z$.

b) By general covering theory ( \cite[section 1]{BG} ) there is  a subgroup $H$ of $\tilde{G}$  and a covering map $\pi': \tilde{\Gamma} \rightarrow \hat{\Gamma}= \tilde{\Gamma}/H$ with 
$\tilde{\pi}= \hat{\pi} \pi'$. 
Moreover $\hat{G}= N/H$. where $N$ is the normalizer of $H$ in $\tilde{G}$. Thus $\tilde{G}= \bf Z$ implies that $H=\{e\}$ and $\pi'$ is an isomorphism. Then $B$ is simply connected as
 a convex subcategory of the simply connected category $\tilde{T}$.\hspace{0.1cm}q.e.d.\vspace{0.2cm}

We are only interested in the case $M(x,y)=DB(y,x)$ where $T(B,M)$ is simply denoted by $T(B)$ and called the trivial extension of $B$. One has $T(B)(x,?)\simeq DT(B)(?,x)$ natural in $x$
whence $T(B)$ is  symmetric and in particular selfinjective.  

\begin{corollary}
 Let $B$ be  a connected algebra such that $T(B)$ is representation-finite. Then $\hat{T}(B)$ is the universal cover of $T(B)$.
\end{corollary}

Proof: Let $\Gamma$ be the Auslander-Reiten quiver of $T(B)$ and  let $\tilde{\pi}:\tilde{\Gamma} \rightarrow \Gamma$ be the universal cover. This induces the universal cover $\tilde{\pi}^{s}$ 
 between the stable parts and the automorphism group of $\tilde{\pi}$ embeds into that of $\tilde{\pi}^{s}$ which is cyclic by theorem 1  . Thus $\tilde{\Gamma}$ is isomorphic to $\hat{\Gamma}$ and 
the claim follows. \hspace{0.1cm}q.e.d.\vspace{0.2cm}

This is not true for any trivial extension: The category $A(1)$ defined at the end of 9.1 has a combinatorial quotient $A(1)//(\Z \times \Z)$ which is a trivial extension, but $A(1)$ is not 
a repetetive category.

A trivial consequence of the corollary is Yamagatas result \cite{Y} that the quiver of $B$ has no oriented cycle. The proof of the corollary in \cite[section 2.6]{HW} is not convincing and this leads there to a wrong statement in proposition 2.7.

\subsection{Trivial extensions of fundamental algebras }

In this subsection  $A(C)$ is a locally representation-finite simply connected  selfinjective category with  Nakayama-automorphism $\nu$.

\begin{proposition} For any fundamental algebra $F$ of $A(C)$ 
the combinatorial quotient  $A(C)// \langle \nu \rangle$ is equivalent to $T(F)$.
 
\end{proposition}
Proof: By the considerations of 9.1 the Schurian category $A(C)$ has a $\nu$-invariant multiplicative basis $B=\{w(x,y)\mid \,A(C)(x,y)\neq 0 \}$. For
 each element $w(x,y)$ with   $x,y \in F$  we define 
$w(x,y)^{\ast} \in DF(x,y)$  by $w(x,y)^{\ast} (w(x,y))=1$. Thus $w(x,y)^{\ast} \in TF(y,x)$ and  $TF(x,y)= F(x,y) \oplus DF(y,x)$ has dimension 1 for $x\neq y$ and 2 for $x=y$. 
The only possibly non-zero compositions are as follows:

$$ w(y,z)w(x,y) \neq 0  \Leftrightarrow w(y,z)w(x,y)=w(x,z)$$
$$ w(z,y)^{\ast}w(x,y) \neq 0  \Leftrightarrow w(x,y)w(z,x)=w(z,y) \Leftrightarrow w(z,y)^{\ast}w(x,y)= w(z,x)^{\ast}$$
$$ w(y,z)w(y,x)^{\ast} \neq 0  \Leftrightarrow w(z,x)w(y,z)=w(y,x) \Leftrightarrow w(y,z)w(y,x)^{\ast}= w(z,x)^{\ast}.$$

Next we construct a $\nu$-invariant functor $L:A(C) \rightarrow TF$. The objects of $A(C)$ are given by $\nu^{i}x$, the morphisms in B  either by $\nu^{i}w(x,y)$ or by $\nu^{j}w(z,\nu^{-1}t)$ for some 
uniquely determined $x,y.z,t$ in $F$ and $i,j \in \Z$.
 We set $L\nu^{i}x=x$, $L\nu^{i}w(x,y)=w(x,y)$ and $L\nu^{j}w(z,\nu^{-1}t)=w(t,z)^{\ast}$ where $w(t,\nu^{-1}t)=w(z,\nu^{-1}t)w(t,z)$ shows $w(t,z)\neq 0$. Then $L$ respects the multiplication.
For instance $$w(y,\nu^{-1}z)w(x,y) \neq 0 \Leftrightarrow w(z,\nu^{-1}z)=w(y,\nu^{-1}z)w(x,y)w(z,y)$$ and then $$Lw(y,\nu^{-1}z)Lw(x,y)=
w(z,y)^{\ast}w(x,y)=w(z,x)^{\ast}=Lw(x,\nu^{-1}z).$$
 Extending $L$ on each morphism space to a $k$-linear map we obtain
a $\nu$-invariant covering functor which induces the wanted equivalence. \hspace{0.1cm}q.e.d.\vspace{0.2cm}

\begin{proposition}Let $A(C)$ be a simply connected locally representation-finite selfinjective category  with Nakayama-automorphism $\nu$. Then we have:
\begin{enumerate}
 \item For any fundamental algebra $F$  with source $x$ the full subcategory $F'= \{\nu x\} \cup (F \setminus \{x\}$ of $A(C)$ is again a fundamental algebra. 
The dual statement holds for a sink and the two operations are inverse to each other.
\item Any two fundamental algebras can be transformed into each other  by applying  these operations. 
\item Any fundamenal algebra can be transformed into a pattern algebra.
\end{enumerate}

\end{proposition}

Proof: Clearly $F'$ is a representative system of the $\nu$-orbits. Let $y \rightarrow \nu x$ be an arrow in $A(C)$. Then $\nu ^{i}y \in F$ holds for some $i$. For $i > 0$ we get by the convexity of $F$ the contradiction $\nu x \in F$ and for $i < 0$ 
the point $x$ is not a source. Thus $y$ belongs to $F$. Similarly, let $x \rightarrow z$ be an arrow of $A(C)$. Then $A(C)(z,\nu x )\,\neq =$ implies $A(C)(z,y) \neq 0$ for some $y$ as before and so
$z \in F$ by the convexity of $F$. Now it follows easily that $F'$ is connected and convex.

Let $F_{1}$ and $F_{2}$ be two  fundamental algebras. We write $F_{1} \leq F_{2}$ if for all $x \in F_{1}$ we have $\nu ^{i(x)} \in F_{2}$ with $i(x) \geq 0$. In that case we define
 $d(F_{1},F_{2})$ as the sum of all $i(x)$. Now applying the operation of part a) often enough to $F_{2}$ we are reduced to the case $F_{1} \leq F_{2}$. If $d(F_{1},F_{2})=0$ then
$F_{1}=F_{2}$ and we are done. If $i(x) > 0$ for some source in $F_{1}$ we  have  $F_{1}' \leq F_{2}$ and $d(F_{1}',F_{2}) < d(F_{1},F_{2})$ so that induction applies. Thus assume $i(x)=0$ for all
 sources in $F_{1}$ i.e. all sources of $F_{1}$ are in $F_{2}$. Any  $y \in F_{1}$ has $\nu ^{i} \in F_{2}$ with $i\geq 0$. There is a path from a source $x \in F_{2}$ through $y$ to 
$\nu^{i}y \in F_{2}$. Since $F_{2}$ is convex $y$ lies in $F_{1}$.

Given a fundamental algebra $F$ take any section $S$ with $F \preceq S$ and $F\cap S \neq \emptyset$ and apply the procedure of the first part to a maximal element $x$ of $F$ if there is one with 
$ x \preceq \nu S$. After finitely many steps we obtain an $F'$ with $\nu S \prec F' \preceq S$ and this comes from a pattern algebra. \hspace{0.1cm}q.e.d.\vspace{0.2cm}

Two fundamental algebras related by finitely many transformations as above are called reflection-equivalent in \cite{HW}. Furthermore the fundamental algebras of tree class $T$ coincide by \cite{ass} with the iterated tilted algebras of the same type.

\subsection{The classification}

\begin{theorem}The map $A(C) \mapsto A(C)//\langle \nu \rangle$ induces a bijection between the isomorphism classes of locally representation-finite simply connected selfinjective categories
 and representation-finite connected trivial extensions.
\end{theorem}
Proof:  $A(C)// \langle \nu \rangle$ is a trivial extension by proposition 12. Since $A(C)$ is the universal cover of $A(C)//\langle \nu \rangle$ it follows  that $A(C)$ and $A(C')$ are 
isomorphic iff their combinatorial quotients  are. Finally,   a representation-finite connected trivial extension $T(B)$ has $\hat{T}(B)$ as its universal cover,
 which is isomorphic to  $A(C)$
for some configuration $C$. Now $B$ is a convex connected subcategory of $\hat{T}(B)$ consisting of representatives of $\nu$-orbits whence a fundamental algebra in $A(C)$. Thus $T(B)$ is isomorphic
 to $A(C)//\langle \nu \rangle
$. \hspace{0.1cm}q.e.d.\vspace{0.2cm}

The next result follows easily from our previous results:
\begin{theorem} Let $A,A'$ be connected algebras. Then:
\begin{enumerate}
 \item $T(A)$ is representation-finite iff $A$ is isomorphic to a fundamental algebra inside some $A(C)$.
\item $T(A)$ and $T(A')$ are isomorphic iff $A$ and $A'$ are isomorphic to fundamental algebras inside the same $A(C)$.
\end{enumerate}

\end{theorem}

\subsection{Brauer-quiver algebras and some  modifications}
For the tree classes $A_{n}$ and $D_{n}$ 
we describe the combinatorial quotients $A(C)//\langle \nu \rangle$  
as trivial extensions of appropriate fundamental algebras. We obtain Brauer-quiver algebras for $A_{n}$ and some modifications thereof for $D_{n}$.

A Brauer-quiver $Q$ is a finite connected quiver such that each arrow occurs in a  cycle, each point belongs to two cycles and two cycles have at most one point 
in common. Here  a cycle of length $n\geq 1$ is a path of length $n$ with the same start and end leading through $n$ points.

To each Brauer-quiver one defines a graph having the cycles as points. Two of these are joined by an edge if the cycles meet each other. This graph is a tree since
 otherwise some points lie on three cycles. Thus the sets of cycles resp. arrows in $Q$ are disjoint unions of two sets 
called $\alpha$-cycles and $\beta$-cycles resp. 
$\alpha$- arrows and  $\beta$-arrows. For a point $x$ in a cycle $Z$ we denote by $\zeta(x,Z)$ the path running once through the cycle starting and ending in $x$.

The Brauer-quiver algebra $B(Q)$ to a Brauer-quiver $Q$ is defined by two sorts of relations: Any path of length two not lying on a cycle vanishes and for any point $x$ one has
 $\zeta(x,Z) = \zeta(x,Z')$ for the two cycles $Z,Z'$ containing $x$. The Brauer quiver with one point is the double loop and we will exclude this trivial case in the following. 
Then there is at most one loop in each point and the reduced Brauer quiver $Q'$ to $Q$ is defined by removing all loops.  It  is the quiver of $B(Q)$. 
The relations that any path not 
lying on a cycle vanishes and $\zeta(x,Z) = \zeta(x,Z')$ for two cycles $Z,Z'$ of lenghts $\geq 2$ containing $x$
are preserved but  one has to add that all $\zeta(x,Z)$ are annihilated by all arrows. Denoting by  $I$ the ideal generated by these relations  we have the presentation $kQ'/I \simeq B(Q)$.

An arrow $\gamma$ in a reduced 
Brauer-quiver is called special if it belongs to a cycle of length $3$ at least and if both end-points have order $2$.

\begin{theorem} The trivial extension $T(A_{L})$ of a pedigree algebra $A_{L}$ is a Brauer-quiver algebra and any Brauer-quiver algebra occurs that way up to isomorphism. 
  
 \end{theorem}

Proof: Any $\alpha$-string from $a$ to $b$ in $L$ induces a complete morphism in $A_{L}$ whence an 
arrow $b \rightarrow \nu a$ in $\hat{T}(A_{L})$ and finally an arrow $\overline{b} \rightarrow \overline{a}$ in the quiver of $T(A_{L})$ which is also called an $\alpha$-arrow. 
One obtains that way an $\alpha$-cycle.
 An analogous statement holds for $\beta$-strings and the quiver $Q'$ of $T(A_{L})$ is a finite connected quiver such that each arrow lies on a cycle, each point belongs to two cycles at most
 and two cycles intersect in at most one point. 
The relations of $T(A_{L})$ are that paths not lying in a  cycle are zero, that $\zeta(x,Z)= \zeta(x,Z')$ for any point in two cycles and that any $\zeta(x,Z)$ is
 in the two-sided socle. By adding a loop in each point not belonging to two cycles we obtain a Brauer-quiver $Q$ such that $B(Q)$ is isomorphic to $T(A_{L})$.

Reversely, given a Brauer-algebra $B(Q)$ with $\alpha$- and $\beta$- arrows one chooses a point $\omega$ in $Q$ and one looks at  all walks in $Q$ that start in $\omega$,
  do not run twice through the same point and contain only  $\alpha$-arrows pointing away from $\omega$  and  $\beta$-arrows pointing to $\omega$.
The vertices of $Q$ and the arrows occurring in one of these walks form a pedigree $L$ with $T(A_{L}) \simeq B(Q)$.\hspace{0.1cm}q.e.d.\vspace{0.2cm}

For the tree $D_{n}$ we have constructed in section 6 fundamental algebras of the simply connected categories $A(C)$ as one-point extensions of pedigree algebras $A_{L}$. 
Therefore the 
trivial extensions of class $D_{n}$ are just some  modifications of Brauer quiver algebras as described carefully in \cite[section 7]{BLR}. 

We explain this only for three-cornered categories with $rst \neq 0$ since we need these in section 11.
For the special case $L=L(r,s,t)$ with $rst \neq 0$ figure 7 shows from 
the left to the right the quivers of $A_{L}$, its one-point exrension $B$ which is a fundamental algebra and $T(B)$. Here $a_{0} \rightarrow a$ is for $r \geq 2$ a short notation for the path 
 $a_{0} \rightarrow a_{1} \rightarrow \ldots  \rightarrow a_{r}=a$ and similarly for $b_{s} \rightarrow b=b_{0}$ and $b=c_{0} \rightarrow c_{t}$. The algebra $B$ has only the
 three complete morphisms $b_{s} \rightarrow b \rightarrow a$, $b \rightarrow a \rightarrow x$ and $a_{0} \rightarrow a \rightarrow x$ inducing three arrows
 connecting the quiver of $B$ with the quiver of $\nu B$.
Therefore 
the quiver $Q$ of $T(B)$ consists of the 
central $\gamma$-triangle formed by $\gamma_{1}:x \rightarrow b, \gamma_{2}:a \rightarrow x, \gamma_{3}:b \rightarrow a$ and three outer cycles $Z_{i}$ containing $\gamma_{i}$ for $i=1,2,3$.
For simplicity  we have denoted the $\nu$-orbit of a point $y$ in $A(C)$ again by $y$.   

\vspace{0.5cm}
\setlength{\unitlength}{0.5cm}
\begin{picture}(23,6)
\put(0.5,5){\circle*{0.2}}
\put(3.5,5){\circle*{0.2}}
\put(2,3){\circle*{0.2}}
\put(2,-1){\circle*{0.2}}
\put(3.5,1){\circle*{0.2}}
\put(0,5.5){$b_{s}$}
\put(3.8,0.7){$a$}
\put(1,-1){$a_{0}$}
\put(1,2.7){$b$}

\put(3.5,5.5){$c_{t}$}

\put(11.5,5){\vector(3,-4){1.5}}
\put(11.5,1){\vector(3,4){1.5}}
\put(8.5,5){\vector(3,-4){1.5}}
\put(10,3){\vector(3,-4){1.5}}
\put(10,3){\vector(3,4){1.5}}
\put(10,-1){\vector(3,4){1.5}}
\put(13,3){\circle*{0.2}}
\put(8.5,5){\circle*{0.2}}
\put(11.5,5){\circle*{0.2}}
\put(10,3){\circle*{0.2}}
\put(10,-1){\circle*{0.2}}
\put(11.5,1){\circle*{0.2}}
\put(8,5.5){$b_{s}$}
\put(11.8,0.7){$a$}
\put(9,-1){$a_{0}$}
\put(9,2.7){$b$}
\put(13.5,2.7){$x$}
\put(11.5,5.5){$c_{t}$}

\put(0.5,5){\vector(3,-4){1.5}}
\put(2,3){\vector(3,-4){1.5}}
\put(2,3){\vector(3,4){1.5}}
\put(2,-1){\vector(3,4){1.5}}

\put(16.5,5){\circle*{0.2}}
\put(19.5,5){\circle*{0.2}}
\put(22.5,5){\circle*{0.2}}
\put(18,3){\circle*{0.2}}
\put(21,3){\circle*{0.2}}
\put(19.5,1){\circle*{0.2}}
\put(16,5.5){$b_{s}$}
\put(19.5,5.5){$a$}
\put(22.5,5.5){$a_{0}$}
\put(17,2.7){$b$}
\put(21.5,2.7){$x$}
\put(19.5,0.5){$c_{t}$}

\put(19.5,5){\vector(-1,0){3}}
\put(22.5,5){\vector(-1,0){3}}
\put(21,3){\vector(-1,0){3}}
\put(16.5,5){\vector(3,-4){1.5}}
\put(18,3){\vector(3,-4){1.5}}
\put(19.5,1){\vector(3,4){1.5}}
\put(21,3){\vector(3,4){1.5}}
\put(18,3){\vector(3,4){1.5}}
\put(19.5,5){\vector(3,-4){1.5}}

\put(10,-3){figure 11}
\end{picture}
\vspace{2cm}

To get the relations for $T(B)$ recall that it is regular with  a cancellative basis inherited from $\hat{T}(B)$. So the zero-relations $b_{1} \rightarrow b \rightarrow c_{1}$,
 $c_{t} \rightarrow x \rightarrow \nu a_{0}$ and $a_{r-1} \rightarrow a \rightarrow \nu b_{s}$ holding upstairs induce downstairs the zero-relations 
$b_{1} \rightarrow b \rightarrow c_{1}$,
 $c_{t} \rightarrow x \rightarrow  a_{0}$ and $a_{r-1} \rightarrow a \rightarrow  b_{s}$. Similarly we get  three commutativity relations 
$a \rightarrow b_{s} \rightarrow b = a \rightarrow x \rightarrow b$, $b \rightarrow c_{t} \rightarrow x = b \rightarrow a \rightarrow x$ and 
$x \rightarrow a_{0} \rightarrow a = x \rightarrow b \rightarrow a$.
Finally we have to add the zero-relations saying that any cyclic path is in the socle of $T(B)$.

For a general pedigree $L$ still with $rst \neq 0$ we  consider the three connected components of $L \setminus \{b\}$. Let $C_{1}$ be the connected component of $c_{1}$, $C_{2}$ that of $b_{1}$
 and $C_{3}$ that of $a_{1}$. Furthermore $Q_{1}$ denotes the full subquiver of $T(B)$ supported by $C_{1} \cup \{b,x\}$, $Q_{2}$ the one supported by $C_{2} \cup \{a\}$ and $Q_{3}$ that supported by 
$C_{3} \cup \{x\}$. Then each $Q_{i}$ is a reduced Brauer-quiver containing $Z_{i}$ and the special arrow $\gamma_{i}$. These three quivers are glued to the quiver of $T(B)$ along the
 $\gamma$-triangle as in the case $L=L(r,s,t)$ and the relations are those gotten above supplemented by the relations in the three Brauer-quiver algebras. Reversely, any algebra $D(Q_{1},Q_{2},Q_{3})$ which
 is obtained from three reduced Brauer-quivers $Q_{i}$ with special arrow $\gamma_{i}$ by glueing the quivers along the $\gamma$-triangle and taking the above relations comes from a uniquely 
determined pedigree $L$.

It is obvious that  $D(Q_{1},Q_{2},Q_{3})$ is isomorphic to $D(Q_{\pi 1},Q_{\pi 2},Q_{\pi 3})$ for any cyclic permutation $\pi$. 
Furthermore the $\gamma$-triangle is the only triangle containing only $\nu$-orbits of $D_{4}$-points. Therefore $ D(Q_{1},Q_{2},Q_{3})$ and $D(Q'_{1},Q'_{2},Q'_{3})$  are isomorphic iff
 - up to cyclic permutation - $Q_{i}$ is isomorphic to $Q'_{i}$ for $i=1,2,3$ as a reduced Brauer quiver with special arrow.
So we have:

\begin{theorem}
 The three-cornered trivial extensions with $rst \neq 0$ are the algebras $D(Q_{1},Q_{2},Q_{3})$ given by the glueing of three reduced Brauer-quivers $Q_{i}$ with 
special arrow $\gamma_{i}$
along a $\gamma$-triangle. Two of them are isomorphic iff the conditions above are satisfied.
\end{theorem}

\section{Large fundamental groups}
By section 7 a large fundamental group $\Pi$ is only possible for the trees  $T= A_{n}$ or $T=D_{3m}$ and for $\Pi= \langle \tau^{re} \rangle $ with  $e \leq re < L(T)$. 
For $\langle \nu \rangle \subset \Pi$  we calculate the combinatorial quotient $A(C)//\Pi$ as the  quotient of the  trivial extension $A(C)//\langle \tau^{L(T)} \rangle$ by the 
 the finite group $\Pi/\langle \tau^{L(T)} \rangle$.

 Recall that a selfinjective algebra is called weakly symmetric provided the Nakayama-permutation $\nu$ 
is the identity. For $A(C)$ we have  $\nu =\tau^{L(T)}$ and so the quotient $A(C)//\langle \tau^{re} \rangle $ is weakly symmetric iff $re$ divides $L(T)$ iff $\langle \nu \rangle \subseteq
 \Pi$. A result of Kupisch in \cite{Kupisch1} says 
that a distributive algebra $A$ is weakly symmetric iff $k\hat{A}$ is symmetric.  Thus by theorems 13 and 14 all regular weakly symmetric representation-finite algebra are even symmetric. This holds also for non-regular algebras as the classification in 11.2 implies.

\subsection{Tree class A: Brauer-quivers with exceptional cycle}

Here we know already that the Auslander-Reiten quiver determines the algebra and so we only have to look at combinatorial quotients. We start with the symmetric ones. Most of the following is already contained in \cite{GRie}.

Let $Z$ be a cycle in a Brauer quiver $Q$ and let $m$ be a natural number. The Brauer-quiver algebra $B_{m,Z}(Q)$ with exceptional cycle $Z$ and multiplicity $m$ has the quiver $Q$. 
Each $x \in Z$ 
belongs to one other cycle $Z'$ and  the relation $\zeta(x,Z)=\zeta(x,Z')$ is replaced by $(\zeta(x,Z))^{m}=\zeta(x,Z')$. All other relations defining $B(Q)$ remain untouched. Note that $Z$ might be a
 loop and that for $m=1$ we get the usual Brauer-quiver algebra.

\begin{proposition}
\begin{enumerate}                  
 \item Let $C$ be a configuration  of $\Z \,A_{n}$ with $n=d\cdot m$,  $m >1$ and  $\tau^{d}C=C$.  Then there is exactly one cycle $Z$ in the Brauer-quiver $Q$ of 
$A(C)//\langle \tau^{n} \rangle $ fixed by $\tau^{d}$. The cycle $Z$ has $f\cdot m$ points for some $d \geq f \geq 1$ and $A(C)//\langle \tau^{d} \rangle $ is isomorphic to 
$B_{m,\overline{Z}}(\overline{Q})$. Here   $\overline{Q}$ resp. $\overline{Z}$ are the quotients of $Q$ resp. $Z$ under the action of $\tau^{d}$.
\item Any Brauer-quiver algebra with exceptional cycle  is obtained in this way. Therefore it is representation-finite. 
\item In the first part $d$ is the period of $C$ iff $B_{m,\overline{Z}}\overline{Q}$ has no non-trivial combinatorial automorphism.

\end{enumerate}

\end{proposition}

Proof: $\tau^{d}$ induces a fixed-point free automorphism $\epsilon$ of order $m$ on the reduced Brauer-quiver algebra $ B(Q) \simeq A(C)//\langle \tau^{n} \rangle $ and therefore also 
an automorphism 
$\epsilon'$ of the tree corresponding to $Q$. Any automorphism of a finite tree fixes an edge or a point. Fixing an edge means fixing the common point of two cycles which is impossible.
Thus $\epsilon'$ fixes a point i.e. there is a cycle $Z$ in $Q$ with $\epsilon (Z)=Z$. Another fixed cycle  would lead to a fixed edge in the tree.

Now we have a Galois covering
$\pi: A(C)//\langle \tau^{n} \rangle  \rightarrow A(C)//\langle \tau^{d} \rangle $ with group $\Z / m\Z = \langle \tau^{d} \rangle/\langle \tau^{n} \rangle$. The cycle $Z$ is mapped to 
a cycle $\overline{Z}$ of length $f$ and any  other cycle $Z'$ is  mapped to a cycle $\overline{Z'}$ of the same length whose preimage is a disjoint union of $m$ cycles. 
Thus $\overline{Q}= \pi\, Q$ is  a union of cycles where two of them intersect in at most one point.  For $x \in Z$ the relation $\zeta(x,Z)=\zeta(x,Z')$ 
becomes the relation $\zeta(\pi x, \pi Z)^{m}= \zeta (\pi x, \pi Z')$ and it follows that the quotient is isomorphic to $B_{m,\overline{Z}}\overline{Q}$.

 In the other direction for a given $B=B_{m,Z}(Q)$ with $d$ points, $m>1$ and exceptional cycle of length $f$ there is  an obvious covering with group $\Z\,/m\Z$
  $\pi:\tilde{Z} \rightarrow Z$. Namely interprete a cycle of length $l$ as $\Z/ l\Z$ with arrows $i \rightarrow i+1$ and take for $\pi$ the projection $\Z/mf\Z \rightarrow \Z/f\Z$. 
This covering can be extended in a unique way to
 covering $\tilde{Q}$ of $Q$ still with group $\Z\,/m\Z$. The relations can be lifted to define a Brauer-quiver algebra $\tilde{B}=B(\tilde{Q})$ with  a fixed point-free combinatorial automorphism 
$\epsilon$ of order $m$ and 
 $B$ is the quotient. Now $\tilde{B}$ is representation-finite. By  \cite{delapenaMartinez} $\epsilon$  acts also freely on $\tilde{B}$-mod whence 
$B$ is representation-finite by theorem 10. We know that $\tilde{B}\simeq A(C)//\langle \tau^{n} \rangle$ for some configuration $C$ of period $e$. 
At the level of Auslander-Reiten quivers we have a chain of coverings 
$\Gamma_{A(C)} \rightarrow \Gamma_{\tilde{B}} \rightarrow \Gamma_{B}$. Since $A(C)$ is  simply connected it is  the universal cover of $\Gamma_{B}$.
 By the table in section 7 the fundamental group $\Pi$  is generated by $\tau^{d'}$ because $\Pi$ is large.
 Furthermore  it contains $\tau^{n}$ and the quotient has $m$ elements whence $\Pi$ is  generated by $\tau^{d}$  
 and so $\tau^{d}C=C$.

Finally, if $d=se$ for some $s>1$ then $\tau^{e}$ induces a non-trivial combinatorial automorphism on $A(C)//\langle \tau^{re} \rangle=B_{m,\overline{Z}}\overline{Q}$. Reversely for 
any such automorphism the proper quotient still has $A(C)$ as its universal cover and so $e<d$ follows.
 \hspace{0.1cm}q.e.d.\vspace{0.2cm}.

The last part of the proposition implies that for each factorization $n=dm$ of $n$ there are configurations with period $d$. One can take the Brauer-quiver consisting of one $\alpha$-cycle 
of length $d$ and one exceptional $\beta$-loop of multiplicity $m$. 

In section 9.4 the algebras of type 2 are special cases of this whereas those of type 1 consist  for $p=2m$ just of one
 exceptional cycle of length $2$ with multiplicity $m$.

If $\Pi=\langle \tau^{re} \rangle$ is large and $re$ does not divide $n$  then $A(C)//\Pi $ is not 
a Brauer-quiver algebra but it is an $r$-fold covering of $A(C)//\langle \tau^{e} \rangle$ which is a Brauer-quiver algebra with exceptional cycle.

\subsection{Typ $D_{3m}$: Non-regular and non-standard algebras}
In the following we use the notations and results of sections 6.3, 7 and 10.4. 
By the table in section 7 large fundamental groups occur only for configurations $C$ of $D_{2m-3}$ and they are generated by $\sigma=\tau^{2m-1}$ or by  $\sigma^{2}$.  Since any  $C$ is  stable under $\nu=\sigma^{3}$ it is stable under $\sigma$ also in the second case. Thus  $A(C)$ contains at least three $\nu$-orbits of 
$D_{4}$-points in $A(C)$. By proposition 9 $A(C)$ is three-cornered with $rst\neq 0$  and by theorem 18 the quotient $A(C)//\langle \nu \rangle$ is  isomorphic to some $D(Q_{1},Q_{2},Q_{3})$. 
 Here  $\sigma$ 
induces a combinatorial automorphism  of order $3$ which respects the $\gamma$-triangle since its  vertices come from the $D_{4}$-points. Thus there is a reduced Brauer-quiver 
$Q$ with $m+1$ points having a special
 arrow and for $i=1,2,3$ isomorphisms $Q_{i} \simeq Q$  respecting the special arrows. 

Reversely, for any reduced Brauer-quiver $Q$ with special arrow and $m+1$ points the algebra $D(Q,Q,Q)$ admits an obvious combinatorial automorphism $\sigma$ of order three and $D(Q,Q,Q)$ 
is isomorphic to $A(C)//\langle \nu \rangle$ for some $C$ whence representation-finite. It follows as in 11.1 that $\sigma$ acts freely  on the indecomposables, that the quotient 
$D(Q)=D(Q,Q,Q)//\langle \sigma \rangle$ is  representation-finite with universal cover $A(C)$ and that $C$ is stable under $\tau^{2m-1}$.
So one gets a bijection beween the isomorphism classes of configurations in $\Z\,D_{3m}$ stable under $\tau^{2m-1}$ and those of reduced Brauer-quivers with $m+1$ points having a special arrow.

For each $C$ stable under $\sigma$ there is one algebra with Auslander-Reiten quiver $(\Z\,D_{3m})_{C}/\langle \sigma \rangle$ namely $D(Q)$, but there might be other 'candidates'. 
By  theorem 12  each  candidate $A$ and $D(Q)$ have isomorphic stem and ray categories whence 
they have the same quiver $\overline{Q}$, the same Cartan-invariants $c(p,q)=dim\, D(Q)(p,q)$ and both $A$ and $D(Q)$ are weakly symmetric.

To formulate our findings and to prove them we work with the reduced Brauer-quiver $Q$ with special arrow $\gamma$. In the notation of 10.5 we choose $Q_{1}$ for $Q$ and $Z_{1}$ for 
the $\alpha$-cycle of length $t+1$ 
containing the points $c_{0}=b,c_{1},\ldots c_{t+1}:=x$ and the arrows $\alpha_{i}:c_{i-1} \rightarrow c_{i}$ for 
$i=1,2,\ldots t$ as well as the special arrow $\gamma=\alpha_{t+2}:x \rightarrow c_{0}$.
The arrows of $Q$ form  a representative system for the arrows of  the quiver $\overline{Q}$  of $D(Q)$ and  we keep the notation.
Moreover only  the points   $c_{0}$ and $x$  have to be identified in $\overline{Q}$ and we denote the contracted point still by $c_{0}$.

Thus the $\alpha$-cycle $Z_{1}$ in $Q$ of length $t+1$ 'decomposes' in $\overline{Q}$ into the loop $\gamma$ at $c_{0}$ and the $\alpha$-cycle $\Z'_{1}$ of length $t$ containing the points $c_{0},c_{1},\ldots c_{t}$.
This is shown in figure 12. 

Note that $\gamma$ is the only loop in $\overline{Q}$ because $Q$ is a reduced Brauer-quiver.
For a point $c_{i}$ with $i\neq 0$ the path $\zeta(c_{i},Z'_{1})$ is defined as   $\alpha_{i-1} \dots \alpha_{1}\gamma\alpha_{t}\ldots \alpha_{i}$ as a path is $k\overline{Q}$
The other $\alpha$- and $\beta$-cycles  and also the definition  of the  path $\zeta(x,Z)$ are preserved in $\overline{Q}$.

Now we define several relations as elements of $k\overline{Q}$.

R1(a): $\alpha_{t}\alpha_{t-1} \ldots \alpha_{1} -\gamma^{2}- a\gamma^{3}$ for a scalar $a$; for $a\neq 0$ we add also  $\gamma^{4}$. 

R2: $\alpha_{1}\alpha_{t}$.

R3: All products in any order of an $\alpha$-arrow and a $\beta$-arrow.

R4: For any $p \neq c_{0}$ all products in any order of a $\zeta(p,Z)$ and an arrow.

R5: $\zeta(p,Z) - \zeta(p,Z')$ for all $p\neq c_{0}$ lying on two cycles $Z,Z'$ of $\overline{Q}$.

Define $A(a)$ as the quotient of $k\overline{Q}$ defined by the relations R1(a),R2,R3,R4 and R5 given above.

\setlength{\unitlength}{1cm}
\begin{picture}(15,8)

\put(8,4){\circle*{0.1}}
\put(7,2){\circle*{0.1}}
\put(5,1){\circle*{0.1}}
\put(0,4){\circle*{0.1}}
\put(3,1){\circle*{0.1}}
\put(1,2){\circle*{0.1}}
\put(1,6){\circle*{0.1}}
\put(3,7){\circle*{0.1}}
\put(5,7){\circle*{0.1}}
\put(7,6){\circle*{0.1}}
\put(8.5,4){\circle{1}}
\put(7,6){\vector(1,-2){1}}
\put(8,4){\vector(-1,-2){1}}
\put(8.0,3.9){\vector(0,1){0.1}}
\put(9.2,4){$\gamma$}
\put(5,7){\vector(2,-1){2}}
\put(7,2){\vector(-2,-1){2}}
\put(5,1){\vector(-2,0){2}}
\put(3,7){\vector(2,0){2}}
\put(3,1){\vector(-2,1){2}}
\put(1,6){\vector(2,1){2}}
\multiput(0,4)(0.1,0.2){10}{\circle*{0.01}}
\multiput(0,4)(0.1,-0.2){10}{\circle*{0.01}}

\put(7,5){$\alpha_{t}$}
\put(7,3){$\alpha_{1}$}
\put(7.5,3.9){$c_{0}$}
\put(6.3,2){$c_{1}$}

\put(6.3,6){$c_{t}$}
\put(8,3){$t\geq 2,\alpha_{1}\alpha_{n}=0$}
\put(8,2.5){$\alpha_{n} \ldots \alpha_{1}=\gamma^{2}$}
\put(7.5,2){$  \alpha_{i}\ldots \alpha_{1}\gamma\alpha_{n}\ldots \alpha_{i+1} \in S$}

\put(5,0){figure 12}

\end{picture}

\begin{theorem} Let $C$ be a $\sigma$-stable configuration in $\Z D_{3m}$ and let $Q$ be the corresponding reduced Brauer-quiver with special arrow $\gamma$. Then the following holds:

\begin{enumerate}
 \item A basic  algebra $A$ has  Auslander-Reiten quiver $(\Z\,D_{3m})_{C}/\langle \sigma \rangle$ iff it is isomorphic to $A(0)$ or $A(1)$.
\item $A(0)$ and $A(1)$ are isomorphic iff the characteristic is not $2$. Otherwise $A(0)$ is the linearization of the ray category of $D(Q)$ whereas $A(1)$ is the linearization of the stem category. 
\item Any algebra with Auslander-Reiten quiver $(\Z\,D_{3m})_{C}/\langle \sigma^{2} \rangle$ is regular.
\end{enumerate}
\end{theorem}

Proof: $D(Q,Q,Q)$ is a trivial extension and so the non-zero entries of the Cartan-matrix are $2$ on the diagonal and $1$ outside. For $D(Q)=D(Q,Q,Q)//\langle \sigma \rangle$ one 
gets by a short calculation the Cartan-numbers  $c(c_{0},c_{0})=4$, $c(c_{0},c_{i})=2=c(c_{i},c_{0})$ for $1 \leq i \leq t$, $c(p,p)=2$ for $p\neq c_{0}$ and finally
 $c(p,q)=1$ if $p \neq q$ are both different from $c_{0}$ and belong to the same cycle. 

Now we take any basic distributive weakly symmetric algebra $A$ with quiver $\overline{Q}$ and the above cartan-numbers.
We start with an arbitrary presentation $k\overline{Q} \rightarrow A$ written as $w \mapsto \overline{w}$ mapping  and modify it step by step until  the wanted relations are in the kernel.

Since $A(c_{0},c_{0})$ is uniserial of dimension 4 and $\overline{\gamma}$ is nilpotent 
the non-zero powers $\overline{\gamma}^{i}$ with $0 \leq i \leq 3$ form a basis of $A(c_{0},c_{0})$ and
 so we have $\overline{\alpha_{t}\alpha_{t-1}\ldots\alpha_{2}\alpha_{1}}= \overline{b \gamma^{2}+c \gamma^{3}}$ for some scalars $b,c$. Similarly the non-zero elements 
 $\overline{\alpha_{1}\gamma^{i}}$ are linearly independent in the twodimensional space $A(c_{0},c_{1})$ and so $\overline{\alpha_{1}\gamma^{2}}=0$.
Multiplication of the last equation from the left with $\overline{\alpha_{1}}$ implies  that $0=\overline{\alpha_{1}\alpha_{t}\alpha_{t-1}
\ldots \alpha_{2}\alpha_{1}}$ and so $A(c_{0},c_{1})$ is not cyclic over $A(c_{1},c_{1})$ but over $A(c_{0},c_{0})$ because $A$ is distributive.
In particular $\overline{\alpha_{1}\gamma}\neq 0$. By  the shape of $\overline{Q}$ and by the vanishing of most $c(p,q)$ the  paths $v=\alpha_{t}\alpha_{t-1} \ldots \alpha_{2}$
and $\gamma v$ give a basis of $A(c_{2},c_{0})$.

 Since $A$ is weakly symmetric we get $\overline{\alpha_{1}\gamma v} \neq 0$  as well as $\overline{v \alpha_{1} \gamma}\neq 0$. The second relation implies $b\neq 0$ and the first 
$\overline{\alpha_{1}\alpha_{t}}= d\overline{\alpha_{1}\gamma\alpha_{t}}$ for some scalar $d$ because 
$0 \neq \overline{\alpha_{1}\gamma\alpha_{t}}$ is a basis of the onedimensional space $A(c_{t},c_{1})$. Replacing $\overline{\alpha}_{1}$ by
 $\overline{\alpha_{1}}
 -d\overline{\alpha_{1}\gamma}$ and $\overline{\gamma}$ by an appropriate scalar multiple we obtain a presentation satisfying 
$\overline{\alpha_{1}\alpha_{t}}=0$ and 
$\overline{\alpha_{t}\alpha_{t-1}
\dots  \alpha_{2}\alpha_{1}}= \overline{\gamma^{2}+ a \gamma^{3}}$. 

If the charakteristic is not 2 we replace $\overline{\gamma}$ by $\overline{\gamma + \frac{a}{2} \gamma^{2}}$ and then $R1(0)$ is in the kernel.
If $a\neq 0$ we multiply $\overline\alpha_{1}$, $\overline{\alpha_{2}}$ and $\overline{\gamma}$ by $a$ and the new presentation satisfies the relations $R1(1)$ and $R2$.

The relations R3 hold for any presentation because $c(p,q)=0$ for points not in the same cycle. Also the relations $R4$ hold for presentations  satisfying $R2$ and some R1(a). 
Namely the first arrow  in   $\zeta(p,Z)$ is not in the kernel of the presentation
 and so it can be prolongated to a path from $p$ to $p$ not in the kernel because $A$ is weakly symmetric. Then $\zeta(p,Z)$ is the only 
possible path with this property because of $c(x,y)=0$ for points in different cycles and because R2 holds. Thus $\overline{\zeta(p,Z)}$ is a basis of the one-dimensional radical of $A(p,p)$.
 We  show that this is annihilated from the left by $\overline{\epsilon}$ for any arrow $\epsilon:p \rightarrow q$. This is clear if $q$ does not belong $Z$ or if $q\neq c_{0}$ because then $c(p,q)=1$. In the remainig case
 we have 
$\zeta(p,Z)=\alpha_{t-1}\alpha_{t-2}\ldots \gamma\alpha_{t}$. Then $\overline{\alpha_{t} \zeta(p,Z)}=0$ by R1(a) since $\overline{\gamma^{2}\alpha_{t}}=0$. Symmetrically one gets that 
$\overline{\zeta(p,Z)}$ is annihilated by all $\overline{\epsilon}$ from the right.

To deal with the relations R5 we consider the tree $G$ associated to $Q$. It has  the cycles as points and an edge between $Z$ and $Z'$ iff the cycles have a point in common. Thus the edges
 are in bijection to the relations in R5. We number the cycles
$Z_{1},Z_{2},\ldots Z_{m}$ in such a way that the points $Z_{1},Z_{2},\ldots,Z_{i}$ form for any $i$ a subtree $G_{i}$ of $G$. Suppose we have already found a presentation satisfying all relations
corrwsponding to edges in $G_{i}$. For $i<m$ there is a cycle $Z_{j}$ with $j \leq i$ and a point $p$ in $Z_{j}$ and $Z_{i+1}$. Then  the elements
$\overline{\zeta(p,Z_{j})}$ and  $\overline{\zeta(p,Z_{i+1})}$ are non-zero and proportinal. There is an arrow $\epsilon$ in $Z_{i+1}$ not occurring in any $Z_{l}$ with $l\leq i$ and we
 multiply it by an
 appropriate non-zero scalar to get $\overline{\zeta(p,Z_{j})} = \overline{\zeta(p,Z_{i+1})}$ in the new presentation. By induction
we get a presentation satisfying all relations  listed in  R5. Since we have not modified the presentation on arrows occurring in R1 or R2 
whereas R3 and R4 always hold all relations are now satisfied.

A simple verification shows that all $A(a)$ have the same  Cartan-numbers as $D(Q)$ and so $A(1)$ or $A(0)$ is isomorphic to $D(Q)$ which has the wanted Auslander-Reiten quiver.
Now the quotients of $A(1)$ and $A(0)$ by $\overline{\gamma^{3}}$ are isomorphic and only the projectives to the point $c_{0}$ are not annihilated by $\overline{\gamma^{3}}$. It follows easily 
that 
$A(1)$ and $A(0)$ have the same Auslander-Reiten quiver namely  $(\Z\,D_{3m})_{C}/\langle \sigma \rangle$. The first statement of the theorem is proven.

That $A(1)$ and $A(0)$ are isomorphic for $char\,k \neq 2$  is already shown in the proof above. In $char\,k=2$ we consider the full subcategories $B(1)$ and $B(0)$ supported by $c_{0}$ and $c_{1}$ which are isomorphic to the algebras $kB$ and $kP$ from section 9.5.  An isomorphism  $A(1)\simeq A(0)$  would induce an isomorphism  $kB\simeq kP$ because $c_{0}$ is the only point in $\overline{Q}$ with a loop. Thus $A(1)$ and $A(0)$ are not isomorphic by  lemma 10. Furthermore $A(0)$ is the linearization of the stem category and $A(0)$ the linearization of the ray category because the restrictions to the subcategory supported by $c_{0}$
  and $c_{1}$ are so  and because the definitions of stem and ray categories are compatible with restrictions to full subcategories. 
  
  Finally  a distributive algebra $A$ with quiver $Q$ is regular iff $A(x,y)$ is cyclic from both sides for all $x\neq y$ where an arrow $x \rightarrow y$ in $Q$ exists.
In $D(Q)$ this condition does not hold only for $\alpha_{1}$ and $\alpha_{t}$. However in the two-fold cover $B=D(Q,Q,Q)//\langle \sigma^{2} \rangle$ of $D(Q)$ this condition hods for the 
lifted arrows and so $B$ is regular even though it has a large fundamental group.\hspace{0.1cm}q.e.d.\vspace{0.2cm}

\section{Comments on the classification}
\subsection{The presented proof}

The classification of the simply connected case treated in sections 4,5 and 6 is taken from \cite{BLR} but with several changes. 
We formulate the main result as a bijection between two subsets of $\Z T$:  the set of configurations  and the set of $S$-patterns. This bijection  is obtained only
by elementary Auslander-Reiten theory via  knitting and knotting without using tilting theory.

In section 5 about the tree class $A$ we include a direct proof that pedigrees are the simply connected algebras with an equioriented 
$A_{n}$ as a section in the Auslander-Reiten quiver. We add the information one needs to go on to tree class $D$ by one-point extensions.
This is somewhat hidden in \cite[section 7]{BLR} which is therefore a bit difficult to read. Also we describe only the fundamental algebras and consider their trivial extensions later.
The fundamental groups listed in section 7 are taken from  \cite[appendix 2]{Rie1} and \cite[section 3]{Rie2} and so the classification of all possible Auslander-Reiten quivers of representation-finite sefinjective algebras is complete at the end of section 7.

 Using  the results of the articles \cite{BrG}, \cite{BGRS} and  \cite{B5} the classification of all algebras with a given Auslander-Reiten quiver follows easily. However,
 we proceed independently since reversely some of the concepts introduced  in these  
articles  e.g. multiplicative bases and ray categories appear here for the first time as do the difficulties in characteristic 2.    
                                                                                 
 In section 8 we consider distributive categories and their  stem and ray categories. We include the result of Kupisch and Waschbüsch from \cite{Kupisch5} that any regular
representation-finite selfinjective algebra is isomorphic to the linearization of its stem category. An easy remark about multiplicative bases
in a distributive category asserts that  a cancellative basis exists only in the linearization of  its ray category. 

In section 9 we recall Galois coverings and we define combinatorial
 quotients. They have always a cancellative basis and so the somewhat mysterious ray categories occur here 'in nature'. We prove  that certain covering functors
 preserve the stem and ray categories and that all algebras are regular except possibly some of tree class $D_{3m}$. Apart from these exceptional cases the classification is given  by taking
 the combinatorial quotients.

Section 10 contains the characterization of representation-finite trivial extensions due to Hughes and Waschbüsch and its application to the classification announced in \cite{HW}.

 In section 11
we consider large fundamental groups and we find the non-standard algebras in characteristic 2. Our proof is still a little bit technical but much simpler than the original sources 
\cite{Wasch1,Rie3}. 
 
\subsection{The approach via coverings of the Auslander-Reiten quiver}

This is the approach of Riedtmann with contributions of Gabriel, Bretscher and Läser. In \cite{Rie1} she introduced universal coverings of ( stable ) translation quivers and she proved theorems 1 and 3
of this article as well as theorem 2 for stable categories. Her program to finish the classification was to determine first for  each tree class $T$ the possible combinatorial configurations $C$ with 
their admissible automorphism groups $\Pi$, then to find out which of these triples $(T,C,\Pi)$ are realized by an Auslander-Reiten quiver $\Z\,T_{C}$ and finally to determine  the
 corresponding algebras by quivers and relations. She carried through the first step of this  program for all trees and all three steps for the tree class $A$ in \cite{Rie2}. 
 In 1979 she discovered  
that in characteristic 2 and for the trees  $D_{3m}$ there are non-isomorphic algebras with the same Auslander-Reiten quiver.

 Bretscher found in the beginning of 1981  the bijection between configurations and $S$-patterns. This was then applied 
 together with old results of Riedtmann on combinatorial configurations  in \cite{BLR} to classify the simply onnected algebras and there is also a general homological argument in \cite{BLR} saying that - 
except for tree class $A$ already treated by Riedtmann - only
for  $D_{3m}$ and fundamental group generated by $\tau^{2m-1}$ there can occur non-isomorphic algebras with the same Auslander-Reiten quiver. This problem was finally solved by Riedtmann in the 
remarkable but difficult 
article \cite{Rie3} that was used later on by Asashiba in \cite{Asa} to determine the derived equivalence classes. 
Thus following only this  approach one has in 1983 a published version of the complete classification.

Only in section 5.3 I cannot follow the arguments of \cite{BLR} and I have to add the general argument of \cite[3.1]{BrG}, but 5.3 plays no role in the proof we present.

\subsection{The approach via coverings of quiver and relations}

Based on some necessary conditions for finite representation type obtained by Jans \cite{Jans} Kupisch introduced already 1965 in \cite{Kupisch1} the stem algebra and he used it to obtain informations
 on the Cartan numbers 
of symmetric representation-finite algebras in \cite{Kupisch1,Kupisch2}. Finally he proved in \cite{Kupisch4} that such an  algebra is isomorphic 
to its stem-algebra. Later  he corrected an error in characteristic 2 he was hinted at by Riedtmann and he generalized the result together with Scherzler to arbitrary representation-finite selfinjective algebras
 in \cite{Kupisch3}.
Here only the Jans-conditions are used to find a uniserial projective before one applies an induction. The proofs are quite involved and long and so it is nice that he found 1983 in \cite{Kupisch5}
 together with 
Waschbüsch theorem 9 that we have used. 

Meanwhile Waschbüsch and Scherzler had studied intensively the article \cite{GRie} of Gabriel and Riedtmann where in the special case of Brauer quiver algebras coverings
 of quivers with relations appear for the first time. They  generalized it in the two articles \cite{Wasch3,Wasch4} and Waschbüsch alone 
classified in \cite{Wasch1} the non-regular symmetric algebras and he introduced  in \cite{Wasch2} the universal cover of regular selfinjective algebras. 
Finally he established  together with Hughes in \cite{HW} a bijection between the isomorphism classes of representation-finite trivial extensions and the reflection classes of tilted algebras
 of Dynkin quivers.

Many interesting ideas and results are contained in these articles  e.g. the stem-algebra and the repetetive categories as coverings of trivial extensions.
These have been used and studied later on by Skowronski, Happel and many others.

However there are also several drawbacks compared to the other approach.
First the proof  that the isomorphism classes  of locally 
representation-finite simply connected self-injective categories are in bijection to those of connected representation-finite trivial extension contains an essential gap. 
Namely that proof is only 
sketched in \cite[section 4.2]{HW} and it uses that the universal cover $\tilde{A}$  constructed in section 2.3 of \cite{Wasch2} for any regular  stem algebra $A$
contains no oriented cycle. The only assumption used in the `proof' of this statement in the lemma contained in 2.3 is the regularity and in this generality  the assertion is wrong as the
 quiver $Q$ with two points $a,b$ and arrows $\lambda:a \rightarrow a$, $\rho:b \rightarrow b$, $\mu:a \rightarrow b$ and $\nu:b \rightarrow a$ and 
the relations $\rho \mu = \mu \lambda$, $\nu \rho = \lambda \nu$,$\rho^{2}=\mu \nu$, $\lambda^{3}=\nu\mu$, $\rho^{3}=0=\lambda^{4}$ shows. This gap is closed only some years later by Fischbachers theorem  in \cite{Fisch2} or else one uses our theorem 15 proven by the first approach. 
 
 Second  the actual classification of the trivial extensions is only done for the tree class $A$ in \cite{HW} where it was already known by \cite{Rie2}.
 Type $D$ is considered later in 1989 by 
Assem and Skowronski in \cite{AS} and  
the exceptional
 types $E$  by Roggon in \cite{Roggon} 
via computer.

Third the second approach is not independent of the first because Riedtmanns structure theorem for stable components is used. 
Also the difficulties in characteristic 2 were  detected by Riedtmann. They caused some delay in her work  and two 
corrections in the work of the Berlin group.

Finally the argumentation and the references are several times  too sketchy e.g. on the beginning of page 88 in \cite{Wasch3}.
On the other hand Waschbüsch introduces in that paper and in his  book \cite{Wasch5} the notion of an $(X,J)$- sequence that can be seen as a forerunner to the cleaving diagrams from \cite{BGRS}.

. 

To sum up,  the classifications both depend on Auslander-Reiten theory and coverings. The main difference is that the Zürich group ( Riedtmann, Gabriel, Bretscher, Läser ) starts  with the Auslander-Reiten quiver and its 
universal cover which always exists having nice properties 
whereas the Berlin group ( Kupisch, Waschbüsch, Scherzler, Hughes ) either needs first some information on the algebras in question and uses then  coverings of  quivers with relations where the required properties  were only 
proven some years 
later or it concentrates on  trivial extensions where the repetetive categories exist.

In 1983 the  classification of all representation-finite selfinjective algebras was only obtained 
via the first approach.

\subsection{Some history related to representation-finite selfinjective algebras}

The classification of the indecomposable representations for  a bloc with cyclic defect obtained independently by Janusz and Kupisch on the basis of Dades results is an early 
highlight in the representation theory 
of selfinjective  algebras  ( see \cite{Dade,Janusz,Kupisch} ). 

The next remarkable article is \cite{M} from 1974  where Müller considers  to any path algebra $kQ$ of a quiver with alternating orientation its trivial extension $T(kQ)$.  Just
 by calculating with certain integer valued vectors he gets an easy proof that $kQ$ is representation-finite iff $Q$ is a Dynkin diagram and that then $T(kQ)$ has twice as much indecomposables
 as $kQ$. So here trivial extensions help to understand representation-finite quivers whereas in \cite{HW} the tilted algebras of Dynkin quivers give all representatin-finite trivial extensions.

Maybe more important than the classification is the introduction of coverings to achieve this goal. In \cite{GRie} coverings occur for the first time in representation theory as coverings 
of Brauer quivers with lifted relations. Then coverings of translation quivers were considered as explained in our sections 2 and 3.
For the  results on Galois coverings or coverings of quivers with relations it is difficult to find out who did something first for various reasons.

First there was no internet and no archive where people could communicate with a precise date their
 new results.  Only  preprints were sent to some other people usually without a date.
 
  Second many important articles appeared only in the proceedings of ICRA II 1979 and ICRA III 1980. 
Quite often the 
articles published in the proceedings differ considerably from the talks given at the conference since they are written some months after the conference and they appear
 without any dates in the proceedings. 
For a concrete example we look at coverings. In volume 903 of the Springer lecture notes on ICRA III 
the titles of the talks by Gabriel resp. Green  resp. Waschbüsch were 'Coverings in representation theory', 'Graded algebras and their representations' and
 ' Representations of selfinjective algebras' whereas 
the titles of the articles are 'The universal cover of a representation-finite algebra', 'Group-graded algebras and the zero-relation problem' and 'Universal coverings of selfinjective algebras'.
If one reads the articles one sees that Galois coverings in our sense occur as a common central part in all three papers. In the references Gabriel mentions a preprint by Gordon and Green about
 gradings in representation theory, 
Green the talk by Gabriel and 
Waschbüsch the paper \cite{GRie} already mentioned.

 Gabriel
 and Green discuss Galois coverings without giving a general criterion when the category upstairs has better properties than that downstairs whereas Waschbüsch considers the special 
case of regular stem-algebras claiming that then such a covering always exists which is not true.
Gabriel and Green give only some examples among which are 
zero-relation algebras where the universal cover $\tilde{A}$ is a tree-algebra which is in general infinite. To decide whether this is locally representation-finite or not they both make 
the same
 mistake we will discuss in section   14.4. 

Thus the main problem with coverings is to find a Galois covering $F:\tilde{A} \rightarrow A$ such that $A$ is representation-finite iff $\tilde{A}$ is locally representation-finite and 
such that this is easy to decide.

To go on to mild categories  and the proofs of the Brauer-Thrall conjectures and the finiteness criterion one needs the results on 
multiplicative bases discussed in the next section.

\section{Further developments}

\subsection{Indecomposables}

Locally representation-finite simply connected categories   cannot be classified if they are not selfinjective but their indecomposable representations can.
Namely, the support $S$ of any indecomposable $U$ is a finite convex subcategory  ( \cite{quadratic} ) whence again simply connected ( \cite{Larrion,BrG,B3} ) and therefore given by a graded tree $(T,g)$ as in  \cite{BG}. As soon
 as the number $n(S)$ of points in $S$ is larger than $1000$
it appears in the list 
 LSS of large sincere simply connected algebras given in \cite{B2}. This  list LSS contains 
 24 infinite 
families of algebras depending on several parameters.
 Each of these contains a 'long line' $L$ as  a convex subcategory i.e.  the path algebra of an  $A_{n}$-quiver with $6n \geq n(S)$.
This fact is crucial later on for the proof of BT2 and for the finiteness criterion.
The algebras with less than 1000 points and their sincere  indecomposables are classified in \cite{Dr,RogatTesche} via computer and  the 
very rough estimate 1000 can be replaced by 13 ( see \cite{Tame} ).

An important   consequence  is:

\begin{theorem} ( \cite{B2} ) A representation-finite basic algebra of dimension $d$ has only indecomposables in 
dimensions smaller than the maximum of $2d$ and $30$. 
\end{theorem}
As a corollary   one gets an explicit upper bound 
for the number of indecomposables of a representation-finite algebra $A$ depending  only on $dim\,A$  (\cite{BB}).

\subsection{Multiplicative bases}
The central article \cite{BGRS} is  discussed with many details and some supplements in \cite{GR,B6}. 
The main result about multiplicative bases says:
\begin{theorem}
 Let $A$ be a distributive  $k$-category such that $k\vec{A}$ is mild. Then $A$ has a multiplicative basis. 
If $char\,k \neq 2$ then  $A$ and $k\vec{A}$ are equivalent.

\end{theorem}

We describe the proof of this result in the general case and and also what it amounts to in the special case of selfinjectives. Observe that a mild selfinjective category is
 automatically locally representation-finite whence distributive.

 We  recall some definitions: 
Let $A$ be a distributive $k$-category with quiver $Q$. For an arrow $\alpha:x \rightarrow y $ in $Q$ we denote by $\vec{\alpha}$ the $A(y,y)-A(x,x)$ 
bimodule $\mathcal{R}A(x,y)$ considered as an element in $\vec{A}$. A presentation $\phi:kQ \rightarrow A$ is just a choice of a generator $\phi{\alpha}$ in each $\vec{\alpha}$. Given a path
$v=\alpha_{n}\alpha_{n-1} \ldots \alpha_{1}$ in $Q$ we denote by $\vec{v}$ the product of the $\vec{\alpha_{i}}$ inside $\vec{A}$. The path is called stable if $\vec{v} \neq 0$ or - equivalently -
if $\phi(v) \neq 0$ for all presentations. A contour is a pair $(v,w)$ of stable paths with $\vec{v}=\vec{w}$. Finally  a zig-zag of length $k$ in a ray category $P$ 
 is  a  sequence of morphisms
 $(\rho_{1},\rho_{2},\ldots ,\rho_{k}  )$  such 
that two consecutive morphisms $\rho_{2i-1}$ and $\rho_{2i}$ 
 have common domain,  $\rho_{2i}$ and $\rho_{2i+1}$  common  codomain and none of any $\rho_{j}, \rho_{j+1}$ factors through  the other.
The ray-category $P$ is zigzag-finite if each morphism occurs  only in finitely many zigzags as the first morphism $\rho_{1}$.

To show $A\simeq k\vec{A}$ one starts with an arbitrary presentation  $\phi:kQ \rightarrow A$ and  modifies it in three steps. After the first step  $\phi$
 annihilates all non-stable paths at the same time. After the second one has  $\phi v= c(v,w)\phi w$ for any contour $(v,w)$ with appropriate scalars $c(v,w)$. In the last step these scalars are
 all made to $1$.

Step 1 is taken  in section 4 of  \cite{BGRS} by analyzing the possible non-stable paths. Such a path does not 
 occur in a regular category and so step 1 reduces for selfinjective categories by theorem 13 to our section 11.2.

In step 2 the equality $\phi v= c(v,w) \phi w$ holds automatically for any contour $(v,w)$ where $\phi v$ is in the two-sided socle of the bimodule $A(x,y)$. 
The remaining contours are called non-deep. 
 In the selfinjective case a non-deep contour is only possible for a large fundamental group and tree class $D_{3m}$ because one has $dim\,A(x,y) \geq 2$ for $x\neq y$ resp. $dim\,A(x,x) \geq 3$
and one gets a $D_{4}$-point in the universal cover. Thus we are again in the situation of 11.2 
where $\alpha_{t}\alpha_{1}$ is a non-stable path and $(\alpha_{t}\alpha_{t-1}\ldots \alpha_{1},\gamma^{2})$ the only non-deep contour. It is called a penny farthing.
 The discussion in 9.5 shows that then in characteristic 2 the first two steps cannot be made simultaneously.

 In the general case there are  three different types of non-deep contours called  dumb-bells, penny fartings and diamonds and each of them can occur several times but without overlappings.
 The proof of this in sections 5,6 and 7 of
\cite{BGRS}  it is the most difficult technical part  ( see \cite{B7} for an easier proof ). Somewhat surprisingly and luckily for us  these new contours cause no new obstacles to take the first two steps in 
characteristic different from  2.

The functions $(v,w) \mapsto c(v,w)$ occurring in the third step define an element in the cohomology group $H^{2}(\vec{A},k^{\ast})$ of \cite{BGRS}. An ingenious inductive argument in \cite[section 8.4]{BGRS}
 based on the next theorem shows 
that this group vanishes.

\begin{theorem}
 Let $A$ be a distributive $k$-category such that $\vec{A}$ is not  zig-zag finite. Then $A$ is not locally representation-finite. It is even strongly unbounded if $A$ 
has only finitely many objects.
\end{theorem}

In the selfinjective case step 3 is for regular algebras somewhat hidden in the proof of theorem  9  and for the exceptional case in 11.2 it is the argument using  the Brauer-tree.

The next result of mine in \cite{B5}  is an important supplement to \cite[section 9]{BGRS} where a weaker statement is proven using covering theory and the truth of  BT2 via Gabriels `Finite representation type is open' from \cite{G3}.

\begin{theorem}
 Let $A$ be a distributive connected $k$-category. Then $A$ is mild iff $k\vec{A}$ is mild and then both categories have isomorphic  Auslander-Reiten quivers. If a mild $A$ has a faithful indecomposable even  
 $A \simeq k\vec{A}$ holds.
\end{theorem}

It is remarkable that all results discussed so far in this subsection are obtained  without using covering theory but - as a substitute - the new technique of cleaving diagrams. According to \cite{BGRS} this  technique is due to Bautista, Larri\'{o}n and Salmer\'{o}n but as already mentioned 
a similar construction is considered earlier by  Waschbüsch in \cite{Wasch1}.

\subsection{Coverings of ray categories}

The universal cover $\pi: \tilde{P} \rightarrow P$ of a ray category - defined in  \cite{BGRS} in the obvious way - is only helpful if it has good properties. This is guaranteed by the next important result of Fischbacher  \cite{Fisch2}.

\begin{theorem}
 Let  $\pi:\tilde{P} \rightarrow P $ be the universal cover of a connected ray category which is  zig-zag finite. 

Then the following holds:
\begin{enumerate}
 \item The fundamental group is a free. 
\item $H^{2}(P,Z)=0$ for all abelian groups $Z$.
\item  $\tilde{P}$ is interval-finite and  each connected finite convex subcategory $C$ is again simply connected.
\end{enumerate}
\end{theorem}

One does not know exactly when $\pi:\tilde{P} \rightarrow P$ has the three properties stated above but one knows by \cite{delapena1} that for an arbitrary ray category $P$ the construction of $\tilde{P}$ 
contains as a subproblem the word problem for finitely presented groups.

\subsection{The finiteness criterion}

Let $C$ be a simply connected directed finite ray category. For  any field $k$ the linearization $kC$ has by \cite{B3} a
simply connected component $K$ in its Auslander-Reiten quiver which is obtained by  the well-known knitting. Since the dimension vectors of the indecomposable projectives and of 
the direct summands
 of their radicals depend only on $C$ the component $K$ has the same shape and the same dimension vectors for any field.  We  say that $C$ is representation-finite if all $kC$ are so which 
means that $K$ is finite.

A directed finite ray category is called critical if it is not representation-finite but any proper convex subcategory is. In that case the simply connected component $K$ from above 
contains all indecomposable
 projective $kC$-modules but no injective ( see e.g. \cite{B10} ). Thus $kC$ is given by a graded tree $(T,g)$ with $K \subseteq \Z\,T$ and $T$ is independent of $k$.

\begin{theorem}\cite{B4}
 Using the above assumptions and notations only Euclidean trees $T$ occur for critical simply connected ray categories.
For a tree  $\tilde{D}_{n}$ only the ray categories $\tilde{D}(p,q,r)$ shown in  figure 13 occur.
\end{theorem}

In the next figure  we have 
 $q\geq 1,p \geq 0, r\geq 0$ and for $p=0$ resp. $r=0$ the points $c_{i}$ 
resp. $f_{i}$ are missing. For $p>0$ resp. $r>0$ there is a commutivity relation starting in  $a$ resp. $d$ and these are the only relations. The edges represent arrows in arbitrary orientation
except for the possible commutativity relations.

\setlength{\unitlength}{0.8cm}
\begin{picture}(13,7)\put(6.5,1){figure 13}
\put(3,6){$a$}
\put(2,5){$c_{1}$}
\put(2,4){$c_{2}$}

\put(2,2){$c_{p}$}
\put(3,1){$b$}
\put(4,3.5){$z_{1}$}

\put(5.5,3.5){$z_{2}$}
\put(10,3.5){$z_{q}$}
\put(11,6){$d$}
\put(12,5){$f_{1}$}
\put(12,4){$f_{2}$}
\put(12,2){$f_{r}$}
\put(11,1){$e$}
\put(3.2,5.8){\line(2,-5){0.6cm}}
\put(10.2,3.3){\line(2,-5){0.6cm}}
\put(3.2,1.3){\line(2,5){0.6cm}}
\put(10.2,3.8){\line(2,5){0.6cm}}
\put(4.6,3.6){\line(1,0){0.6cm}}
\put(2.9,5.85){\vector(-1,-1){0.5cm}}
\put(11.9,1.85){\vector(-1,-1){0.5cm}}
\put(2.2,1.8){\vector(1,-1){0.5cm}}
\put(11.2,5.8){\vector(1,-1){0.5cm}}
\put(2.1,4.8){\vector(0,-1){0.3cm}}
\put(12.1,4.8){\vector(0,-1){0.3cm}}

\multiput(6.5,3.6)(0.2,0){14}{\circle*{0.03}}
\multiput(2.1,3.6)(0,-0.2){6}{\circle*{0.03}}
\multiput(12.1,3.6)(0,-0.2){6}{\circle*{0.03}}
\end{picture}

\begin{theorem}\cite{B3,BGRS,B8} Let $A$ be a basic connected distributive $k$-category of dimension $d$ with ray category $P$.
  Then $A$ is representation-finite iff the following three conditions are satisfied:
\begin{enumerate}
 \item $P$ contains no zigzag of length 2d.
\item $\tilde{P}$ contains no $\tilde{D}(p,q,r)$-category with at most $2d+2$ objects as a convex subcategory.

\item All convex subcategories of $\tilde{P}$ with at most $9$ points are representation-finite.
\end{enumerate}
Moreover $A$ is strictly unbounded if it  is not representation-finite.
\end{theorem}

As a consequence the representation-finiteness of a basic distributive algebra of
 dimension $d$ can be decided in 'polynomial time' depending only $d=dim\,A$ ( see \cite{delapena1} ) whereas the calculation of the Auslander-Reiten
quiver needs more efforts just because the number of its points can grow exponentially as in \cite{B3}.

The three conditions on $P$ given in the theorem do not refer to a field. Thus $kP$ is representation-finite for all fields iff  it is so for one field and the Auslander-Reiten quivers are all isomorphic
to the same translation quiver $\Gamma_{P}$ which is always the quotient of $\Gamma_{k\tilde{P}}$ by the fundamental group. Thus one gets the next result which is related to \cite{Brenner}:

\begin{theorem}
 Using these notations the map  $$P \mapsto \Gamma_{P}$$ induces a bijection between the isomorphism classes of finite representation-finite ray categories and those of
finite translation-quivers that are Auslander-Reiten quivers.   .

\end{theorem}

\subsection{The Brauer-Thrall conjectures}
Recall that a $k$-linear functor $F:mod\,B \rightarrow mod\,A$ is a representation embedding if it is exact, maps indecomposables to indecomposables and if $FU\simeq FV$ implies $U\simeq V$.
A  $k$-algebra $A$ is $k$-split if $k$ is  the endomorphism algebra of all simple $A$-modules. Any finite dimensional algebra over an algebraically closed field $k$ is $k$-split by Schurs lemma.

The next theorem is the strongest form of the second Brauer-Thrall conjecture that I know:
\begin{theorem} \cite{B9}
 Let $A$ be a $k$-split basic representation-infinite algebra  of dimension $d$. Then there is an $k[T]-A$ bimodule $_{k[T]}M_{A}$ which is free as a $k[T]$-module
of rank $r\leq max(4d,30)$ such that tensoring with $M$ is a representation embedding $mod\,k[T] \rightarrow mod\,A$.
\end{theorem}

 Thus over any infinite field the Jordan-blocs of size $n$ to  different 
eigenvalues produce  infinitely many pairwise non-isomorphic indecomposable $A$-modules in dimension $rn$.

I call the next theorem  Brauer-Thrall 0 even though it was not conjectured by Brauer and Thrall since it
 is not at all  obvious for group representations.
\begin{theorem}\cite{B6}
 There is no gap in the lengths of the indecomposables  over a finite dimensional  $k$-split algebra. 
\end{theorem}

The proofs of both theorems depend on the next theorem which has so far a quite involved proof:
\begin{theorem}\cite{B6}
 Let $A$ be a distributive finite dimensional mild $k$-algebra which is $k$-split.  Then the universal cover $\pi:\tilde{P} \rightarrow P$ of $P=\vec{A}$ 
has the three properties listed in theorem  24.
\end{theorem}

\subsection{Minimal representation-infinite algebras}
There are some types of minimal representation-infinite algebras where a full classification is known: The tame concealed algebras are  described   by their frames  in the well-known  HV-list of Happel and Vossieck in\cite{Listebild},   the minimal representation-infinite special biserial algebras are classified by Ringel in \cite{Rinmin} and  the non-distributive ones by myself in \cite{B10}. There  I explain also that  the classification of all minimal representation-infinite algebras will lead to a finite but completely unreadable list.

\section{Comments on the further developments}
\subsection{Indecomposables}
 Given any representation-finite algebra $A$ any indecomposable $A$-module $U$ is a faithful 
indecomposable over the quotient $B$ of $A$ by the annihilator of $U$. By theorem  24 the algebra $B$ is isomorphic to $k\vec{B}$. Then there is some indecomposable 
$\tilde{U}$ over $\tilde{B}$ such that $U\simeq F_{\bullet}\tilde{U}$  where $F:\tilde{B} \rightarrow B$ is the universal covering whence a combinatorial quotient ( see theorems 10 and 24 ). Since we `know' all possible $\tilde{U}$ we also `know' all possible indecomposables $U$.

\subsection{Multiplicative bases}
In the spring 1981 Roiter visited Zürich for a month and he wanted to discuss his preprint \cite{Bongo} on multiplicative bases with Gabriel and myself.
His proof consisted in a local part which he called semi-normalization and a global part called normalization.  Roiter knew that the preprint contained gaps and errors. He claimed that 
the global part 
is almost the same as in my trivial proof \cite{B1} for directed algebras, but much to my frustration I did not understand his explanations given to me  only  in  discussions  during 
longer walks along the lake. So I continued to work on the list LSS and the proof of theorem 21 presented in July 1981 at Oberwolfach.

Some years later I looked at section 8 of \cite{BGRS} and on the German translation of the Russian preprint. The lexikographic ordering on zig-zags defined there was never mentioned by Roiter 
in Zürich. This ordering  is used
 in the key lemma 8.4 leading to the proof of the global part and also to Fischbachers result.

The  definition of the lexicographic ordering in turn requires  our theorem 22   which is the unproven lemma 5' in Roiters preprint. As mentioned  in \cite[section 3.3]{BGRS} the important tool of 
cleaving diagrams originated in a 
rigorous proof of
lemma 5'. 

Another important notion not contained in Roiters preprint are the ray categories predecessors of which can be found in \cite{G5,delapena2,BrG}. Also the beautiful structure theorems
 for non-deep contours are not in the preprint but only some weaker statements.

 However the general strategy of Roiters proof  was made to work thanks to the efforts of Bautista, Gabriel and Salmeron.

\subsection{Coverings of ray categories}

Section 10 of \cite{BGRS} has no counterpart in the preprint but it is very important. Here coverings of ray categories and their topological properties are studied  and a version of the finiteness criterion is formulated.
Most of this is simplified and generalized by Fischbacher in \cite{Fisch2} to theorem 24 of the  present survey  using only the key lemma 8.4 of \cite{BGRS} mentioned above.

\subsection{The finiteness criterion }

In the literature there are some  crooked statements  about the finiteness criterion whose version given in 13.4  seems to be widely unknown.  So 
I  explain   here  the proof in several steps always saying which results are needed. After that I make some critical remarks.

 Thus we consider  a basic distributive algebra $A$ of dimension $d$ with ray category $P$ and we prove theorem 26.

\underline{Step 1}: We can reduce to $A=kP$.

\underline{Proof of step 1}: Indeed theorem 23 says that $A$ is representation-finite iff $kP$ is so and even that $A\simeq kP$ if $A$ is minimal representation-infinite.

\underline{Step 2}: $A$ is strictly unbounded if condition 1 is not satisfied.
 
\underline{Proof of step 2}: Let $(\rho_{1},\rho_{2}, \ldots ,\rho_{2d})$ be a zig-zag in $P$. Since the $\rho_{i}$ are not identities there exist indices $i <j <k$ with $\rho=\rho_{i}=\rho_{j}=\rho_{k}$ and also $p<q$ such that $\rho=\rho_{p}=\rho_{q}$ point into the same direction inside of the zig-zag.
 Then we get the infinite periodic zig-zag
$(\rho_{p},\rho_{p+1},\ldots \rho_{q}=\rho_{p},\rho_{p+1},\ldots \rho_{q}=\rho_{p}, \ldots )$. Thus $P$ is not zig-zag finite and step 2 follows from theorem 22. 

\underline{Step 3}: The conditions are necessary.

\underline{Proof of step 3}: Condition 1 holds by step 2. If $kP$ is representation-finite then $k\tilde{P}$ is locally representation-finite by  theorem 10 whence conditions 2 and 3 are satisfied.

\underline{Step 4}: If condition 1 holds  we have:
\begin{enumerate}
\item There is no line of length $2d$ as a convex subcategory of $\tilde{P}$.
\item $k\tilde{P}$ is locally representation-finite iff all of its finite convex subcategories are representation-finite.
\item BT 2 holds for $kP$.
\end{enumerate}

\underline{Proof of step 4}: As seen above condition 1 implies that $P$ is zig-zag finite and so the universal cover $\pi: \tilde{P} \rightarrow P$ has all the nice properties listed in theorem 25.
Thus one can apply lemma 7.2 of \cite{B8} to get a periodic zig-zag from a line of length $2d$.

In the second statement the implication from the left to the right is trivial. To prove the other direction we lift a spanning tree $T$ of the quiver of $P$ to a tree $\tilde{T}$ in  the quiver $\tilde{Q}$ of $\tilde{P}$. Because at most three arrows start or end in any point of $\tilde{Q}$  the set of  points   that can be  reached from a point in $\tilde{T}$ by a walk of length $\leq m$  is finite for any $m$ and so is its convex hull  $C(m)$ by the  interval-finiteness of $\tilde{P}$.  Furthermore $C(m)$ is  representation-finite. Here  the points  in a walk should be pairwise different and their number is the length of the walk e.g. a point is a  walk of length $1$.

 Now we  take a point $x\in \tilde{T}$ and  an indecomposable representation $U$ of $\tilde{P}$ with $U(x)\neq 0$. We claim that the support $S$ of $U$ is contained in  $C(12d)$.  If not there is a point  $y \in S\setminus C(12d)$. Since $S$ is connected $y$  can be joined to $x$  by a walk $w$  within $S$ say of length $l$. Then $l>12d \geq 12$   and $S$  has at least  $l$ elements.  $S$ is a representation-finite simply connected convex subcategory with $U$ as a sincere indecomposable and therefore $S$ occurs in the list LSS. By the properties mentioned in 13.1 there is a line of length $2d$  in $S$ and so in $\tilde{P}$ contradicting the first assertion in step 4.
  We have shown that $U$ is in fact a representation of the representation-finite algebra $kC(12d)$. Since any point of $\tilde{P}$ is conjugate to a point $x\in \tilde{T}$ any indecomposable over $\tilde{P}$ is conjugate to an indecomposable over $kC(12d)$ and so $k\tilde{P}$ is locally representation-finite.

Finally  if $kP$ is not representation-finite  $k\tilde{P}$ is not locally representation-finite by theorem 10. Then there is a finite convex subcategory $C$ in $k\tilde{P}$ which is not representation-finite. We can assume that $C$ is critical whence a one-point extension of a product of representation-finite simply connected algebras. By \cite[section 6]{BG} the second Brauer-Thrall conjecture holds for $C$ and then also for $kP$ again by theorem  10.

\underline{Step 5}: The conditions are sufficient.

\underline{Proof of step 5}: By step 4 all finite convex subcategories are representation-finite iff there is no critical subcategory $C$. By theorem 25  such a category is given by a graded tree $(T,g)$ with a Euclidean tree $T$ and this is excluded by the conditions 2 and 3. Observe here that any algebra $\tilde{D}(p,q,r)$ with  more than $2d+2$ points contains a line of length $2d$ which is excluded by the first assertion in step 4.\vspace{0.5cm}

Note that the proof of  BT2 in steps 2 and 4  uses only the list LSS and not any characterization of critical simply connected algebras or tame concealed algebras let alone the HV-list.

Here are some more historical and mathematical remarks.

\begin{itemize}
\item 
In september 1982 at a conference in Luminy I presented the first version of the criterion. 
Theorems 21,22,23,24 and ray categories with their universal coverings did not yet exist. So I considered a kind of 'standard' algebra $A$ of dimension $d$ with a simply connected Galois-covering
$\pi:\tilde{A} \rightarrow A$ that generalized the representation-finite situation. The main points of my  theorem are  that such a 'standard' algebra $A$ is representation-finite iff four conditions are satisfied  each of which can be checked in finitely many steps 
 and that it is strictly unbounded otherwise. The last two  of these four conditions are exactly  the
 last two in theorem 26
whereas the first two follow now from theorem 22 and the zig-zag finiteness.

At the same conference in Luminy Happel gave a talk about tame algebras including the  classification  of the tame concealed algebras i.e. the Happel-Vossieck list of frames. Happel did not mention any characterization of the tame concealed algebras.

 Since a critical simply connected algebra is tame concealed  ( e.g. by \cite{B10} ) the last condition of the theorem  says that no algebra of the HV-list occurs as a convex subcategory of $\tilde{P}$. This is extremely helpful  if one calculates  many examples  by hand  as in \cite{Fisch3}.

\item The next version of the criterion appears in section 10 of \cite{BGRS}. For several reasons I do not like this version which is formulated for a not necessarily finite ray category $P$. 

First, the main points of my version are  a finite algorithm to decide whether a finite $P$  is  locally representation-finite and  a proof of  BT2 if it is not. For an infinite $P$ there is no finite algorithm and BT2 does not hold.

Second, the HV-list is displayed and explained in full detail whereas the list LSS is not even mentioned. 
Furthermore the formulation  of the finiteness-criterion in \cite{BGRS} is reproduced in \cite{G6,Rie5,GR} and BT2 is then proven with the help of the criterion giving the wrong impression that the HV-list is essential for the proof of BT2.

\item The final version of the criterion in theorem 27 is based on  contributions of many mathematicians. So it should not be called 'Bongartz-criterion'. However already the first version needs more than the`first Jans-condition and preprojective tilting' as stated by Ringel in an article about canonical algebras.

\item At many places in the literature e.g in \cite{G5,Green3,Br} the subtle point in step 4 is overseen: The equivalence of the two assertions  '$\tilde{P}$ is locally representation-finite'  and 'all finite subcategories of $\tilde{P}$ are representation-finite' is  only true if $P$ contains no line of length $2d$.
It is here where the list LSS is used and where the proof of BT 2 is not obvious. 

\item If the fundamental group $\Pi$ of the universal cover $\pi:\tilde{P} \rightarrow P$ is $\Z$ the method of truncating $\tilde{P}$ and knitting `until one gets periodicity '  can be applied to decide the representation type and to compute the Auslander-Reiten quiver if it is finite  ( see  \cite{G5} resp. \cite{BrG}  ). Observe that here one has to know when to stop the knitting because $P$ is not representation-finite. This `halteproblem' is solved for any fundamental group by the following result  depending only on theorem 21 and   the definition of DTr and TrD.
\begin{theorem} 
 Using  the notations introduced in the proof above  $kP$ is representation-finite iff $kC(N+2d)$ is representation-finite with indecomposables of dimension $N=max(2d,30)$ at most.
Moreover in that case any indecomposable in $mod\,kP$ is the push-down of an indecomposable in $mod \,kC(N+2d)$ and any almost split sequence in $mod \,kP$ is the push-down of an almost 
split sequence in $mod\,kC(N+2d)$.
\end{theorem}

\end{itemize}

\subsection{My contributions}

I am a very limited mathematician working on rather special questions. Between 1984 and 2009 I wrote some  articles about geometric properties of module varieties 
 but the basic problems discussed in this article is what really interests me. 

To conclude  I say which of my results mentioned in this survey I like most.
\begin{itemize}

\item My favourite is theorem 29  saying that there are no gaps in the lengths of the indecomposables. The statement is elementary but my proof uses almost all of the theory of mild algebras, but no lists. Furthermore it answers in 2009 a question I stumbled on already in my diploma thesis in 1973.
\item The most important result is  the list LSS of large sincere simply connected representation-finite algebras and its various applications: An explicit bound on the dimension of indecomposables over a representation-finite algebra of dimension $d$ ( theorem 21 ),   
the solution of the halteproblem
in the computation of finite Auslander-Reiten quivers ( theorem 31 ), the finiteness-criterion  and all known proofs of BT2.
\item From a  lyrical point of view my best result is in the common article \cite{BR} with Ringel: For each  indecomposable over a representation-finite tree algebra there is a vertex $p$ such that each arrow pointing to $p$ is represented by a mono, each arrow pointing away by an epi.

In German this can be formulated shortly as: \begin{center}{  Unzerlegbare über einer darstellungsendlichen Baumalgebra haben  Wipfel.}
\end{center}
\end{itemize}
\begin{center}{\itshape Über allen Gipfeln, \\ist Ruh\\In allen Wipfeln\\spürest Du
\\kaum einen Hauch.\\Die Vögelein schweigen im Walde.\\Warte nur: Balde\\ruhest Du auch. }
 
\end{center} 
\hspace*{6cm} {\itshape  Johann Wolfgang von Goethe}

\end{document}